\documentclass[preprint,12pt]{elsarticle}

\usepackage{amssymb}
\usepackage{amsmath}
\usepackage{amsthm}
\usepackage{amssymb}
\usepackage{amsfonts}
\usepackage{hyperref}
\usepackage{epsfig,anysize}
\usepackage{enumerate}
\usepackage{newlfont}
\usepackage{latexsym}

\renewcommand{\proof}{\noindent{\bf Proof. \quad}}
 \newcommand{\eproof}{\hfill \mbox{${\square}$} \vspace{3mm}}

\newcommand{\R}{\mathbb{R}}

\newtheorem {teo} {Theorem\,} [section]
\newtheorem {prop} [teo] {Proposition}
\newtheorem {cor} [teo]{Corollary}
\newtheorem {lema} [teo] {Lemma}

\newtheorem{rem}[teo]{Remark}

\renewcommand{\proof}{{\noindent \bf Proof. \ }}

\begin{document}

\begin{frontmatter}



\title{Continuity of attractors for a highly oscillatory family of  perturbations of the square}

\author[BL1]{Bianca Paolini Lorenzi\fnref{label1}}
\ead{bianca.lorenzi@alumni.usp.br}
\fntext[label1]{Partially supported by CAPES-Brazil grant 88882.377938/2019-01.}
\affiliation[BL1]{organization={Universidade de São Paulo},
             city={São Paulo},
             country={Brazil}}

\author[AP1]{Antônio L. Pereira \corref{cor2}\fnref{label2}}
\ead{alpereir@ime.usp.br}
\cortext[cor2]{Corresponding author.}
\fntext[label2]{Partially supported by FAPESP-SP Brazil grant 2020/14075-6.}
\affiliation[AP1]{organization={Instituto de Matemática e Estatística, Universidade de São Paulo},
   addressline={Rua do Matão, 1010},
             city={São Paulo},
             postcode={05508-090},
             country={Brazil}}

\begin{abstract}
\indent Consider the family of semilinear parabolic problems
\begin{equation*}
\left\{
\begin{array}{lll}
u_{t}(x,t) = \Delta u(x,t) - au(x,t) + f(u(x,t)), \,\,\, x \in \Omega_{\epsilon}, t > 0, \\
\frac{\partial u}{\partial N} (x,t) = g(u(x,t)), \,\,\, x \in \partial \Omega_{\epsilon}, t > 0,
\end{array}
\right.
\end{equation*}

\noindent where $a > 0$, $\Omega$ is the unit square, $\Omega_{\epsilon} = h_{\epsilon}(\Omega)$, $h_{\epsilon}$ is a family of $C^{m}$ - diffeomorphisms, $m \geq 1$, which  converge to the identity of $\Omega$ in $C^{\alpha}$ norm, if $\alpha <1$  but do not converge in the $C^{1}$ - norm and, $f,g: \mathbb{R} \rightarrow \mathbb{R}$ are real functions.
 We show that  a weak version of this problem, transported to the fixed domain $\Omega$ by a ``pull-back''  procedure,  is well posed for $0 <\epsilon \leq \epsilon_{0}$, $\epsilon_{0} > 0$, in a suitable phase space, the associated semigroup has a global attractor $\mathcal{A}_{\epsilon}$ and the family $\{ \mathcal{A}_{\epsilon} \}_{0 \, < \, \epsilon \, \leq \, \epsilon_{0}}$ converges as  $\epsilon \to  0$
to the attractor of the limiting problem:

\begin{equation*}\
\left\{
\begin{array}{lll}
u_{t}(x,t) = \Delta u(x,t) - au(x,t) + f(u(x,t)), \,\,\, x \in \Omega, t > 0, \\
\frac{\partial u}{\partial N} (x,t)  = g(u(x,t))\mu, \,\,\, x \in \partial \Omega, t > 0,
\end{array}
\right.
\end{equation*}
\noindent where $\mu$ is essentially the limit of the Jacobian determinant of the diffeomorphism ${h_{\epsilon}}_{| \partial \Omega} : \partial \Omega \rightarrow \partial h_{\epsilon}(\Omega)$ (but does not depend on the particular family $h_{\epsilon})$.
\end{abstract}


 
\begin{keyword}
parabolic problem \sep perturbation of the domain \sep Lipschitz domain \sep global attractor \sep continuity of attractors

MSC: Primary: 35B41; Secondary: 35K91, 58D25
\end{keyword}

\end{frontmatter}

\section{Introduction} \label{sec_intro}

\indent   Consider the family of semilinear parabolic problems
\begin{equation} \label{pert_bianca}
\left\{
\begin{array}{lll}
u_{t}(x,t) = \Delta_{\Omega_{\epsilon}} u(x,t) - au(x,t) + f(u(x,t)), \,\,\, x \in \Omega_{\epsilon}, t > 0, \\
\frac{\partial u}{\partial N_{\Omega_{\epsilon}}} (x,t) = g(u(x,t)), \,\,\, x \in \partial \Omega_{\epsilon}, t > 0,
\end{array}
\right.
\end{equation}

\noindent where $a > 0$, $\Omega$ is the unit square, $\Omega_{\epsilon} = h_{\epsilon}(\Omega)$, $h_{\epsilon}$ is a family of $C^{m}$ - diffeomorphisms, $m \geq 1$,   and $f,g: \mathbb{R} \rightarrow \mathbb{R}$ are real functions.

The same problem was considered in \cite{pricila1}, with the family of diffeomorphisms given by
 \begin{equation} \label{pert_pricila}
 h _{\epsilon}(x_{1}, x_{2}) = (x_{1}, x_{2} + x_{2} \epsilon \sin(x_{1}/ \epsilon^{\alpha}))
 \end{equation}
\noindent with $0 < \alpha < 1$ and $\epsilon$
small. Under these hypotheses, the authors prove in \cite{pricila1}   that the problem is well posed for $0 \leq \epsilon \leq \epsilon_{0}$, $\epsilon_{0} > 0$, in a suitable phase space, the associated semigroup has a global attractor $\mathcal{A}_{\epsilon}$ and the family $\{ \mathcal{A}_{\epsilon} \}_{0 \,\leq \, \epsilon \, \leq \, \epsilon_{0}}$  is continuous at $\epsilon = 0$.
   In order to obtain this result, the authors used, in a essential way,  the convergence of the  family of diffeomorphisms \eqref{pert_pricila}  to the identity in $\Omega $ in the $C^1$-norm.  Our aim here is to investigate what happens in the `critical case' $\alpha =1$,   when  the diffeomorphisms oscillate more wildly, so  this hypothesis does not hold.

One of the first difficulties encountered in problems of perturbation of the domain is that the spaces change with the change of the region. One way to overcome this difficulty is to perform a change of variables in order to bring the problem back to a fixed region and   use the tools  developed by D. Henry in \cite{henry}.
More precisely, if $ \Omega_{\epsilon}= h_{\epsilon}(\Omega) $,  we consider the family of problems in the fixed region $\Omega$,
\begin{equation} \label{pert_bianca_fixed}
\left\{
\begin{array}{lll}
u_{t}(x,t) = h_{\epsilon}^{*}\Delta_{\Omega_{\epsilon}} {h_{\epsilon}^{*}}^{-1} u(x,t) - au(x,t) + f(u(x,t)), \,\,\, x \in \Omega_{h_{\epsilon}}, t > 0, \\
h_{\epsilon}^{*}\frac{\partial }{\partial N_{\Omega_{\epsilon}}} {h_{\epsilon}^{*}}^{-1}u(x,t) = g(u(x,t)), \,\,\, x \in \partial \Omega_{\epsilon}, t > 0,
\end{array}
\right.
\end{equation}
where $ h_{\epsilon}^{*} u = u \circ h_{\epsilon}$ is the `pullback' operator.
It is easy to show that, if $h_{\epsilon}$ is smooth enough,
then $v$ is a solution of \eqref{pert_bianca} if and only if
$ u = v \circ h $ is a solution of \eqref{pert_bianca_fixed}.
 The solutions of
\eqref{pert_bianca_fixed} lie then in the same function space for any $\epsilon$  and, therefore, can be directly compared.
This is the approach taken, for instance, in \cite{pricila1}, \cite{pricila2} and  \cite{marcone}. In \cite{marcone}   the authors considered the problem
\eqref{pert_bianca} in the case of \emph{regular perturbations} of the domain, that is, when
 $\Omega$ is a $C^{2}$ bounded region in   $\mathbb{R}^{n}$,  $\Omega_{h_{\epsilon}} = h_{\epsilon}(\Omega)$, and $h_{\epsilon}$ is a family of smooth diffeomorphisms that converge to the identity in $\Omega$ in the  $C^{2}$-norm.
  In \cite{pricila2} these   results were extended to the case where the convergence of the diffeomorphisms is assumed only  in the $C^1$-norm.

  An essential step in these works is the proof that the  family of operators
 $ - h_{\epsilon}^{*}\Delta {h_{\epsilon}^{*}}^{-1} $ in a weak form, given by
 \begin{equation} \label{def_A_fixed_weak_intro}
 \langle A_{\epsilon} {u},  \psi \rangle_{-1,1} =
 \
\int_{\Omega} h_{\epsilon}^{*}\nabla_{\Omega_{\epsilon}}h_{\epsilon}^{*-1} (u)(x,t) \cdot h_{\epsilon}^{*}\nabla_{\Omega_{\epsilon}}h_{\epsilon}^{*-1} \left(\frac{\psi}{|Jh_{\epsilon}|} \right)(x)
|Jh_{\epsilon}(x)|\, dx
 \end{equation}
where $\langle  \cdot, \cdot \rangle_{-1,1} $ is the dual pairing between $H^1(\Omega) $ and its dual (which we denote here by $H^{-1}(\Omega)) $
  converge, in some sense, to the  (negative of) Neumann Laplacian operator.

   However, in the present case, due to the fact that    our diffeomorphisms do not converge to the identity  in a sufficiently strong norm (for instance, in the  $C^1$-norm), this does not hold anymore. In fact, these operators converge to different ones, depending on the particular family of diffeomorphisms,

 To obtain the `correct' weak form, we first multiply the whole equation \eqref{pert_bianca}  by a test function $\psi$, obtaining the integral form

\begin{equation}\label{PI_intro}
\left\{
\begin{array}{lll}
\int_{\Omega_{\epsilon}} v_{t}(y,t)\psi(y)dy = \int_{\Omega_{\epsilon}} \Delta_{\Omega_{\epsilon}} v(y,t)\psi(y) dy  - a \int_{\Omega_{\epsilon}} v(y,t) \psi(y) dy \\
+ \int_{\Omega_{\epsilon}} f(v(y,t)) \psi(y) dy , \,\,\, y \in \Omega_{\epsilon}, t > 0, \\
\int_{\partial \Omega_{\epsilon}} \frac{\partial v}{\partial N_{\Omega_{\epsilon}}} (y,t) \psi(y) d\sigma(y) = \int_{\partial \Omega_{\epsilon}} g(v(y,t)) \psi(y) d\sigma(y), \,\,\, y \in \partial \Omega_{\epsilon}, t > 0.
\end{array}
\right.
\end{equation}\\
Now, integrating by parts in the right side and then  applying the pull-back operator to the whole equation, we obtain the weak form in the fixed domain (see section \ref{sec_weak_form2}. for details):

\begin{eqnarray} \label{weak_correct_intro}
\int_{\Omega} u_{t}(x,t) \phi(x) |Jh_{\epsilon}(x)|dx & = & -\int_{\Omega} (h_{\epsilon}^{*} \nabla_{\Omega_{\epsilon}} h_{\epsilon}^{*-1} u) (x,t) \cdot (h_{\epsilon}^{*} \nabla_{\Omega_{\epsilon}} h_{\epsilon}^{*-1} \phi) (x) |Jh_{\epsilon} (x)| \, dx \nonumber  \\
& - & a \int_{\Omega} u(x,t) \phi(x) |Jh_{\epsilon}(x)| \, dx  \nonumber \\
& + & \int_{\Omega} f(u(x,t)) \phi(x) |Jh_{\epsilon}(x)| \, dx \nonumber \\
& + & \int_{\partial \Omega} g(u(x,t)) \phi(x) |J_{\partial \Omega_{\epsilon}} h_{\epsilon} (x)| \, d\sigma(x)
\end{eqnarray}
and we now have to deal with the family of operators
\begin{equation} \label{def_A_fixed_weak_intro_correct}
 \langle A_{\epsilon} {u},  \psi \rangle_{-1,1} =
 \
\int_{\Omega} h_{\epsilon}^{*}\nabla_{\Omega_{\epsilon}}h_{\epsilon}^{*-1} (u)(x,t) \cdot h_{\epsilon}^{*}\nabla_{\Omega_{\epsilon}}h_{\epsilon}^{*-1} \psi (x)
|Jh_{\epsilon}(x)| \, dx
 \end{equation}
which we will show to converge in the sense of resolvents  to the (weak form)  of the minus  Neumann Laplacian.

      We will show then  that  the family of attractors of \eqref{weak_correct_intro}
        converge   as $\epsilon \to 0$
 to the attractor of a `limiting problem'
\begin{equation*}\
\left\{
\begin{array}{lll}
u_{t}(x,t) = \Delta u(x,t) - au(x,t) + f(u(x,t)), \,\,\, x \in \Omega, t > 0, \\
\frac{\partial u}{\partial N} (x,t)  = g(u(x,t))\mu, \,\,\, x \in \partial \Omega, t > 0,
\end{array}
\right.
\end{equation*}
\noindent where $\mu$ is essentially the limit of the jacobian determinant of the diffeomorphism ${h_{\epsilon}}_{| \partial \Omega} : \partial \Omega \to \partial h_{\epsilon}(\Omega)$.

A different technique, based on extension operators, is used in \cite{simone_rapidly} and \cite{simone_very_rapidly} where the authors consider the elliptic equation with nonlinear boundary conditions
\begin{equation} \label{eqarrieta}
\begin{array}{rlr}
\left\{
\begin{array}{lll}
-\Delta u + u = f(x,u) \,\,\,\text{em} \,\,\,\Omega_\epsilon \,,\\
\displaystyle\frac{\partial u}{\partial n} + g(x,u) = 0 \,\,\,\text{em}\,\,\,\partial\Omega_\epsilon\,,
\end{array}
\right.
\end{array}
\end{equation}
where $\Omega_{\epsilon} \subset \mathbb{R}^{n}, \, 0 \leq \epsilon \leq \epsilon_{0}$, is a family of smooth domains with  $\Omega_{\epsilon} \to \Omega$ and  $\partial \Omega_{\epsilon} \rightarrow \partial \Omega$ in the sense of  Hausdorff and $f,g$ are real functions continuous in the first and of class  $C^{2}$ in the second variable in a bounded domain containing  $\Omega_{\epsilon}$,  for all $0 \leq \epsilon \leq \epsilon_{0}$.

When the boundary of the domain oscillates  rapidly as the parameter $\epsilon$ converges to $0$,  the authors obtain the convergence of solutions to the ones of the limiting problem:
 \begin{equation} \label{eqlimitearrieta}
\begin{array}{rlr}
\left\{
\begin{array}{lll}
-\Delta u + u = f(x,u) \,\,\,\text{em} \,\,\,\Omega,\\
\displaystyle\frac{\partial u}{\partial n} + \gamma(x)g(x,u) = 0 \,\,\,\text{em}\,\,\,\partial\Omega,
\end{array}
\right.
\end{array}
\end{equation}
\noindent where $\gamma \in L^{\infty}(\partial \Omega)$  and  $\gamma \geq 1$.

A similar problem was considered in \cite{simone_very_rapidly}, with  faster oscillations, so that the Lipschitz constant of the  local  charts  of their boundaries  are not uniformly bounded.
 In this case,  they show that the limit problem has Dirichlet type boundary conditions.

  We observe that one additional difficulty here is that, since  our  region is only Lipschtiz regular, we have not been able to use some standard  results established in this case. In particular, we could not use  comparison results to obtain uniform bounds for the solutions, as done, for instance,  in \cite{marcone}.
 Using the results in \cite{grisvard},  most of our results could be easily extended to more general convex domains, but we have chosen to restrict  to an specific setting, for the sake of clarity.

 This paper is organized as follows: in Chapter 2, we show that the family of operators defined by \eqref{def_A_fixed_weak_intro} does not converge to the minus Laplacian. In Chapter 3, we obtain a weak form that is shown to define a sectorial operator in Chapter 4 and to converge in the resolvent sense to the minus Laplacian in Chapter 5.
 Chapter 6 is dedicated to the formulation of the problem in an abstract form in a convenient scale of Banach spaces, for which the well posedness is then shown in Chapters 7 and 8. Finally, the existence and continuity of global attractors to the one of our limit problem is proved in the remaining Chapters.

 \section{The weak form  in the canonical metric} \label{section_weak_form1}

 \subsection{Reduction to a fixed domain}
\indent \par Let  $\Omega = ]0,1[ \times ]0,1[ \subset \mathbb{R}^{2}$ be the unit square in $\R^2$,  $a$ a positive constant and  $f,g : \mathbb{R} \rightarrow \mathbb{R}$ real functions. Consider the  problem
 \eqref{pert_bianca} written in a slightly different notation,
 for convenience:
 \begin{equation}\label{pert_bianca2}
\left\{
\begin{array}{lll}
v_{t}(y,t) = \Delta_{\Omega_{\epsilon}} v(y,t) - av(y,t) + f(v(y,t)), \,\,\, y \in \Omega_{\epsilon}, t > 0, \\
\frac{\partial v}{\partial N_{\Omega_{\epsilon}}} (y,t) = g(v(y,t)), \,\,\, y \in \partial \Omega_{\epsilon}, t > 0,
\end{array}
\right.
\end{equation}
where $\Omega_{\epsilon} = h_{\epsilon} (\Omega)$ and  $h_{\epsilon} : \Omega \rightarrow \mathbb{R}^{2}$  is  a family of  diffeomorphisms
converging to the identity in  the $\mathcal{C}^{\alpha} $ norm  as
 $\epsilon \to 0$ (with $\alpha$ to be specified later).

   As observed in the introduction, we can   bring the problem back to the fixed region $\Omega$, by  means of the pullback operator
    $h_{\epsilon}^{*} : C^m(\Omega_\epsilon) \to C^m(\Omega)$, $  v \to  u = h_{\epsilon}^{*}v$, giving raise to the problem \eqref{pert_bianca_fixed}.
\begin{equation} \label{pert_bianca_fixed2}
\left\{
\begin{array}{lll}
u_{t}(x,t) = {h_{\epsilon}}^{*}\Delta_{\Omega_{\epsilon}} {{h_{\epsilon}}^{*}}^{-1} u(x,t) - au(x,t) + f(u(x,t)), \,\,\, x \in \Omega, t > 0, \\
{h_{\epsilon}}^{*}\frac{\partial }{\partial N_{\Omega_{\epsilon}}} {{h_{\epsilon}}^{*}}^{-1}u(x,t) = g(u(x,t)), \,\,\, x \in \partial \Omega, t > 0,
\end{array}
\right.
\end{equation}

 To obtain a weak form of \eqref{pert_bianca_fixed2}, we multiply by a test function $\psi: \bar{ \Omega } \to \R$ and integrate,
 \begin{eqnarray}
 \int_{\Omega} u_{t}(x,t) \psi(x) \, d\,x   & = &
 \int_{\Omega} h_{\epsilon}^{*}\Delta_{\Omega_{\epsilon}} {{h_{\epsilon}}^{*}}^{-1} u(x,t)\psi(x) \, d\,x  -
 \int_{\Omega} au(x,t) \psi(x) \, d\,x  \nonumber \\
  & + &
 \int_{\Omega}  f(u(x,t)) \psi(x) \, d\,x
\end{eqnarray}

 Integrating by parts, we get
\begin{align} \label{weak_form_operator}
 \displaystyle\int_\Omega ({h_{\epsilon}}^*\Delta_{\Omega_{\epsilon}}{h_{\epsilon}}^{*-1}u)(x,t) &
 \,\psi(x)\,dx
 =  \displaystyle\int_\Omega \Delta_{\Omega_{\epsilon}} (u \circ {h_{\epsilon}}^{-1})({h_{\epsilon}}(x),t) \psi(x)\,dx  \nonumber \\
&=  \displaystyle\int_{\Omega_{\epsilon}}
\Delta_{\Omega_{\epsilon}} v(y,t) \psi({h_{\epsilon}}^{-1}(y)) \frac{1}{|\,J{h_{\epsilon}}({h_{\epsilon}}^{-1}(y))\,|}\,dy  \nonumber \\
&=  \int_{\partial \Omega_{\epsilon}} \frac{\partial v}{\partial N_{\Omega_{\epsilon}}}(y,t) \psi({h_{\epsilon}}^{-1}(y)) \frac{1}{|J{h_{\epsilon}}({h_{\epsilon}}^{-1}(y))|}\,d\sigma(y)  \nonumber \\
& -  \int_{\Omega_{\epsilon}} \nabla_{\Omega_{\epsilon}} v(y,t)\cdot \nabla_{\Omega_{\epsilon}} {\left(\psi({h_{\epsilon}}^{-1}(y))  \frac{1}{|J{h_{\epsilon}}({h_{\epsilon}}^{-1}(y))|}\right)} dy
\nonumber \\
&= -\int_{\Omega_{\epsilon}} \nabla_{\Omega_{\epsilon}} v(y,t)\cdot \nabla_{\Omega_{\epsilon}} {\left(\psi(h_\epsilon^{-1}(y))  \frac{1}{|Jh_\epsilon(h_\epsilon^{-1}(y))|}\right)}dy  \nonumber  \\
 &  +  \int_{\partial \Omega_{\epsilon}} g(v(y,t)) \psi(h_\epsilon^{-1}(y)) \frac{1}{|Jh_\epsilon(h_\epsilon^{-1}(y))|}\,d\sigma(y)
 \nonumber\\
&=  -
\displaystyle\int_{\Omega} (h_\epsilon^*\nabla_{\Omega_{\epsilon}}h_\epsilon^{*-1}u)(x,t) \cdot\left(h_\epsilon^* \nabla_{\Omega_{\epsilon}} h_\epsilon^{*-1} \displaystyle\frac{\psi}{|Jh_\epsilon|}\right)(x)|Jh_\epsilon(x)| \, dx  \nonumber \\
 & +
  \int_{\partial \Omega} g(u(x,t)) \psi (x)  \left(\frac{|J_{\partial \Omega}h_{\epsilon} (x)|}{|Jh_{\epsilon}|} \right)  d\sigma(x) .
\end{align}

\indent Therefore we obtain the weak form  of the problem
  \eqref{pert_bianca_fixed2}   in the original domain $\Omega$,
\begin{eqnarray} \label{pert_bianca_fixed2_weak}
\int_{\Omega} u_{t}(x,t) \psi(x)\, d\,x & = & - \int_{\Omega} h_{\epsilon}^{*}\nabla_{\Omega_{\epsilon}}h_{\epsilon}^{*-1} (u)(x,t) \cdot h_{\epsilon}^{*}\nabla_{\Omega_{\epsilon}}h_{\epsilon}^{*-1} \left(\frac{\psi }{|Jh_{\epsilon}|} \right) (x)
|Jh_{\epsilon}(x)|\, dx \nonumber \\
& - & \int_{\Omega} au(x,t) \psi(x)\, d\,x + \int_{\Omega} f(u(x,t)) \psi(x)\, d\,x \nonumber \\
& + & \int_{\partial \Omega} g(u(x,t)) \psi (x)  \left(\frac{|J_{\partial \Omega}h_{\epsilon} (x)|}{|Jh_{\epsilon}|} \right)  d\sigma(x).
\end{eqnarray}

  It was shown in \cite{pricila1} and \cite{pricila2} that, if $ h_{\epsilon} \to Id_{\Omega}$ in the $C^1$-norm and satisfies some mild additional conditions, the principal part of the equation above, given by the opposite of the  operator
  defined in \eqref{def_A_fixed_weak_intro}, with domain $H^1(\Omega)$ converges to the  weak Neumann Laplacian in the reference region as
   $\epsilon \to 0$.     However, as we show below this is not true anymore if the convergence holds only in the $C^{\alpha}$ norm, with
    $\alpha < 1$ and  in fact, its limit depends on the particular family of perturbations
 $ h_{\epsilon} : \Omega \to \Omega_{\epsilon}  $
 and not only  on its image $\Omega_{\epsilon} $.

 \subsection{ The limit of the pull-backed  operator}

  We take $a=0$ in this section to simplify the notation.

 Let's  consider the family
 \begin{equation} \label{defhepsilon2}
h_{\epsilon} (x_{1},x_{2}) = (x_{1}, x_{2} + x_{2} \epsilon \sin (x_{1}/\epsilon)),
\end{equation}
 which is the family considered in \cite{pricila1}, except that now, we are dealing with the critical exponent $\alpha =1$.

  The inverse transpose of the Jacobian matrix of   $h_{\epsilon}$ is
\begin{equation*}
\begin{pmatrix}
1 & \dfrac{-x_{2} cos(x_{1}/ \epsilon)}{1 + \epsilon sin(x_{1}/ \epsilon)} \\
0 & \dfrac{1}{1 + \epsilon sin(x_{1}/ \epsilon)}
\end{pmatrix}.
\end{equation*}

Using the expression (4) in \cite{pricila1}, the perturbed operator in the weak form is then given by
\begin{eqnarray}
\langle A_{\epsilon} u, \, \psi \rangle_{-1,1} & = & \int_{\Omega} \frac{\partial u}{\partial x_{1}} (x) \frac{\partial \psi}{\partial x_{1}} (x) dx \label{div1}   \\
& + & \int_{\Omega} \frac{-x_{2} cos(x_{1}/ \epsilon)}{1 + \epsilon sin(x_{1}/ \epsilon)}\frac{\partial u}{\partial x_{2}} (x) \frac{\partial \psi}{\partial x_{1}} (x) dx \label{div2}\\
& + & \int_{\Omega} \frac{-x_{2} cos(x_{1}/ \epsilon)}{1 + \epsilon sin(x_{1}/ \epsilon)} \frac{\partial u}{\partial x_{1}} (x) \frac{\partial \psi}{\partial x_{2}} (x) dx \label{div3}\\
& + & \int_{\Omega} \frac{x_{2}^{2} cos^{2}(x_{1}/ \epsilon)}{(1 + \epsilon sin(x_{1}/ \epsilon))^{2}} \frac{\partial u}{\partial x_{2}} (x) \frac{\partial \psi}{\partial x_{2}} (x) dx \label{div4} \\
& + & \int_{\Omega} \frac{1}{(1 + \epsilon sin(x_{1}/ \epsilon))^{2}}\frac{\partial u}{\partial x_{2}} (x) \frac{\partial \psi}{\partial x_{2}} (x) dx. \label{div5} \\
 & + & \int_{\Omega} - \frac{cos(x_{1}/ \epsilon)}{1 + \epsilon sin(x_{1}/ \epsilon)}\frac{\partial u}{\partial x_{1}} (x) \psi (x)dx \label{l1}\\
& + & \int_{\Omega} \frac{x_{2} cos^{2}(x_{1}/ \epsilon)}{(1 + \epsilon sin(x_{1}/ \epsilon))^{2}}\frac{\partial u}{\partial x_{2}} (x) \psi (x)dx.\label{l2}
\end{eqnarray}

\noindent Let  $C_{\epsilon}$
be the operator given by \eqref{div1} + \eqref{div5}  and
 $\Delta_{\Omega}$, the Laplacian operator in  $H^{-1}(\Omega)$, with domain  $H^{1}(\Omega)$.
 We have

\begin{eqnarray*}
\langle (C_{\epsilon} - (-\Delta_{\Omega}))u, \psi \rangle_{-1,1} & = & \int_{\Omega} \frac{\partial u}{\partial x_{1}} (x) \frac{\partial \psi}{\partial x_{1}} (x) dx +\int_{\Omega} \frac{1}{(1 + \epsilon sin(x_{1}/ \epsilon))^{2}}\frac{\partial u}{\partial x_{2}} (x) \frac{\partial \psi}{\partial x_{2}} (x) dx \\
& - & \int_{\Omega} -\Delta_{\Omega} u(x) \psi(x)dx \\
& = & \int_{\Omega} \left( \frac{1}{(1 + \epsilon sin(x_{1}/ \epsilon))^{2}} - 1 \right) \frac{\partial u}{\partial x_{2}} (x) \frac{\partial \psi}{\partial x_{2}} (x) dx
\end{eqnarray*}
 and  $$\int_{\Omega} \left| \left( \frac{1}{(1 + \epsilon sin(x_{1}/ \epsilon))^{2}} - 1 \right) \right| \left| \frac{\partial u}{\partial x_{2}} (x) \frac{\partial \psi}{\partial x_{2}} (x) \right| dx \leq C(\epsilon) \left| \left| \frac{\partial u}{\partial x_{2}} \right| \right|_{L^{2}(\Omega)} \, \left| \left| \frac{\partial \psi}{\partial x_{2}} \right| \right|_{L^{2}(\Omega)},$$ where  $C(\epsilon)$ is a constant which converges to $0$  as  $\epsilon \to 0 $.

Now, the integrals  \eqref{div2}, \eqref{div3} and \eqref{l1}
go to zero, since $cos(x_{1}/\epsilon)$ converges weakly* to zero in  $L^{\infty}(\mathbb{R})$ and the remainder terms in each integrand are bounded by an integral function  $\eta(x_1,x_2)$,

\indent Now, using that
$ cos^{2}(x_{1}/ \epsilon) =
\frac{1}{2}
\left(1 + cos(2x_{1}/ \epsilon) \right) \to \frac{1}{2} $
 weakly* in  $L^{\infty}(\mathbb{R})$, we obtain that
$
\eqref{div4} \to  \frac{1}{2} \int_{\Omega} {x_{2}^{2}} \frac{\partial u}{\partial x_{2}} (x) \frac{\partial \psi}{\partial x_{2}} (x) dx $.

Similarly, we obtain that
$ \eqref{l2} \to    \frac{1}{2} \int_{\Omega} x_{2} \frac{\partial u}{\partial x_{2}} (x) \psi (x) dx$.

\indent We, therefore obtain
 \begin{equation}\langle A_{\epsilon} u, \psi \rangle_{-1,1} \to \langle (-\Delta_{\Omega} + aI)u, \psi \rangle_{-1,1} + \frac{1}{2} \int_{\Omega} x_{2}^{2} \frac{\partial u}{\partial x_{2}} (x) \frac{\partial \psi}{\partial x_{2}} (x) dx +  \frac{1}{2} \int_{\Omega} x_{2} \frac{\partial u}{\partial x_{2}} (x) \psi (x) dx, \label{lim_operator_1}
\end{equation}
\noindent and the additional terms in the right-hand side are not null.

\vspace{3mm}

We now consider another family of diffeomorphisms given by
 \begin{equation}
 h_{\epsilon} (x_{1}, x_{2}) = (x_{1}, x_{2} + x_{2}\epsilon^{2 - x_{2}} sin(x_{1}/ \epsilon))
 \end{equation}
which
have the same image as
 \eqref{defhepsilon2}.

 Using similar arguments, we now obtain that
 \begin{align*}
 \langle A_{\epsilon} u, \psi \rangle_{-1,1} & = \langle (-\Delta_{\Omega} + a)u , \psi \rangle  \nonumber  \\
 &  +
\int_{\Omega} -x_{2} \epsilon^{1 - x_{2}} ln(\epsilon) cos (x_{1}/ \epsilon) \frac{\partial u}{\partial x_{1}}(x) \psi (x) dx  \nonumber \\
 & +
\int_{\Omega} x_{2}^{2} \epsilon^{2 - 2x_{2}} ln(\epsilon) cos^{2}(x_{1}/ \epsilon)\frac{\partial u}{\partial x_{2}}(x) \psi (x) dx \nonumber \\
 & + o(\epsilon).
\end{align*}
 and the right-hand side has also an additional term, which is different than the one in \eqref{lim_operator_1}

\section{The weak form in an adapted metric} \label{sec_weak_form2}

As we have seen in the previous section, the pull-backed operator  $h_{\epsilon}^{*} \Delta_{\Omega_{\epsilon}}h_{\epsilon}^{-1*}$ in the weak form, associated with the canonical inner product, does not converge to the weak Laplacian in the reference region and, in fact, may converge to different operators, depending on the  specific  family used to perturb the original region. As we will see,  this can be fixed by  using an inner product adapted to the problem.
This new inner product  appears naturally, if we pull-back the whole equation, instead of the operator.

More precisely, if $\psi$ is a smooth function in $\bar{\Omega}_{\epsilon}$, we can rewrite (\ref{pert_bianca2}) in the integral form,

\begin{equation}\label{PI}
\left\{
\begin{array}{lll}
\int_{\Omega_{\epsilon}} v_{t}(y,t)\psi(y)dy = \int_{\Omega_{\epsilon}} \Delta_{\Omega_{\epsilon}} v(y,t)\psi(y) dy  - a \int_{\Omega_{\epsilon}} v(y,t) \psi(y) dy \\
+ \int_{\Omega_{\epsilon}} f(v(y,t)) \psi(y) dy , \,\,\, y \in \Omega_{\epsilon}, t > 0, \\
\int_{\partial \Omega_{\epsilon}} \frac{\partial v}{\partial N_{\Omega_{\epsilon}}} (y,t) \psi(y) d\sigma(y) = \int_{\partial \Omega_{\epsilon}} g(v(y,t)) \psi(y) d\sigma(y), \,\,\, y \in \partial \Omega_{\epsilon}, t > 0.
\end{array}
\right.
\end{equation}\\
Now, integrating by parts in the right hand side,
\begin{eqnarray*}
\int_{\Omega_{\epsilon}} v_{t}(y,t)\psi(y)dy & = & \int_{\Omega_{\epsilon}}  \Delta_{\Omega_{\epsilon}} v(y,t)\psi(y) dy  - a \int_{\Omega_{\epsilon}} v(y,t) \psi(y) dy  \\
& + & \int_{\Omega_{\epsilon}} f(v(y,t)) \psi(y) dy \\
& = & - \int_{\Omega_{\epsilon}} \nabla_{\Omega_{\epsilon}} v(y,t)\cdot \nabla_{\Omega_{\epsilon}}  \psi(y) dy  - a \int_{\Omega_{\epsilon}} v(y,t) \psi(y) dy  \\
& + & \int_{\Omega_{\epsilon}} f(v(y,t)) \psi(y) dy + \int_{\partial \Omega_{\epsilon}} g(v(y,t)) \psi(y) d\sigma(y). \\
\end{eqnarray*}

Then, if $u(x, t) = v(h_{\epsilon}(x),t), \phi (x) = \psi(h_{\epsilon}(x)), \, x \in \Omega$, we obtain
\begin{eqnarray} \label{pert_bianca_fixed2_weak_new}
\int_{\Omega} u_{t}(x,t) \phi(x) |Jh_{\epsilon}(x)| \, dx & = & -\int_{\Omega} (h_{\epsilon}^{*} \nabla_{\Omega_{\epsilon}} h_{\epsilon}^{*-1} u) (x,t) \cdot (h_{\epsilon}^{*} \nabla_{\Omega_{\epsilon}} h_{\epsilon}^{*-1} \phi) (x)  |Jh_{\epsilon} (x)| \, dx  \nonumber \\
& - & a \int_{\Omega} u(x,t) \phi(x) |Jh_{\epsilon}(x)| \, dx  \nonumber \\
& + & \int_{\Omega} f(u(x,t)) \phi(x) |Jh_{\epsilon}(x)| \, dx \nonumber \\
& + & \int_{\partial \Omega} g(u(x,t)) \phi(x) | J_{\partial \Omega_{\epsilon}} h_{\epsilon} (x)| \, d\sigma(x)
\end{eqnarray}
\begin{rem}
A simple computation shows that the metric $\mu = \mu(\epsilon) =
 |Jh_{\epsilon}| \, d\,x $ is exactly the metric $h^{*} \, d \, x $ obtained by pull-backing the canonical metric in $ \Omega_{\epsilon}$. Also, the first term in the operator
$ A_{\epsilon} : H^{1}(\Omega, \mu) \rightarrow H^{-1}(\Omega, \mu)$ defined by
\begin{eqnarray}\label{def_A_fixed_weak_new}
 <A_{\epsilon}u, \phi>_{-1,1, \mu} = \int_{\Omega} (h_{\epsilon}^{*} \nabla_{\Omega_{\epsilon}} h_{\epsilon}^{*-1} u) (x,t) \cdot (h_{\epsilon}^{*} \nabla_{\Omega_{\epsilon}} h_{\epsilon}^{*-1} \phi) (x) | Jh_{\epsilon} (x)| \, dx  \nonumber \\
 + a \int_{\Omega} u(x,t) \phi(x) |Jh_{\epsilon}(x)| \, dx, \, u, \phi \in H^{1}(\Omega, \mu)
\end{eqnarray}
where $\langle  \cdot, \cdot \rangle_{-1,1, \mu} $ is the dual pairing between $H^1(\Omega, \mu) $ and its dual (which we denote here by $H^{-1}(\Omega, \mu )) $,
is the Laplace-Beltrami operator in $\Omega$, with respect to the metric $\mu(\epsilon)$.
\end{rem}

The equation \eqref{pert_bianca_fixed2_weak_new} can be written as an abstract equation in $H^{-1}(\Omega, \mu)$:

\begin{equation*}
{u}_{t} = - A_{\epsilon} u +  H_{\epsilon}{u}.
\end{equation*}

  where

\[H ({u})   = F(u) + G(u), \]

 \begin{align*}
  F(u)(\psi)   & =      \int_{\Omega} f(u(x,t)) \psi(x)\, |Jh_{\epsilon}(x)| \,  d\,x \\
  G(u) (\psi) & =\int_{\partial \Omega} g(u(x,t)) \psi (x)    |J_{\partial \Omega} h_{\epsilon} (x)| \,  d\sigma(x)
  \end{align*}

\noindent with $\psi  \in H^{1} (\Omega, \mu).$

We finally observe that the metrics $\mu(\epsilon)$ generated by the new inner product are uniformly equivalent to the usual euclidean metric, so the spaces   $H^{1}(\Omega, \mu)$ can be identified as metric spaces with the (fixed) usual Sobolev  space   $H^{1}(\Omega)$. This identification will be used henceforth, without further comments in our estimates.

\section{Sectoriality of the linear operators} \label{sec_sectorial}
From now on, we suppose that
 $h_{\epsilon} : \Omega \rightarrow \mathbb{R}^{2}$  is  a family of  diffeomorphisms  given by
\begin{equation} \label{defhepsilon}
h_{\epsilon} (x_{1},x_{2}) = (x_{1}, x_{2} + \varphi (x_{2}, \epsilon) \sin (x_{1}/\epsilon)), \,\, 0 \leq \epsilon \leq \epsilon_{0}, \,\, \epsilon_{0} > 0,
\end{equation}
 where $\varphi (x_{2}, \epsilon)$ is  of class $C^2$ and satisfies the following hypotheses:

\begin{align}\label{hypotheses_varphi}
  \textrm{i)} \ & \varphi (0, \epsilon) = 0, \nonumber \\
 \textrm{ii)} \ &   \varphi (1, \epsilon) = \epsilon,  \nonumber \\
  \textrm{iii)} \ &  \frac{\partial \varphi}{\partial x_{2}} \to 0 \textrm{ uniformly as }  \epsilon \to 0, \nonumber \\
  \textrm{iv)} \ & \left| \left| \varphi (x_{2}, \epsilon)^{2} \frac{1}{\epsilon^{2}} \right| \right|_{L^{2}(\Omega)} \to 0  \textrm{  when } \epsilon \to 0.
\end{align}

\noindent (For instance, the function $\varphi (x_{2}, \epsilon) = x_{2} \epsilon^{2 - x_{2}}$ satisfies these requirements).

\begin{figure}[!h] \label{perturb_square}
\begin{center}
	\includegraphics[scale=4.0]{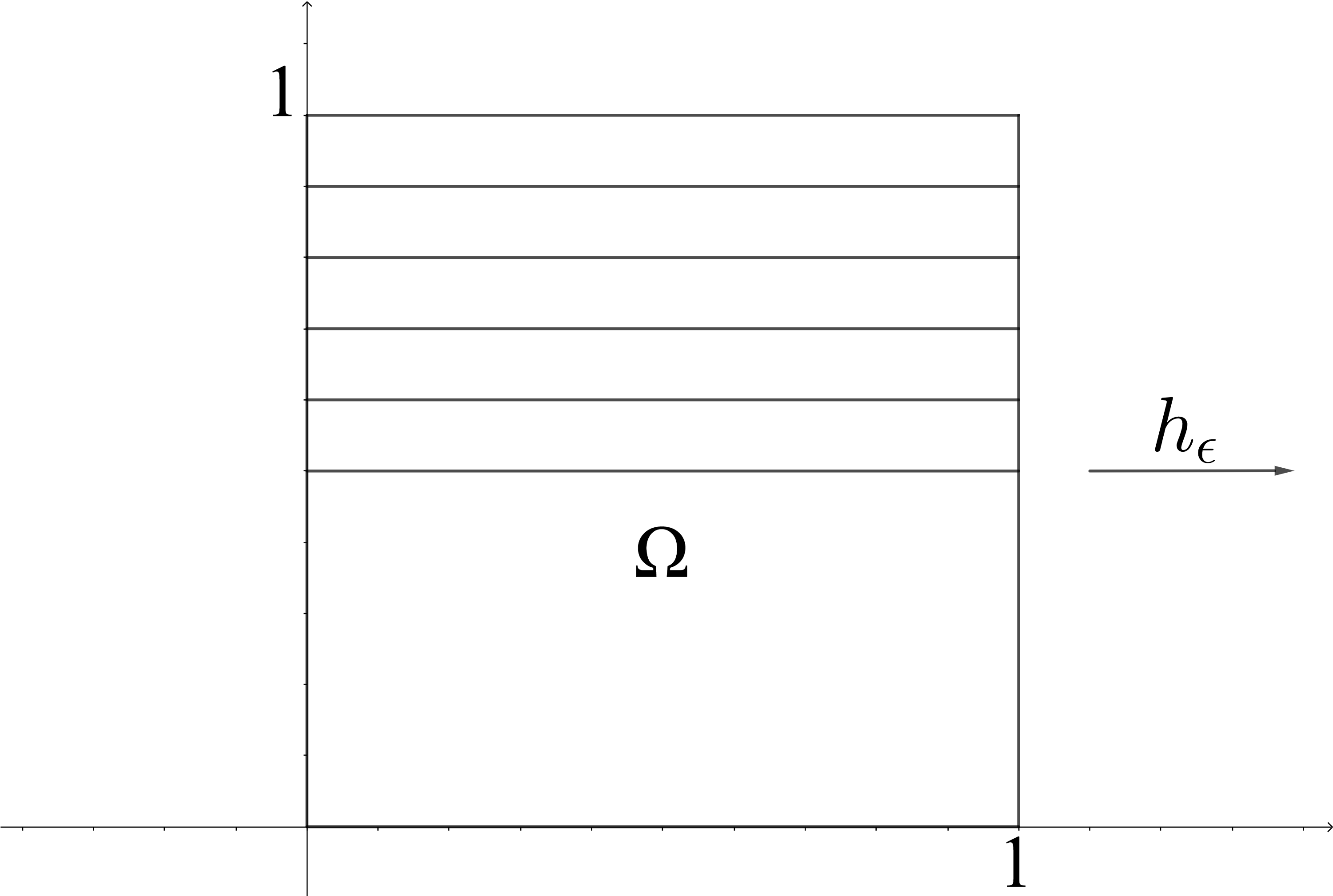} \quad
	\includegraphics[scale=4.0]{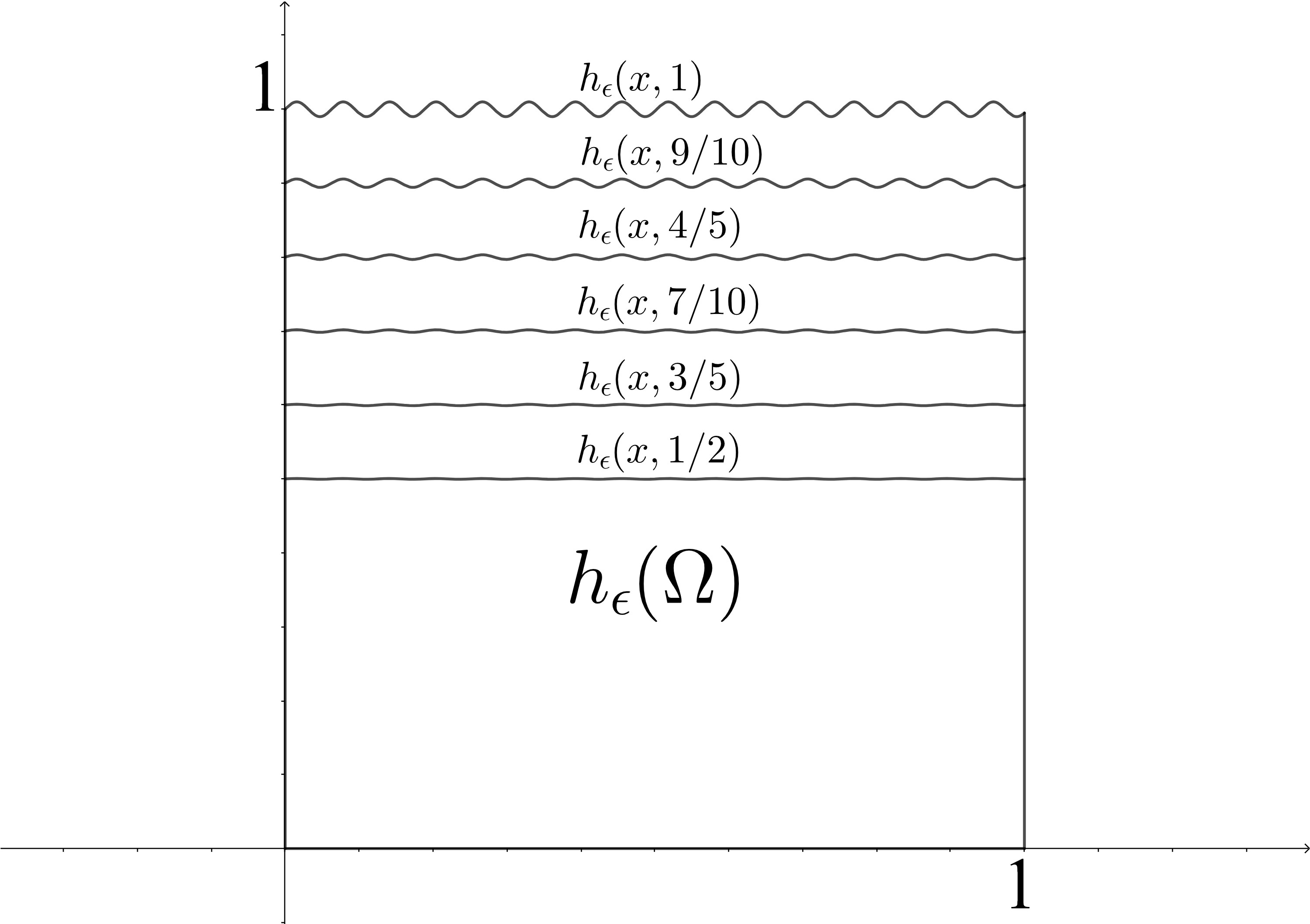}
\caption{The family $h_{\epsilon}$.}
\end{center}
\end{figure}

\indent We observe that the family  $h_{\epsilon}$ in \eqref{defhepsilon}  converges to  $id_{\Omega}$
 as $\epsilon \to 0$ in the
  $C^{\alpha}(\Omega)$ if $\alpha < 1$, but not in the
  $C^{1}(\Omega)$ norm since
  $ h_{\epsilon}(x_{1}, x_{2}) - id_{\Omega} (x_{1}, x_{2})  = (0, \varphi (x_{2}, \epsilon) \sin (x_{1}/\epsilon)) $ is uniformly Lipschitz bounded but this Lipschitz constant does not converge to $0$.

 Consider, for $0< \epsilon \leq \epsilon_0$, the operator
   $-A_{\epsilon}: H^1(\Omega) \subset H^{-1}(\Omega)  \to H^{-1}(\Omega, \mu) $  given in \eqref{def_A_fixed_weak_new}
   with the family $h_{\epsilon}$ satisfying hypotheses (\ref{hypotheses_varphi}).

  We will show that
 $-A_{\epsilon} $ as well as its restriction to  $L^2(\Omega)$ are sectorial operators.

\subsection{Sectoriality in $H^{-1}$}

\begin{lema}\label{self_positive}  The operators $-A_{\epsilon}$
 defined by \eqref{def_A_fixed_weak_new} are  self-adjoint and
  $\left\langle -A_\epsilon u, u \right\rangle_{-1,1,\mu}   \geq C \parallel u \parallel_{H^{1}(\Omega)}^2$ where  $C$ is a positive constant independent of $\epsilon$.
\end{lema}

\proof

\indent Let  $\epsilon > 0$,  $a_{\epsilon} : H^{1}(\Omega) \times H^{1}(\Omega) \rightarrow \mathbb{R}$ be  the bilinear form given by

\begin{equation*}
a_{\epsilon} (u,v) = \int_{\Omega} (h_{\epsilon}^{*} \nabla_{\Omega_{\epsilon}} h_{\epsilon}^{*-1} u) (x,t) \cdot (h_{\epsilon}^{*} \nabla_{\Omega_{\epsilon}} h_{\epsilon}^{*-1} v) (x) |Jh_{\epsilon}(x)| \, dx + a \int_{\Omega} u(x,t) v(x) |Jh_{\epsilon}(x)| \, dx.
\end{equation*}

The bilinear  form  $a_{\epsilon}$ is symmetric,   continuous and
\begin{equation*}
< -A_{\epsilon}u, v>_{-1,1} = a_{\epsilon}(u,v), \,\, \forall v \in H^{1}(\Omega).
\end{equation*}
 Thus $A_{\epsilon}$ is also symmetric.
  We claim that $a_{\epsilon}$ is uniformly coercive, that is, there exists
a constant $C > 0$ such that  $a_{\epsilon} (u,u) \geq C \parallel u \parallel^{2}_{H^{1}(\Omega)}$ for all $u \in H^{1}(\Omega)$.

 In fact, we have
 \begin{align*}
   a_{\epsilon} (u,u) & = \, \int_{\Omega} \sum_{i,j,k = 1}^{2} b_{ij}b_{ik} \frac{\partial u}{\partial x_{j}} \frac{\partial u}{\partial x_{k}} |Jh_{\epsilon}| \, dx + \, a \int_{\Omega}  u^2  |Jh_{\epsilon}|  \,dx \nonumber \\
    & \geq \int_{\Omega} \lambda(\epsilon)
 \sum_{i} {\frac{\partial u}{\partial x_i}}^2 \, |Jh_{\epsilon}|  \, dx
   + a   u^2 \,  |Jh_{\epsilon}| \, d\,x ,
  \end{align*}
  where  the  $b_{ij} = b_{ij}^{\epsilon}$  are the entries of the inverse transposed of the Jacobian matrix of  $h_{\epsilon}$,

\noindent that is,
\begin{equation} \label{bij}
b_{11}^{\epsilon} = 1, \, b_{12}^{\epsilon} = \frac{-\varphi(x_{2}, \epsilon) cos(x_{1}/ \epsilon) \frac{1}{\epsilon}}{1 + \frac{\partial \varphi}{ \partial x_{2}}(x_{2}, \epsilon) sin(x_{1}/ \epsilon)}, \, b_{21}^{\epsilon} = 0, \, b_{22}^{\epsilon} = \frac{1}{1 + \frac{\partial \varphi}{ \partial x_{2}}(x_{2}, \epsilon) sin(x_{1}/ \epsilon)},
\end{equation}
and $\lambda(\epsilon)$ is the smallest of the eigenvalues of the matrix
\begin{equation*}
\begin{pmatrix}
1 & b_{12} \\
b_{12} &  b_{12}^{2} + b_{22}^{2}
\end{pmatrix}.
\end{equation*}

 It is easy to show, using the hypotheses in (\ref{hypotheses_varphi}), that
 $ \lambda (\epsilon) \geq \lambda >0$ and then,
  if $ C = \min \{a,  \lambda \} \cdot \inf\{Jh_{\epsilon}(x) \, : \, x \in \Omega\} $, we obtain

 \begin{align} \label{coercive}
   a_{\epsilon} (u,u) &  \geq C    \|u\|_{H^1(\Omega)}^2,
  \end{align}
as claimed.

 It follows from  Lax-Milgram Theorem, that  $A_{\epsilon}$ is surjective and, therefore, self-adjoint.  The claimed estimate for $A_{\epsilon}$ follows immediately from \eqref{coercive}. \eproof

\begin{cor}\label{sectorialh1} The operators  $\{A_{\epsilon}\}_{0 \, \leq \, \epsilon \,\leq \, \epsilon_{0}}$ in $H^{-1}(\Omega)$  are sectorial and there exist $b>0$ and $M(=\sqrt{2})$ independent of $\epsilon$, such that
$\parallel (A_{\epsilon} - \lambda)^{-1} \parallel \leq \frac{M}{|\lambda - b|}, \, \forall \lambda \in S_{b, \frac{\pi}{4}} =
\{ \lambda \in \mathbb{C} : \frac{\pi}{4} \leq |\arg(\lambda-b)| \leq \pi    \}$. The associated semigroups   $\{T_{\epsilon}(t)= e^{-A_{\epsilon} t} \}_{0 \, \leq \, \epsilon \, \leq \, \epsilon_{0}}$ satisfy $\parallel T_{\epsilon}(t) \parallel \leq Me^{-b t}$, \  $\parallel  A_{\epsilon} e^{-A_{\epsilon}t} \parallel \, \leq \frac{M}{t}e^{-bt}$ and  $\parallel  A_{\epsilon}^{\alpha} e^{-A_{\epsilon}t} \parallel \, \leq M_ {\alpha} \frac{1}{t^\alpha }e^{-bt}$, with $M_ {\alpha}$ bounded in a compact interval of $[0, \infty[ $.
\end{cor}
\proof From Lemma \ref{self_positive}, $A_{\epsilon}$ is self-adjoint and there is a positive constant   $b$ independent  of  $\epsilon$, for  $0 \leq \, \epsilon \,\leq \epsilon_{0}$, such that
$ \|A_{\epsilon} u\| \geq b \|u\|_{H^1 (\Omega)}$ and is densely defined.
 A simple computation then shows that
 $A_{\epsilon}$ is sectorial,  $S_{b, \frac{\pi}{4}} \subset \rho(A_{\epsilon} )$ and  $\parallel (A_{\epsilon} - \lambda)^{-1} \parallel \leq \frac{M}{|\lambda - b|}, \, \forall \lambda \in S_{b, \frac{\pi}{4}}$, with $M= \sqrt{2}$ (see \cite{bianca}).\\
\indent From Theorem 1.3.4 in \cite{henry_semilinear}, $-A_{\epsilon}$ is the infinitesimal generator of an analytic semigroup   $\{ e^{-A_{\epsilon}t}\}_{t \, \geq  \, 0}$ with $\parallel e^{-A_{\epsilon}t} \parallel \, \leq Me^{-bt}$, $\parallel  A_{\epsilon} e^{-A_{\epsilon}t} \parallel \, \leq Mt^{-1}e^{-bt}$ where  $M$ is independent of $\epsilon$. The last inequality is proved in \cite{henry_semilinear}, Theorem
1.4.3.
\eproof
\begin{rem} \label{A_epsilon_equiv_norms}
From Lemma \ref{self_positive}, it follows that
$\|A_{\epsilon} u \|_{H^{-1}(\Omega)} \geq C \|u\|_{H^1(\Omega)} $, where $C$ is a constant independent of $\epsilon$. On the other hand,
the operators  $A_{\epsilon}$ are also bounded as operators from  $H^{1}(\Omega)$ to  $H^{-1}(\Omega)$, uniformly in $\epsilon$
since
\begin{align*}
  \langle  & A_{\epsilon}u, \Phi\rangle_{-1,1,\mu}  \leq \left|  \int_{\Omega}  \sum_{i,j,k = 1}^{2}b_{ij}b_{ik} \frac{\partial u}{\partial x_j}
   \frac{\partial \Phi}{\partial x_k} \, |Jh_{\epsilon}| \, d\,x \right| +
  \left| \int_{\Omega}  a   u \Phi \, |Jh_{\epsilon}| \, d\,x \right|  \nonumber \\
    & \leq B \sum_{jk} \left(\int_{\Omega}
  \frac{\partial u}{\partial x_j}^2 d\,x \right)^{1/2}
   \left(\int_{\Omega}
  \frac{\partial \Phi}{\partial x_k}^2 d\,x \right)^{1/2}
   + A  \left(\int_{\Omega}     u^2  \, d\,x \right)^{1/2}
   \left(\int_{\Omega}     \Phi^2  \, d\,x \right)^{1/2} \\
   & \leq   B \left(\int_{\Omega}   \sum_{j}
  \frac{\partial u}{\partial x_j}^2 d\,x \right)^{1/2}
   \left(\int_{\Omega}   \sum_{k}
  \frac{\partial \Phi}{\partial x_k}^2 d\,x \right)^{1/2}
   + A  \left(\int_{\Omega}     u^2  \, d\,x \right)^{1/2}
   \left(\int_{\Omega}     \Phi^2  \, d\,x \right)^{1/2} \\
   & \leq C'  \|u\|_{H^1(\Omega)}.
  \end{align*}
  where  the  $b_{ij}$  are the entries of the inverse transposed of the Jacobian matrix of  $h_{\epsilon}$, \newline
  $ B= sup \{ |Jh_{\epsilon}(x)| : x \in \Omega \} \cdot \max \{ b_{ij} b_{ik} \}$, and we have used Cauchy-Schwartz inequality.

  Therefore, the norms   $ \| u \|_{1, \epsilon} =
   \|A_{\epsilon} u\|_{H^{-1}(\Omega)} $
   are  equivalent to
   the norm   $\|u \|_{H^{1}(\Omega)} $ with constants of equivalence  independent of  $\epsilon$.

\end{rem}

\begin{cor}\label{limitresolunif} If  $\lambda $ is in the sector
$
S_{b, \frac{\pi}{4}} =
\{ \lambda \in \mathbb{C} : \frac{\pi}{4} \leq |\arg(\lambda-b)| \leq \pi    \} $  and $f \in H^{-1}(\Omega)$, then
  $\parallel (A_{\epsilon} - \lambda)^{-1} f \parallel_{H^{1}(\Omega)} \leq  C  (\frac{|\lambda|}{|\lambda-b|} )\parallel f \parallel_{H^{-1}(\Omega)}$, with $C$ independent of $\epsilon$.
\end{cor}

\proof
Denoting by $ \| \cdot  \|$ the norm in $H^{-1}(\Omega)$ we have
 from Corollary  \ref{sectorialh1},

$$\parallel (A_{\epsilon} - \lambda)^{-1} \parallel \leq \frac{M}{|\lambda - b|}, \, \forall \lambda \in S_{b, \frac{\pi}{4}}.$$ Thus,

 \begin{align}
 \parallel A_{\epsilon} (A_{\epsilon} - \lambda)^{-1} f \parallel & \leq
  \| (A_{\epsilon} - \lambda ) (A_{\epsilon} - \lambda)^{-1} f    \|
  + \|  \lambda  (A_{\epsilon} - \lambda)^{-1} f    \| \nonumber \\
  & \leq \left(1 + \frac{|\lambda | M}{ |\lambda-b|} \right)\|f\|.
  \end{align}
 Using the equivalence of norms of Remark  \ref{A_epsilon_equiv_norms}, the result follows.  \eproof

\subsection{Sectoriality in $L^{2}$}

We now show that the operators $\{A_{\epsilon}\}_{0 \, \leq \, \epsilon \,\leq \, \epsilon_{0}}$  are sectorial as operators in $L^2(\Omega)$. More precisely,

\begin{prop} \label{setoriall2}  The operators
 $\{A_{\epsilon}\}_{0 \, \leq \, \epsilon \,\leq \, \epsilon_{0}}$
  in $L^2(\Omega)$  with domain  \newline
 $\ D(A_{\epsilon}) = \left\{ u \in H^{2} (\Omega) \, | \, h_{\epsilon}^{*} \frac{\partial}{\partial N_{\Omega_{\epsilon}}}h_{\epsilon}^{*-1} u = 0 \, \text{on} \, \partial \Omega \right\}$   are sectorial.
\end{prop}
\proof   The operators are restrictions of the ones in \eqref{def_A_fixed_weak_new} and are thus symmetric and lower bounded.
 It remains to prove that they are also surjective. This follows from Theorem 3.2.1.3 in \cite{grisvard}. \eproof

\section{Convergence of resolvents}

 \noindent In this section, we prove results on the convergence of our resolvents, which will be   central in this work.
  In what follows, we denote by $H^{-s}(\Omega)$ the dual space
of $H^{s}(\Omega)$.

\begin{teo}\label{convop} If $s> 2$ then
$(A_{\epsilon} - A)u \to 0$ in $H^{-s} (\Omega)$, uniformly for  $u$ in bounded subsets of   $H^{1}(\Omega)$.
\end{teo}
\proof  If  $u, \psi \in H^{1}(\Omega)$, we obtain

\begin{align*}
\langle (A_{\epsilon} - A)u, \psi \rangle_{-1,1} & =  \int_{\Omega} \frac{\partial u}{\partial x_{1}}(x) \frac{\partial \psi}{\partial x_{1}} (x) (|Jh_{\epsilon} (x)| - 1)dx \\
& +  \int_{\Omega} \frac{-\varphi(x_{2}, \epsilon) cos(x_{1}/ \epsilon) \frac{1}{\epsilon}}{1 + \frac{\partial \varphi}{ \partial x_{2}}(x_{2}, \epsilon) sin(x_{1}/ \epsilon)} \frac{\partial u}{\partial x_{1}}(x) \frac{\partial \psi}{\partial x_{2}} (x) | Jh_{\epsilon} (x)| \, d\,x  \\
& +  \int_{\Omega} \frac{-\varphi(x_{2}, \epsilon) cos(x_{1}/ \epsilon) \frac{1}{\epsilon}}{1 + \frac{\partial \varphi}{ \partial x_{2}}(x_{2}, \epsilon) sin(x_{1}/ \epsilon)} \frac{\partial u}{\partial x_{2}}(x) \frac{\partial \psi}{\partial x_{1}}(x) |Jh_{\epsilon}(x)| \, d\,x \\
& +  \int_{\Omega} \frac{\left(\varphi(x_{2}, \epsilon) cos(x_{1}/ \epsilon) \frac{1}{\epsilon}\right)^{2}}{\left(1 + \frac{\partial \varphi}{ \partial x_{2}}(x_{2}, \epsilon) sin(x_{1}/ \epsilon)\right)^{2}} \frac{\partial u}{\partial x_{2}}(x) \frac{\partial \psi}{\partial x_{2}}(x)|Jh_{\epsilon}(x)| \, d\,x  \\
& +  \int_{\Omega} \frac{\partial u}{\partial x_{2}}(x) \frac{\partial \psi}{\partial x_{2}} (x) \left( \frac{1}{\left(1 + \frac{\partial \varphi}{ \partial x_{2}}(x_{2}, \epsilon) sin(x_{1}/ \epsilon)\right)} - 1\right)dx  \\
& +  a\int_{\Omega} u(x) \psi (x) (|Jh_{\epsilon} (x)| - 1)dx.
\end{align*}

\indent Now, if we choose $ s > 2  $, we have  the continuous embedding  $H^{s-1}(\Omega) \hookrightarrow C^{0,\mu}(\overline{\Omega})$ for some  $\mu > 0$ and, therefore
 $ \| \psi \|_{H^s(\Omega) } \leq 1  \to
 \frac{\partial \psi}{\partial x_j} (x) \leq C $, for some constant $C$ and $x\in \Omega$.  The result then follows from simple estimates in the above integrals. \eproof

\begin{lema}\label{convresol}
  If  $\lambda $ is in the sector
$
S_{b, \frac{\pi}{4}} =
\{ \lambda \in \mathbb{C} : \frac{\pi}{4} \leq |\arg(\lambda-b)| \leq \pi  \} $, $\lambda \neq b$ and $ f \in H^{-1}(\Omega)$, $(A_{\epsilon} - \lambda)^{-1} f$ converges to  $(A - \lambda)^{-1} f$ in $H^{s}(\Omega), \, -1 \leq s < 1$, as  $\epsilon$ goes to zero.
\end{lema}
\proof  Let  $ u_{\epsilon} =(A_{\epsilon} - \lambda)^{-1} f $,  $u = (A - \lambda)^{-1} f, \, u, \, u_{\epsilon} \in H^{1}(\Omega)$. From Corollary \ref{limitresolunif}, it follows that   $\{u_{\epsilon}\}_{0 \, \leq \, \epsilon \, \leq \, \epsilon_{0}}$   is uniformly bounded in  $H^{1}(\Omega)$ and, therefore, admits a convergent subsequence in  $H^{s}(\Omega)$, for $-1 \leq s < 1$, $\{u_{\epsilon}\}_{0 \, \leq \, \epsilon \, \leq \, \epsilon_{0}}$, which we still denote by  $\{u_{\epsilon}\}_{0 \, \leq \, \epsilon \, \leq \, \epsilon_{0}}$; $u_{\epsilon} \to w$ in  $H^{s}(\Omega)$ when  $\epsilon \to 0$.\\
\indent  We claim that $w = u$. It suffices to show that  $(A - \lambda)w = f$ since $(A - \lambda)$ is injective. We know that  $(A - \lambda)u_{\epsilon} \to (A - \lambda)w$ in  $H^{s - 2}(\Omega)$. Now,  $(A - \lambda)u_{\epsilon} = (A_{\epsilon} - \lambda)u_{\epsilon} + (A - A_{\epsilon})u_{\epsilon} = f + (A - A_{\epsilon})u_{\epsilon} \to f$ in $H^{-s}(\Omega)$ as  $\epsilon \to 0$ by Theorem \ref{convop}. From uniqueness of the limit  $(A - \lambda)w = f$ and thus $w = u$, as claimed.
\eproof
\begin{rem}\label{convresolfrac}
Let $\tilde{X}_{\epsilon}^{\alpha}$ denote  the fractional powers of the
 operators $A_{\epsilon} $, for $ \alpha \geq 0$.
Since the embeddings $H^1(\Omega) = \tilde{X}_{\epsilon}^{1} \hookrightarrow \tilde{X}_{\epsilon}^{\eta}$  and
$   \tilde{X}_{\epsilon}^{\beta} \hookrightarrow H^{-1}(\Omega) = \tilde{X}_{\epsilon}^{0}$  are compact if $ 0\leq \beta,  \eta <1$, we obtain, under the same hypotheses of Lemma \ref{convresol} and exactly the same argument, that
$(A_{\epsilon} - \lambda)^{-1} f$
converges to  $(A - \lambda)^{-1} f$ in
   $\tilde{X}_{\epsilon}^{\eta}, \,$ as
   $\epsilon$ goes to zero, for any $f \in \tilde{X}^{\beta}
   $.
\end{rem}

\begin{teo}\label{convresolnorm} If  $\lambda $ is in the sector
$
S_{b, \frac{\pi}{4}} =
\{ \lambda \in \mathbb{C} : \frac{\pi}{4} \leq |\arg(\lambda-b)| \leq \pi    \} $, $\lambda \neq b,  \, \beta < 1, \, 0 \leq s < 1$, then  $\parallel (A_{\epsilon} - \lambda)^{-1} - (A - \lambda)^{-1} \parallel_{\mathcal{L}(H^{-\beta}(\Omega), H^{s}(\Omega))} \to 0$ as  $\epsilon \to 0$.
 The convergence is uniform for $\lambda $ in a compact set of the
 sector.
\end{teo}
\proof  Since $H^{-\beta}(\Omega) $ is compactly embedded in
 $H^{-1}(\Omega) $, the result, for each $\lambda$ follows from Lemma \ref{convresol} above and from Lemma 4.5 in
\cite{marcone}. Now, if $\lambda$ is in a compact set of the sector, we  also have  $\frac{|\lambda|}{|\lambda-b|} $ bounded by a constant. Thus,
from Corollary  \ref{limitresolunif}  and   the first resolvent identity, we obtain

  \begin{align*} \parallel  &  \left((A_{\epsilon} - \lambda')^{-1}   - (A - \lambda')^{-1} \right) -
 \left((A_{\epsilon} - \lambda)^{-1} - (A - \lambda)^{-1} \right) \parallel_{\mathcal{L}(H^{-\beta}(\Omega), H^{s}(\Omega))}   \\
  & \leq \parallel  \left((A_{\epsilon} - \lambda')^{-1} - (A_\epsilon - \lambda)^{-1} \right) \|_{\mathcal{L}(H^{-\beta}(\Omega), H^{s}(\Omega))} \\
    &     + \|
  \left((A - \lambda)^{-1} - (A - \lambda')^{-1} \right) \parallel_{\mathcal{L}(H^{-\beta}(\Omega), H^{s}(\Omega))}   \\
  & \leq |\lambda'-\lambda| ( \|(A_\epsilon - \lambda')^{-1}    \|_{\mathcal{L}(H^{-\beta}(\Omega), H^{s}(\Omega))}
   \|(A_{\epsilon} - \lambda)^{-1}    \|_{\mathcal{L}(H^{-\beta}(\Omega), H^{s}(\Omega))}  \\
   & + \|(A - \lambda)^{-1}    \|_{\mathcal{L}(H^{-\beta}(\Omega), H^{s}(\Omega))} \|(A - \lambda')^{-1}    \|_{\mathcal{L}(H^{-\beta}(\Omega), H^{s}(\Omega))} )
   \\
   & \leq \left( \frac{2C^{2}|\lambda'||\lambda|}{|\lambda'-b||\lambda-b|}
     \right) |\lambda- \lambda'| \leq
     K |\lambda- \lambda'|,
  \end{align*}
  where $K$ is a constant independent of $\epsilon$.  Using the compacity of the set, the uniformity of the convergence with respect to $\lambda$ follows by a standard argument.
  \eproof
\begin{rem}   \label{convresolnormfrac} As in Remark \ref{convresolfrac}, we can replace, in Theorem \ref{convresolnorm}, the norm in \newline
$  {\mathcal{L}(H^{-\beta}(\Omega), H^{s}(\Omega))} $ by the norm in
$  {\mathcal{L}(\tilde{X}_{\epsilon}^{\beta}, \tilde{X}_{\epsilon}^{\eta} }) $, if
$0 < \beta,\eta < 1. $
\end{rem}
\begin{cor} \label{convauto}
 If   $\sigma(A)$  and $ \sigma(A_{\epsilon})$ denote  the spectral sets of $A$  and $A_{\epsilon}$, respectively,  then $ \lambda \in  \sigma(A)$, if and only if,  there exists a sequence $\lambda_{\epsilon} \in \sigma (A_{\epsilon})$, with     $\lambda_{\epsilon} \to \lambda$.
\end{cor}
\proof Follows from Theorem 6.38 in \cite{teschl}.

\section{Abstract formulation and  local well posedness}

In this section, we  show that problem \eqref{pert_bianca_fixed2_weak_new}  can be formulated in an abstract form  in a fixed scale of Banach spaces and has a unique local solution for any initial data.

\subsection{A scale of Banach spaces associated to the operator
$A_{\epsilon}$}

\indent \par Denote by  $\tilde{X}^{\alpha}_{\epsilon}$ (respectively  $X^{\alpha}_{\epsilon}$) the domain of the fractional powers of  $A_{\epsilon}$ in  $H^{-1}(\Omega)$ (respectively in $L^{2}(\Omega)$).  Then
 $\tilde{X}^{0}_{\epsilon} = H^{-1}(\Omega)$
 (respectively  $X^{0}_{\epsilon} = L^{2}({\Omega})$), \newline
  $\tilde{X}^{1}_{\epsilon} = H^{1}(\Omega)$ (respectively $X^{1}_{\epsilon} = \left\{ u \in H^{2} (\Omega) \, | \, h_{\epsilon}^{*} \frac{\partial}{\partial N_{\Omega_{\epsilon}}}h_{\epsilon}^{*-1} u = 0 \, \text{em} \, \partial \Omega \right\}$) and $X^{\alpha - \frac{1}{2}}_{\epsilon} = \tilde{X}^{\alpha}_{\epsilon}$ for  $\frac{1}{2} \leq \alpha \leq 1$.

 By an abuse of notation, we write $X^{\alpha - \frac{1}{2}}_{\epsilon}$ instead of  $\tilde{X}^{\alpha}_{\epsilon}$ if $0 \leq \alpha \leq \frac{1}{2}$.
  In this way, we obtain an scale of Banach spaces, for $ -1/2 \leq \alpha \leq 1$, with
  $X^{\alpha}_{\epsilon} = H^{2\alpha}(\Omega)$ when $2\alpha$
  is an integer and the inclusion  $ X^{\alpha}_{\epsilon} \hookrightarrow
 X^{\beta}_{\epsilon}$ is compact if $\alpha > \beta$, since $A_{\epsilon}$ has compact resolvent.
We also write simply $\tilde{X}^{\alpha}$ and $X^{\alpha}$ if
$\epsilon = 0$.

\begin{teo} \label{boundresolbeta}  Let  $(A_{\epsilon})_{\beta}: X^{\beta + 1}_{\epsilon} \subset X^{\beta}_{\epsilon} \rightarrow X^{\beta} $
Then  $(A_{\epsilon})_{\beta}$ is sectorial, $0 \leq \epsilon \leq \epsilon_{0}$, and there is a common sector and a common constant in the sectorial inequality for the resolvent operators of $(A_{\epsilon})_{\beta}, \, 0 \leq \epsilon \leq \epsilon_{0}$, $\parallel (\lambda - (A_{\epsilon})_{\beta})^{-1}  \parallel_{\mathcal{L} (X^{\beta}_{\epsilon})} \leq \frac{M}{|\lambda - b|}, \parallel (\lambda - (A_{\epsilon})_{\beta})^{-1}  \parallel_{\mathcal{L} (X^{\beta}_{\epsilon},X^{\beta+1}_{\epsilon}) } \leq \frac{M |\lambda|}{|\lambda - b|}, \, \ \parallel (\lambda - (A_{\epsilon})_{\beta})^{-1}  \parallel_{\mathcal{L} (X^{\beta}_{\epsilon},X^{\beta+\alpha}_{\epsilon}) } \leq \frac{M |\lambda|^{\alpha}}{|\lambda - b|}, \,\forall \lambda \in S_{b, \phi}$, with $M $, $b > 0$, independent of $\epsilon$ and  $\phi \in (0, \frac{\pi}{2})$ .
The associated semigroups   $T_{\epsilon}(t)= \{e^{-(A_{\epsilon})_{\beta} t} \}$ satisfy
$\parallel T_{\epsilon}(t) \parallel_{X^{\beta}_{\epsilon}} \leq Me^{-b t}$,
$\parallel T_{\epsilon}(t) \parallel_{X^{\beta+ 1}_{\epsilon}} \leq
 M \frac{1}{t} e^{-b t}$ and
 $\parallel T_{\epsilon}(t) \parallel_{X^{\beta+ \alpha}_{\epsilon}} \leq
 M_ {\alpha} \frac{1}{t^{\alpha}}  e^{-b t}$, with
  $M_ {\alpha}$ bounded in a compact interval of $[0, \infty[ $.

\end{teo}
\proof
 Let $\delta = 1/2+ \beta$. Then, by the definition of the norms,
    $(A_{\epsilon})^{\delta}$ is an isometry from  $X^{\beta }_{\epsilon} $  into $X^{-1/2}_{\epsilon}$. Also, since
     $(A_{\epsilon})^{\delta} =  (A_{\epsilon})^{-1/2} (A_{\epsilon})^{\beta +1 }$,  $(A_{\epsilon})^{\delta}$ is an isometry from $X^{\beta+1}_{\epsilon}$ into $X^{1/2}_{\epsilon} $.
      Therefore, we have
      $(A_{\epsilon})_{\beta} = (A_{\epsilon})^{-\delta}  A_ {\epsilon} (A_{\epsilon})^{\delta}$.

 \begin{eqnarray*}
\parallel (\lambda - (A_{\epsilon})_{\beta})^{-1} x \parallel_{X^{\beta}_{\epsilon}} & = & \parallel  (\lambda - A_{\epsilon})^{-1} x \parallel_{X^{\beta}_{\epsilon}} \\
& = & \parallel (A_{\epsilon})^{-\delta } (\lambda - A_{\epsilon})^{-1} (A_{\epsilon})^{\delta }
 x \parallel_{X^{\beta}_{\epsilon}} \\
& = & \parallel  (\lambda - A_{\epsilon})^{-1} (A_{\epsilon})^{\delta} x \parallel_{X^{-\frac{1}{2}}_{\epsilon}} \\
& \leq & \frac{M}{|\lambda - b|} \parallel (A_{\epsilon})^{\delta} x \parallel_{X^{-\frac{1}{2}}_{\epsilon}} \\
& = & \frac{M}{|\lambda - b|} \parallel x \parallel_{X^{\beta}_{\epsilon}}.
\end{eqnarray*}
   where $M$ and $b$ are the constants (independent of $\epsilon)$ given
    by Corollary \ref{sectorialh1}.
Applying $  (A_{\epsilon})_{\beta} $   to this estimate, we obtain
 the estimate in $ {\mathcal L}(X^{\beta} , X^{\beta+1}) $.  The estimate  in  $ {\mathcal L}(X^{\beta} , X^{\beta+\alpha}) $ follows by interpolation (see Theorem 1.4.4 in \cite{henry_semilinear}).

   From Theorem 1.3.4 in \cite{henry_semilinear}, $-(A_{\epsilon})_{\beta}$ is
    the infinitesimal generator of an analytic semigroup
    $\{ e^{-(A_{\epsilon})_{\beta}t}\}_{t \, \geq  \, 0}$
    with $\parallel e^{-(A_{\epsilon})_{\beta}t} \parallel \, \leq Me^{-bt}$, $\parallel  (A_{\epsilon})_{\beta} e^{-(A_{\epsilon})_{\beta}t} \parallel \, \leq M\frac{1}{t}e^{-bt}$ where  $M$ is independent of $\epsilon$. The last inequality is proved in \cite{henry_semilinear}, Theorem
1.4.3.
    \eproof

\begin{teo} \label{equiv_norms_alpha} For $ - 1/2 \leq \alpha \leq 1/2$, $\parallel u \parallel_{{X}^{\alpha}_{\epsilon}} \leq K_{1} \parallel u \parallel_{{X}^{\alpha}} \leq K_{2} \parallel u \parallel_{{X}^{\alpha}_{\epsilon}}$, with $K_{1}$ and  $K_{2}$ constants independent of  $\epsilon$. In other words,  ${X}^{\alpha}_{\epsilon} = {X}^{\alpha}$  with equivalent norms uniformly in  $\epsilon$.
\end{teo}

\proof
 A simple estimate shows that  $\parallel A_{\epsilon} u\parallel_{H^{-1}} \leq C' \parallel u \parallel_{H^{1}}$ with  $C$ independent of  $\epsilon$. On the other hand, by Lemma \ref{self_positive}, $\parallel A_{\epsilon} u \parallel_{H^{-1}} \geq C \parallel u \parallel_{H^{1}}$, where $C$ is  independent
of   $\epsilon$.
%

\indent Now, by Theorem  16.1 in \cite{yagi},  the domains
 $ \tilde{X}_{\epsilon}^{\alpha} = {X}_{\epsilon}^{\alpha-1/2}$ of the fractional powers of order  $\alpha$, $0 < \alpha < 1$, of $A_{\epsilon}$  coincide isometrically with the interpolation spaces $[H^{-1}(\Omega), D(A_{\epsilon})]_{\alpha}$.
Thus, from Theorem  1.15 in \cite{yagi}, if  $I$ denotes the inclusion operator, we have, for  $ 0 < \alpha < 1$,
\begin{equation*}
\parallel I \parallel_{\mathcal{L}([ \tilde{X}^{0}_{\epsilon}, \tilde{X}^{1}_{\epsilon}]_{\alpha}, [\tilde{X}^{0}, \tilde{X}^{1}]_{\alpha})} \, \leq \, \parallel I \parallel^{1 - \alpha}_{\mathcal{L}(\tilde{X}^{0}_{\epsilon}, \tilde{X}^{0})} \parallel I \parallel^{\alpha}_{\mathcal{L}(\tilde{X}^{1}_{\epsilon}, \tilde{X}^{1})}.
\end{equation*}

\indent Since  $\tilde{X}^{0}_{\epsilon} = \tilde{X}^{0} = H^{-1}(\Omega)$, we only have to estimate  $\parallel I \parallel_{\mathcal{L}(\tilde{X}^{1}_{\epsilon}, \tilde{X}^{1})}$, where  $\tilde{X}^{1}_{\epsilon} = D(A_{\epsilon})$  and  $\tilde{X}^{1} = D(A)$. We have
\begin{eqnarray*}
\parallel Iu \parallel_{\tilde{X}^{1}}\, &  = \, \parallel u \parallel_{D(A)} \,= \, \parallel Au \parallel_{H^{-1}}\, \leq \, C' \parallel u \parallel_{H^{1}}\, \leq \frac{C'}{C} \parallel A_{\epsilon} u \parallel_{H^{-1}} = \frac{C'}{C} \parallel u \parallel_{D(A_{\epsilon})} \\
 & = \frac{C'}{C} \parallel u \parallel_{\tilde{X}^{1}_{\epsilon}}.
\end{eqnarray*}

\noindent Thus, $\parallel I \parallel_{\mathcal{L}(\tilde{X}^{1}_{\epsilon}, \tilde{X}^{1})} \, \leq \frac{C'}{C}$  and

\begin{equation*}
\parallel u \parallel_{[\tilde{X}^{0}, \tilde{X}^{1}]_{\alpha}} \leq \left( \frac{C'}{C} \right)^{\alpha} \parallel u \parallel_{[ \tilde{X}^{0}_{\epsilon}, \tilde{X}^{1}_{\epsilon}]_{\alpha}},
\textrm{ for }  \   0 \leq \alpha \leq 1.
\end{equation*}

\indent The reverse inequality is obtained in the same way.
Since
$ \tilde{X}_{\epsilon}^{\alpha} = {X}_{\epsilon}^{\alpha-1/2}$,
the result follows immediately.  \eproof

 We can now write the problem  \eqref{pert_bianca_fixed2_weak_new} as an abstract problem in
  the fixed scale of Banach (Hilbert)  spaces  $\{X^{\beta}, -\frac{1}{2} \leq \beta \leq 0\}$,

\begin{equation} \label{problem_abstract_scale}
\left\{
\begin{array}{lll}
u_{t} + (A_{\epsilon})_{\beta}u = (H_{\epsilon})_{\beta} u, \,\,\, t > t_{0}, \\
u (t_{0}) = u_{0} \in X^{\eta},
\end{array}
\right.
\end{equation}

\noindent where $(H_{\epsilon})_{\beta} = H( \cdot, \epsilon) = (F_{\epsilon})_{\beta} + (G_{\epsilon})_{\beta} : X^{\eta} \rightarrow X^{\beta}, \, 0 \, \leq \, \eta \, \leq \, \beta + 1$,\\
$\bullet (F_{\epsilon})_{\beta} = F( \cdot, \epsilon ) : X^{\eta} \rightarrow X^{\beta}$, $0 \leq \epsilon \leq \epsilon_{0}$,
\begin{equation}\label{deff}
<F(u,\epsilon), \phi >_{\beta, - \beta} = \int_{\Omega} f(u) \, \phi\, |Jh_{\epsilon}| \, d\,x, \, \forall \phi \in X^{-\beta}.
\end{equation}

$\bullet (G_{\epsilon})_{\beta} = G( \cdot, \epsilon) : X^{\eta} \rightarrow X^{\beta}$, $0 < \epsilon \leq \epsilon_{0}$,
\begin{equation}\label{defg}
<G(u, \epsilon), \phi >_{\beta, - \beta} = \int_{\partial \Omega} g(\gamma(u)) \, \gamma(\phi) \, |J_{\partial \Omega} h_{\epsilon}| \, d\sigma(x), \, \forall \phi \in X^{-\beta},
\end{equation}

\noindent and

\begin{eqnarray*}
<G(u,0), \phi>_{\beta, - \beta} & = & \int_{\partial \Omega} g(\gamma(u))\gamma(\phi)M_{\pi}(p) d\sigma(x), \, \forall \phi \in X^{-\beta},
\end{eqnarray*}

\noindent where  $\gamma$ is the trace function and 
$M_{\pi}(p) = \frac{1}{\pi} \int_{0}^{\pi} \sqrt{1 + cos^{2}(y)}dy$ -  the average of the  $\pi$ - periodic function $p(y) = \sqrt{1 + cos^{2}(y)}$, in the oscillating part of the boundary and equal to $1$ otherwise. We observe that $J_{\partial \Omega} h_{\epsilon}$ converges in average to $p$ in the oscillatory part of the boundary and converges uniformly to 1 in the remaining parts of it.

 \section{Local Existence}

In order to show the problem (\ref{problem_abstract_scale}) is locally well posed, we need some regularity and growth conditions on the nonlinearities.

(i) $f \in C^{1} (\mathbb{R}, \mathbb{R})$  
and there are real numbers  $\lambda_{1} > 0$ and  $L_{1} > 0$ such that
\begin{equation} \label{hipf}
|\, f(u_{1}) - f(u_{2}) \, | \, \leq \, L_{1} \, ( \, 1 + |\, u_{1} \,|^{\lambda_{1}} + |\, u_{2} \,|^{\lambda_{1}} \,) \, | \, u_{1} - u_{2} \,|, \, \forall u_{1}, \, u_{2} \in \mathbb{R}.
\end{equation}

(ii) $g \in C^{2}(\mathbb{R}, \mathbb{R})$  and there are real numbers  $\lambda_{2} > 0$ and  $L_{2} > 0$ such that:
\begin{equation} \label{hipg}
|\, g(u_{1}) - g(u_{2}) \,| \, \leq \, L_{2} \, ( \, 1 + | \, u_{1} \, |^{\lambda_{2}} + | \, u_{2} \,|^{\lambda_{2}} \,) \, | \, u_{1} - u_{2} \, |, \, \forall u_{1}, \, u_{2} \in \mathbb{R}.
\end{equation}

\begin{lema}\label{bemdeff} Suppose that  $f$ satisfies  (\ref{hipf})  and  $\beta + 1 > \eta > \frac{1}{2} - \frac{1}{2( \lambda_{1} +1)}$. Then, the map $(F_{\epsilon})_{\beta} = F : X^{\eta} \rightarrow X^{\beta}$ is well defined for  $0 \leq \epsilon \leq \epsilon_{0}$  and is bounded in bounded sets of  $X^{\eta}$, uniformly in  $\epsilon$.
\end{lema}

\noindent \proof

\indent Let  $u \in X^{\eta}$ and  $\phi \in X^{-\beta}$,
\begin{eqnarray}
  |< F(u, \epsilon), \phi >_{\beta, -\beta}  |  &  \leq &  \int_{\Omega} |\, f(u) - f(0)\,| \, |\, \phi \,| \, |Jh_{\epsilon}| \, d\,x + \int_{\Omega} |\, f(0) \,| \, |\, \phi \,| \,  | Jh_{\epsilon}| \, d\,x \nonumber \\
& \leq &  K L_{1} \int_{\Omega} |\, u \,| \, |\, \phi \,|\, d\,x  +  K L_{1} \int_{\Omega} |\, u \,|^{\lambda_{1} + 1} |\, \phi \,|\, d\,x \nonumber \\
 & + &   K \int_{\Omega} | \, f(0) \,| \, |\, \phi \,|\, d\,x \nonumber \\
& \leq & K L_{1} \parallel u \parallel_{L^{2}(\Omega)} \, \parallel \phi \parallel_{L^{2}(\Omega)} + K \, L_{1} \parallel u^{\lambda_{1} + 1} \parallel_{L^{2}(\Omega)} \, \parallel \phi \parallel_{L^{2}(\Omega)} \nonumber \\
& + & K \parallel f(0) \parallel_{L^{2}(\Omega)} \, \parallel \phi \parallel_{L^{2}(\Omega)} \nonumber \\
& \leq &  K L_{1} \parallel u \parallel_{L^{2}(\Omega)} \, \parallel \phi \parallel_{L^{2}(\Omega)} + \, K  L_{1} \parallel u \parallel^{\lambda_{1} + 1}_{L^{2(\lambda_{1} + 1)}(\Omega)} \, \parallel \phi \parallel_{L^{2}(\Omega)} \nonumber \\
& + & K  \parallel f(0) \parallel_{L^{2}(\Omega)} \, \parallel \phi \parallel_{L^{2}(\Omega)} \nonumber \\
& \leq & K L_{1} K_{1} K_{2} \parallel u \parallel_{X^{\eta}} \parallel \phi \parallel_{X^{-\beta}} + K L_{1} K_{1} K_{3}^{\lambda_{1} + 1} \parallel u \parallel^{\lambda_{1} + 1}_{X^{\eta}} \parallel \phi \parallel_{X^{-\beta}} \nonumber \\
& + & K  K_{1} \parallel f(0) \parallel_{L^{2}(\Omega)} \parallel \phi \parallel_{X^{-\beta}}, \nonumber
\end{eqnarray}

\noindent where $K = \{ \sup |Jh_{\epsilon}(x)| \, | \, x \in \Omega, \, 0 \leq \epsilon \leq \epsilon_{0} \}$,
 $K_{1}, K_{2} \, \text{e} \, K_{3}$ are the constants of the embeddings  $X^{-\beta} \subset L^{2}(\Omega), \, X^{\eta} \subset L^{2} (\Omega) \, \text{and} \, X^{\eta} \subset L^{2(\lambda_{1} + 1)} (\Omega)$, respectively.
%

\indent Therefore, if  $u \in X^{\eta}$ and  $\phi \in X^{-\beta}$,
\begin{equation*}
\parallel F(u, \epsilon) \parallel_{X^{\beta}} \,  \leq \, K L_{1} K_{1} K_{2} \parallel u \parallel_{X^{\eta}} + \, K L_{1} K_{1} K_{3}^{\lambda_{1} + 1} \parallel u \parallel^{\lambda_{1} + 1}_{X^{\eta}} + \, K K_{1} \parallel f(0) \parallel_{L^{2}(\Omega)}.
\end{equation*}

\eproof

\begin{lema}\label{flipschitz} Suppose the conditions of Lemma \ref{bemdeff} are met  and let  $p$ and $q$ conjugated numbers with  $\frac{1}{\lambda_{1}} < p < \infty$. Then, if  $\beta + 1 > \eta > \max \left\lbrace \frac{1}{2} - \frac{1}{2p \lambda_{1}} , \frac{1}{2} - \frac{1}{2q} \right\rbrace$, the map $F$ is  Lipschitz continuous in  bounded subsets  of $ X^{\eta}$    , with   Lipschitz constant independent of $\epsilon$, for $0 \leq \epsilon \leq \epsilon_{0}$, and some positive $\epsilon_{0}$.
\end{lema}

\proof
\indent Let  $u_{1}, \, u_{2} \in X^{\eta}$ and  $\phi\in X^{-\beta}$. We have
\begin{eqnarray*}
  & & | <  F (u_{1}, \epsilon)  -  F(u_{2}, \epsilon), \phi >_{\beta, -\beta} | \\
  & \leq  & \int_{\Omega} L_{1} \, (\, 1 + |\, u_{1} \,|^{\lambda_{1}} + |\, u_{2} \,|^{\lambda_{1}} \, ) \, |\, u_{1} - u_{2} \,| \, |\, \phi \,|\, |Jh_{\epsilon}| \, d\,x\\
& \leq &  K L_{1} \, \left(\int_{\Omega} ( \, 1 + | \, u_{1} \,|^{\lambda_{1}} + |\, u_{2} \,|^{\lambda_{1}} \,)^{2} \, | \, u_{1} - u_{2} \, |^{2} dx \right)^{\frac{1}{2}} \, \parallel \phi \parallel_{L^{2} (\Omega)} \\
& \leq &  K L_{1} \, \left(\int_{\Omega} (\, 1 + | \, u_{1} \, |^{\lambda_{1}} + |\, u_{2} \,|^{\lambda_{1}} \, )^{2p}\, d\,x \right)^{\frac{1}{2p}} \, \left(\int_{\Omega} |\, u_{1} - u_{2} \,|^{2q}\, d\,x \right)^{\frac{1}{2q}} \, \parallel \phi \parallel_{L^{2}(\Omega)} \\
& \leq  & K L_{1} \, ( \, | \, \Omega \,|^{\frac{1}{2p}} + \parallel u_{1}^{\lambda_{1}} \parallel_{L^{2p}(\Omega)} + \parallel u_{2}^{\lambda_{1}} \parallel_{L^{2p}(\Omega)}) \, \parallel u_{1} - u_{2} \parallel_{L^{2q}(\Omega)} \, \parallel \phi \parallel_{L^{2}(\Omega)} \\
& \leq & K L_{1} \, ( \, |\, \Omega \,|^{\frac{1}{2p}} + \parallel u_{1} \parallel^{\lambda_{1}}_{L^{2p \lambda_{1}}(\Omega)} + \parallel u_{2} \parallel^{\lambda_{1}}_{L^{2p \lambda_{1}}(\Omega)} ) \, \parallel u_{1} - u_{2} \parallel_{L^{2q}(\Omega)} \, \parallel \phi \parallel_{L^{2}(\Omega)} \\
& \leq & K L_{1} \, ( \, |\, \Omega \,|^{\frac{1}{2p}} + K_{4}^{\lambda_{1}} \parallel u_{1} \parallel^{\lambda_{1}}_{X^{\eta}} + \, K_{4}^{\lambda_{1}} \parallel u_{2} \parallel^{\lambda_{1}}_{X^{\eta}} ) \, K_{5} \parallel u_{1} - u_{2} \parallel_{X^{\eta}} K_{1} \parallel \phi \parallel_{X^{-\beta}},
\end{eqnarray*}

\noindent where $K = \sup \{ \, |Jh_{\epsilon} (x)| \, | \, x \in \Omega, 0 \leq \epsilon \leq \epsilon_{0}\}$,
 $K_{1}$ is the constant of the embedding  $X^{-\beta} \subset L^{2}(\Omega)$ and  $K_{4} \, \text{e} \, K_{5}$ are given by the embeddings  $X^{\eta} \subset L^{2p \lambda_{1}}(\Omega) \, \text{and} \, X^{\eta} \subset L^{2q} (\Omega)$, respectively.\\
\indent  Thus, if  $u_{1} \, \text{and} \, u_{2}$ belong to a bounded subset   $U$ of  $X^{\eta}$, $ \parallel u \parallel_{X^{\eta}} \, \leq \, L, \, \forall u \in U$, then
\begin{equation*}
\parallel F(u_{1}, \epsilon) - F(u_{2}, \epsilon) \parallel_{X^{\beta}} \, \leq \, K L_{1} \, ( \, | \,\Omega \,|^{\frac{1}{2p}} + 2 K_{4}^{\lambda_{1}} L^{\lambda_{1}} \,) \, K_{5} K_{1} \parallel u_{1} - u_{2} \parallel_{X^{\eta}}.
\end{equation*}

\eproof

\indent  To obtain the regularity properties of  $G$, we first need to compute  $|J_{\partial \Omega} h_{\epsilon}|$ in each of the segments of  $\partial \Omega$.
We have (see \cite{bianca_doctoral}, for details):

\begin{equation*}
|J_{\partial \Omega} h_{\epsilon}(x_{1}, x_{2})| =
 \begin{cases}
  \sqrt{1 + cos^{2}(x_{1}/ \epsilon)} & (x_{1}, x_{2}) \in I_{1}= \{ \, (x_{1},1) \, | \, 0 \leq x_{1} \leq 1 \, \},\\
  1 + \frac{\partial \varphi}{\partial x_{2}}(x_{2}, \epsilon) sin(1/\epsilon) & (x_{1}, x_{2}) \in I_{2}= \{ \, (1,x_{2}) \, | \, 0 \leq x_{2} \leq 1 \, \},\\
  1 & (x_{1}, x_{2}) \in I_{3}= \{ \, (x_{1},0) \, | \, 0 \leq x_{1} \leq 1 \, \},\\
  1 & (x_{1}, x_{2}) \in I_{4} = \{ \, (0,x_{2}) \, | \, 0 \leq x_{2} \leq 1 \, \},
 \end{cases}
\end{equation*}

\noindent where  $\epsilon$ is sufficiently small, so that  $1 + \frac{\partial \varphi}{\partial x_{2}}(x_{2}, \epsilon) sin(1/\epsilon) \geq 0$.

\begin{lema} \label{bemdefg}
Suppose that  g satisfies  (\ref{hipg}), $\beta + 1 > \eta > \frac{1}{2} - \frac{1}{4(\lambda_{2} + 1)}$ and $-\frac{1}{2} \leq \beta < - \frac{1}{4}$. Then, for  $0 \leq \epsilon $  sufficiently small,  the map $(G_{\epsilon})_{\beta} = G : X^{\eta} \rightarrow X^{\beta}$  is  well defined and  bounded  in bounded subsets of $X^{\eta}$, uniformly in $\epsilon$.
\end{lema}

\proof Let  $u \in X^{\eta}$, $\phi \in X^{-\beta}$  and  $0 < \epsilon \leq \epsilon_{0}$. Then
\begin{eqnarray*}
 & & |< G(u, \epsilon), \phi >_{\beta, - \beta} |  \leq  \int_{\partial \Omega} |\, g(\gamma (u))\,| \, |\, \gamma (\phi)\,| \, | J_{\partial \Omega}h_{\epsilon}| \, d\sigma (x) \\
 & \leq & \overline{K} \left( \int_{\partial \Omega} |\, g(\gamma (u)) - g(\gamma (0))\,| \, | \, \gamma (\phi) \,| \, d\sigma (x) + \int_{\partial \Omega} | \, g(\gamma (0)) \, | \, | \, \gamma (\phi) \, | \, d\sigma (x) \right) \\
 & \leq  &\overline{K} L_{2} \int_{\partial \Omega} | \, \gamma (u) \,| \, | \, \gamma (\phi) \,| \, d\sigma(x) \, + \overline{K} L_{2} \int_{\partial \Omega} | \, \gamma (u) \,|^{\lambda_{2} + 1} | \, \gamma (\phi)\, | \, d\sigma (x) \\
& +  &\overline{K} \int_{\partial \Omega} | \, g(\gamma (0)) \,| \, |\, \gamma (\phi) \, | \, d\sigma (x) \\
& \leq & \overline{K} L_{2} \parallel \gamma (u) \parallel_{L^{2}(\partial \Omega)} \, \parallel \gamma (\phi) \parallel_{L^{2} (\partial \Omega)}  \\
& + & \overline{K} L_{2} \parallel \gamma (u)^{\lambda_{2} + 1} \parallel_{L^{2}(\partial \Omega)} \, \parallel \gamma (\phi) \parallel_{L^{2} (\partial \Omega)} \\
&  +  &\overline{K} \parallel g(\gamma (0)) \parallel_{L^{2} (\partial \Omega)} \, \parallel \gamma (\phi) \parallel_{L^{2} (\partial \Omega)} \\
 &\leq  &\overline{K} L_{2} \, ( \, \parallel \gamma (u) \parallel_{L^{2}(\partial \Omega)} + \parallel \gamma (u) \parallel^{\lambda_{2} + 1}_{L^{2(\lambda_{2} + 1)}(\partial \Omega)} \, ) \, \parallel \gamma (\phi) \parallel_{L^{2} (\partial \Omega)} \\
& +  & \overline{K} \parallel g(\gamma (0)) \parallel_{L^{2} (\partial \Omega)} \, \parallel \gamma (\phi) \parallel_{L^{2} (\partial \Omega)} \\
& \leq & L_{2} \overline{K} \, \overline{K_{3}} \, \overline{K_{1}} \parallel u \parallel_{X^{\eta}} \, \parallel \phi \parallel_{X^{-\beta}} + \, L_{2} \overline{K} \, \overline{K_{2}}^{\lambda_{2} + 1} \, \overline{K_{1}} \parallel u \parallel_{X^{\eta}} \parallel \phi \parallel_{X^{-\beta}} \\
& +  &\overline{K} \, \overline{K_{1}} \,\parallel g(\gamma (0)) \parallel_{L^{2} (\partial \Omega)} \, \parallel \phi \parallel_{X^{-\beta}},
\end{eqnarray*}

\noindent where $\overline{K} = \{ \sup |J_{\partial \Omega} h_{\epsilon}(x) | \, |  \,x \in \partial \Omega, \, 0 \leq \epsilon \leq \epsilon_{0} \}$ and the other constants  $\overline{K_{1}} , \,  \overline{K_{2}} , \, \overline{K_{3}}$ exist because of the continuity of the embedding  $X^{\eta} \subset H^{s} (\Omega)$, and the trace  $ \gamma : H^{s} (\Omega) \rightarrow L^{2} (\partial \Omega)$ and
$ \gamma : H^{s} (\Omega) \rightarrow L^{2(\lambda_2 + 1)} (\partial \Omega)$,
 when $s = 2 \eta$, $\beta < -\frac{1}{4}$, and $\eta > \frac{1}{4}, \, \eta > \frac{1}{2} - \frac{1}{4(\lambda_{2} + 1)}$.
Thus

\begin{equation*}
\parallel G(u,\epsilon) \parallel_{X^{\beta}} \leq L_{2} \overline{K} \, \overline{K_{3}} \, \overline{K_{1}} \parallel u \parallel_{X^{\eta}}+ \, L_{2} \overline{K} \, \overline{K_{2}}^{\lambda_{2} + 1} \, \overline{K_{1}} \parallel u \parallel_{X^{\eta}} + \overline{K} \, \overline{K_{1}} \,\parallel g(0) \parallel_{L^{2} (\partial \Omega)},
\end{equation*}

\noindent so $(G_{\epsilon})_{\beta}$ is well defined and bounded in bounded sets of  $X^{\eta}$ for  $0 < \epsilon \leq \epsilon_{0}$ and some $\epsilon_0 >0 $.
 The case $\epsilon = 0$ is similar and simpler.  \eproof


\begin{lema} \label{glipschitz} Suppose the conditions of Lemma  \ref{bemdefg} are satisfied, and let  $p$ and  $q$ be conjugated exponents with  $\frac{1}{2 \lambda_{2}} < p < \infty$. Then, if  $\beta + 1 > \eta > \max \left\lbrace \frac{1}{2} - \frac{1}{4p \lambda_{2}}, \frac{1}{2} - \frac{1}{4q} \right\rbrace$ and $-\frac{1}{2} \leq \beta < - \frac{1}{4}$,  the map  $G : X^{\eta} \rightarrow X^{\beta}$ is  Lipschitz continuous in  bounded subsets  of $ X^{\eta}$    , with  Lipschitz constant independent of  $\epsilon$ for $0 \leq \epsilon \leq \epsilon_{0}$, and some positive $ \epsilon_{0}$.
\end{lema}

\proof Let  $u_{1}, \, u_{2} \in X^{\eta}$,  $\phi \in X^{-\beta}$  and  $0 < \epsilon \leq \epsilon_{0}$. Then
\begin{eqnarray*}
& &  | < G(u_{1}, \epsilon) - G(u_{2}, \epsilon), \phi >_{\beta, -\beta} |
  \leq  \int_{\partial \Omega} |\, g(\gamma (u_{1})) - g(\gamma(u_{2}))\,| \, |\, \gamma (\phi)\,| \, | J_{\partial \Omega}h_{\epsilon}| \, d\sigma (x) \\
 & \leq   & \overline{K} L_{2} \int_{\partial \Omega} (\, 1 + | \, \gamma (u_{1})\,|^{\lambda_{2}} + |\, \gamma (u_{2}) \,|^{\lambda_{2}} \,) \, | \, \gamma (u_{1}) - \gamma (u_{2}) \, | \, | \, \gamma (\phi) \,| \, d \sigma (x) \\
& \leq & \overline{K} L_{2} \left( \int_{\partial \Omega} ( \, 1 + |\, \gamma (u_{1}) \,|^{\lambda_{2}} + | \, \gamma (u_{2}) \, |^{\lambda_{2}} \,)^{2} \, | \, \gamma (u_{1}) - \gamma (u_{2}) \, |^{2} d\sigma (x) \right)^{\frac{1}{2}} \parallel \gamma (\phi) \parallel_{L^{2}(\partial \Omega)} \\
& \leq  & \overline{K} L_{2}\parallel \gamma (\phi) \parallel_{L^{2}(\partial \Omega)} \left( \int_{\partial \Omega} (\, 1 + | \, \gamma (u_{1}) \, |^{\lambda_{2}} + |\, \gamma (u_{2}) \,|^{\lambda_{2}} \,)^{2p} d\sigma (x) \right)^{\frac{1}{2p}} \cdot  \\
& &\left(\int_{\partial \Omega} |\, \gamma (u_{1}) - \gamma (u_{2}) \, |^{2q} d\sigma (x) \right)^{\frac{1}{2q}} \\
& \leq & \overline{K} L_{2}\parallel \gamma (\phi) \parallel_{L^{2}(\partial \Omega)} \parallel \gamma (u_{1}) - \gamma (u_{2}) \parallel_{L^{2q} (\partial \Omega)} ( \, |\, \partial \Omega \,|^{\frac{1}{2p}} \\
 & +  &\parallel  \gamma (u_{1}) \parallel^{\lambda_{2}}_{L^{2p \lambda_{2}} (\partial \Omega)}  + \parallel \gamma (u_{2}) \parallel_{L^{2p \lambda_{2}} (\partial \Omega)}^{\lambda_{2}} \,) \\
 &\leq & \overline{K} L_{2}  \overline{K_{1}} \parallel \phi \parallel_{X^{-\beta}} \overline{K_{5}} \parallel u_{1} - u_{2} \parallel_{X^{\eta}} (\, | \, \partial \Omega \,|^{\frac{1}{2p}} + \overline{K_{4}}^{\lambda_{2}} \parallel u_{1} \parallel^{\lambda_{2}}_{X^{\eta}} + \, \overline{K_{4}}^{\lambda_{2}} \parallel u_{2} \parallel^{\lambda_{2}}_{X^{\eta}} \,),
\end{eqnarray*}
\noindent where we have used  continuity of the  embedding and trace function in the appropriate spaces as in Lemma \ref{bemdefg}
 to obtain the inequalities:  $ \parallel \gamma (\phi) \parallel_{L^{2}(\partial \Omega)} \, \leq \, \overline{K_{1}} \parallel \phi \parallel_{X^{-\beta}}$, if  $-\beta > \frac{1}{4}$, $ \parallel \gamma
(u) \parallel_{L^{2p \lambda_{2}} (\partial \Omega)} \, \leq \, \overline{K_{4}} \parallel u \parallel_{X^{\eta}}$ if $ \eta > \frac{1}{2} - \frac{1}{4p \lambda_{2}}$ and  $\parallel \gamma (u) \parallel_{L^{2q} (\partial \Omega)} \, \leq \, \overline{K_{5}} \parallel u \parallel_{X^{\eta}}$, if  $\eta > \frac{1}{2} - \frac{1}{4q}$.\\
\indent Thus, if  $u_{1}, \, u_{2}$ are in  a bounded subset  $U$ of  $X^{\eta}$, $ \parallel u \parallel_{X^{\eta}} \, \leq \, L, \, \forall u \in U$, we have
\begin{equation*}
\parallel G(u_{1}, \epsilon) - G(u_{2}, \epsilon) \parallel_{X^{\beta}} \, \leq \, \overline{K} L_{2}  \overline{K_{1}} (\, |\, \partial \Omega \,|^{\frac{1}{2p}} + 2\overline{K_{4}}^{\lambda_{2}} L^{\lambda_{2}} \,) \overline{K_{5}} \parallel u_{1} - u_{2} \parallel_{X^{\eta}},
\end{equation*}
\noindent with  $\overline{K} = \{ \sup |J_{\partial \Omega} h_{\epsilon}(x) | \, | x \in \partial \Omega, \, 0 \leq \epsilon \leq \epsilon_{0} \}$.
 The case $\epsilon=0$ is similar and simpler. \eproof

\begin{teo} \label{propfg}  Suppose that $f$ and $g$ satisfy (\ref{hipf}) and (\ref{hipg}), respectively, $p$ and  $q$ conjugated exponents with   $\max \left\lbrace \frac{1}{2 \lambda_{2}}, \frac{1}{\lambda_{1}} \right\rbrace < p < \infty$,
$\beta$ and $\eta$ satisfying $-\frac{1}{2} \leq \beta < - \frac{1}{4}$, $\beta + 1 > \eta > \max \left\lbrace \frac{1}{2} - \frac{1}{2( \lambda_{1} +1)}  ,\frac{1}{2} - \frac{1}{2p \lambda_{1}}, \frac{1}{2} - \frac{1}{4(\lambda_{2} + 1)}, \frac{1}{2} - \frac{1}{4p \lambda_{2}}, \frac{1}{2} - \frac{1}{4q}  \right\rbrace$ and  $ 0 \leq \epsilon \leq \epsilon_{0}$. Then, the map $(H_{\epsilon})_{\beta} = (F_{\epsilon})_{\beta} + (G_{\epsilon})_{\beta} :  X^{\eta} \rightarrow X^{\beta}$  takes bounded sets of $ X^{\eta}$ into bounded sets   of $ X^{\beta}$, uniformly in $\epsilon$, and is  locally Lipschitz in $X^{\eta} $ with Lipschitz constant  independent of  $\epsilon$, for $0 \leq \epsilon \leq \epsilon_0$.
\end{teo}

\proof
Follows immediately from Lemmas \ref{bemdeff}, \ref{flipschitz}, \ref{bemdefg} and \ref{glipschitz}. \eproof

\begin{lema} \label{contfgepsilon}
 Suppose the conditions on Lemmas \ref{bemdeff}  and  \ref{bemdefg} are satisfied. Then, the map  $(H_{\epsilon})_{\beta} = (F_{\epsilon})_{\beta} + (G_{\epsilon})_{\beta} :  X^{\eta} \rightarrow X^{\beta}$ is  continuous in $\epsilon$, at the point  $\epsilon = 0$, uniformly for $u$ in bounded sets of $X^{\eta}$.
\end{lema}

\proof Let's analyse separately the functions $ (F_{\epsilon})_{\beta} $ and $ (G_{\epsilon})_{\beta}$.
For each $\phi \in X^{- \beta}$ and $0 < \epsilon \leq \epsilon_{0}$,
 \begin{eqnarray*}
\left| \int_{\Omega} f(u) \, \phi \,|Jh_{\epsilon}| \, dx  - \int_{\Omega} f(u) \, \phi \, dx \, \right| & \leq & \int_{\Omega} |f(u)| \, |\phi| \, ||Jh_{\epsilon}| - 1| \, dx \\
& \leq & ||Jh_{\epsilon}| - 1| \int_{\Omega} |f(u)| \, |\phi| \, dx \\
& \leq & ||Jh_{\epsilon}| - 1| \, (L_{1} K_{1} K_{2} \parallel u \parallel_{X^{\eta}} \, \parallel \phi \parallel_{X^{-\beta}} \\
& + & \, L_{1} K_{1} K_{3}^{\lambda_{1} + 1} \parallel u \parallel^{\lambda_{1} + 1}_{X^{\eta}} \, \parallel \phi \parallel_{X^{-\beta}} \\
& + & K_{1} \parallel f(0) \parallel_{L^{2}(\Omega)} \, \parallel \phi \parallel_{X^{-\beta}}),
\end{eqnarray*}
\noindent where we have used the computations in Lemma \ref{bemdeff}. Thus, for $u$ in bounded sets of $X^{\eta}$, since $Jh_{\epsilon} \rightarrow 1$ uniformly as $\epsilon \rightarrow 0$,
\begin{equation*}
\parallel F(u,\epsilon) - F(u,0) \parallel_{X^{\beta}} \rightarrow 0, \, \text{when} \, \epsilon \rightarrow 0,
\end{equation*}

\noindent with $< F(u,0), \phi>_{\beta, - \beta} \, =  \int_{\Omega} f(u) \, \phi \, dx$.

Now, considering $(G_{\epsilon})_{\beta}$,
\begin{eqnarray*}
<G(u, \epsilon), \phi >_{\beta, - \beta} & = & \int_{\partial \Omega} g(\gamma(u)) \, \gamma(\phi) \, |J_{\partial \Omega} h_{\epsilon}| \, d\sigma(x) \\
& = & \int_{I_{1}} g(\gamma(u)) \, \gamma(\phi) \, |J_{\partial \Omega} h_{\epsilon}| \, d\sigma(x) + \int_{I_{2}} g(\gamma(u)) \, \gamma(\phi) \, |J_{\partial \Omega} h_{\epsilon}| \, d\sigma(x) \\
& + & \int_{I_{3}} g(\gamma(u)) \, \gamma(\phi) \, |J_{\partial \Omega} h_{\epsilon}| \, d\sigma(x) + \int_{I_{4}} g(\gamma(u)) \, \gamma(\phi) \, |J_{\partial \Omega} h_{\epsilon}| \, d\sigma(x) \\
& = & <(G_{1} + G_{2} + G_{3} + G_{4})(u, \epsilon), \phi>_{\beta, - \beta},
\end{eqnarray*}
where

\begin{equation}\label{expressionJbound}
|J_{\partial \Omega} h_{\epsilon}(x_{1}, x_{2})| =
 \begin{cases}
  \sqrt{1 + cos^{2}(x_{1}/ \epsilon)} & (x_{1}, x_{2}) \in I_{1} = \{ \, (x_{1},1), 0 \leq x_{1} \leq 1 \, \},\\
  Jh_{\epsilon} & (x_{1}, x_{2}) \in I_{2} = \{ \, (1,x_{2}), 0 \leq x_{2} \leq 1 \, \},\\
  1 & (x_{1}, x_{2}) \in I_{3} = \{ \, (x_{1},0), 0 \leq x_{1} \leq 1 \, \},\\
  1 & (x_{1}, x_{2}) \in I_{4} = \{ \, (0,x_{2}), 0 \leq x_{2} \leq 1 \, \}.
 \end{cases}
\end{equation}


 The result is immediate in  $I_{3}$ and $I_{4}$, since there is no dependence on $\epsilon$ at these portions of the boundary.

  In $I_2$, we have

\begin{eqnarray*}
 & &\left|   \int_{0}^{1} g(\gamma(u(1,x_{2}))) \gamma(\phi) Jh_{\epsilon} \, dx_{2} - \int_{0}^{1} g(\gamma(u(1,x_{2}))) \gamma(\phi) dx_{2} \right| \\
 & \leq &  \int_{0}^{1} |g(\gamma(u(1,x_{2})))| \, | \gamma(\phi)| \,  |Jh_{\epsilon} - 1| dx_{2}\\
& \leq &  |Jh_{\epsilon} - 1| \int_{0}^{1} |g(\gamma(u(1,x_{2})))| \, | \gamma(\phi)|\, d\,x_{2} \\
& \leq  & |Jh_{\epsilon} - 1| (L_{2} \, \overline{K_{3}} \, \overline{K_{1}} \parallel u \parallel_{X^{\eta}} \, \parallel \phi \parallel_{X^{-\beta}}  +  \, L_{2} \, \overline{K_{2}}^{\lambda_{2} + 1} \, \overline{K_{1}} \parallel u \parallel_{X^{\eta}} \parallel \phi \parallel_{X^{-\beta}} \\
& +  & \, \overline{K_{1}} \,\parallel g(\gamma (0)) \parallel_{L^{2} (I_{2})} \, \parallel \phi \parallel_{X^{-\beta}}),
\end{eqnarray*}

\noindent where we have used the computations in Lemma  \ref{bemdefg}.
Thus, for $u$ in bounded sets of
 $X^{\eta}$, since  $Jh_{\epsilon} \to 1$
 uniformly  as  $\epsilon \to 0$,
\begin{equation*}
\parallel G_{2}(u,\epsilon) - G_{2}(u,0) \parallel_{X^{\beta}} \to 0, \, \text{as} \, \epsilon \to 0,
\end{equation*}

\noindent with $< G_{2}(u,0), \phi>_{\beta, - \beta} =  \int_{0}^{1} g(\gamma(u(1,x_{2}))) \gamma(\phi) dx_{2}$.

\indent It only remains to prove the result in  $I_{1}$.
We know that the function
 $p(x_{1}) = \sqrt{1 + cos^{2}(x_{1})}$
 is periodic with period  $\pi$. From Theorem 2.6 in  \cite{cionarescu},  $p_{\epsilon} (x_{1}) = p(x_{1}/\epsilon) \rightharpoonup M_{\pi} (p)$ weakly* in  $L^{\infty}(0,1)$, that is,
\begin{equation*}
\int_{0}^{1} p_{\epsilon}(x_{1}) \, \varphi(x_{1})\, d\,x_{1} \to \int_{0}^{1} M_{\pi}(p) \, \varphi(x_{1}) \, d{x_{1}}, \, \forall \varphi \in L^{1}(0,1),
\end{equation*}

\noindent where  $M_{\pi} (p) = \frac{1}{\pi} \int_{0}^{\pi} \sqrt{1 + cos^{2}(y)}dy$.  Therefore,
\begin{equation*}
<G_{1}(u, \epsilon), \phi>_{\beta, -\beta} \, \to \, <G_{1}(u,0), \phi>_{\beta, -\beta}, \, \forall \, \phi \in X^{-\beta},
\end{equation*}

\noindent with $<G_{1}(u, 0), \phi>_{\beta, -\beta} = \int_{0}^{1} M_{\pi}(p) g(\gamma(u (x_{1},1))) \gamma(\phi) dx_{1}$.

\indent Now $(G_{\epsilon})_{\beta}$  is bounded in bounded sets of   $X^{\eta}$, uniformly in  $\epsilon$, $-\frac{1}{2} \leq \beta < - \frac{1}{4}$. Thus, if   $B \subset X^{\eta}$ is bounded, $\{ G_{1}(u,\epsilon) : u \in B, 0 \leq \epsilon \leq \epsilon_{0} \}$ is bounded in  $X^{\beta'}$ with  $\beta' > \beta$ and  $-\frac{1}{2} \leq \beta' < - \frac{1}{4}$.
 Since the embedding  $X^{\beta'} \subset X^{\beta}$
  is compact, for any sequence    $\epsilon_{n_{k}} \to 0$, there is a subsequence, which we still denote by $\epsilon_{n_{k}} $, such that  $G_{1}(u,\epsilon_{n_{k}}) \to z $ in  $X^{\beta}$ for each $u \in B$.
%
%

It follows that
 $z = G_{1}(u, 0)$ and  $G_{1}(u, \epsilon_{n_{k}}) \rightarrow G_{1}(u, 0) \, \text{in} \, X^{\beta} \, \text{for each} \, u \in B$.

\indent  Fix a  $\rho > 0$.  If   $\eta' < \eta$,  $B$ is a compact subset of  $X^{\eta'}$ and then, for each  $u \in i(B)$, $i: X^{\eta} \rightarrow X^{\eta'}$, there exists a neighborhood  $V_{u}$  of  $u$ in  $X^{\eta'}$ and  $\delta_{u} > 0$ such that ,
\begin{equation*}
\parallel G_{1}(v,\epsilon) - G_{1}(v,0) \parallel_{X^\beta} \,\leq \, \rho,
\end{equation*}
if $0 \leq \epsilon \leq \delta_{u} \leq \epsilon_{0}$  and  $v \in V_{u}$ since $(G_{\epsilon})_{\beta}$ is locally  Lipschitz continuous in  $u$, uniformly in  $\epsilon$.


\indent Since  $B$ is compact in   $X^{\eta'}$ , there exist $u_{1}, ..., u_{n} \in B$ such that $\cup_{i = 1}^{n} V_{u_{i}} \supset i(B)$. Let  $\delta = \min_{1 \, \leq \, i \, \leq \, n} \delta_{u_{i}}$. Then, if  $0 \leq \epsilon \leq \delta$,
\begin{equation*}
\parallel G_{1}(v, \epsilon) - G_{1}(v,0) \parallel_{X^{\beta}} \leq \rho, \, \forall \, v \in  B,
\end{equation*}

%
\noindent and we conclude that
 $(H_{\epsilon})_{\beta}$ is continuous in  $\epsilon = 0$ uniformly in  $u$ for  $u$ in bounded sets of  $X^{\eta}$.
 \eproof

\begin{teo} \label{existsollocal} Suppose that  $f$ and  $g$ satisfy (\ref{hipf})  and  (\ref{hipg}), respectively, let  $\beta$ and  $\eta$ be as in Lemmas  \ref{flipschitz} and  \ref{glipschitz},  $ 0 \leq \epsilon \leq \epsilon_{0}$. Then, for any  $(t_{0}, u_{0}) \in \mathbb{R} \times X^{\eta}$, there is a unique  local solution  of the problem (\ref{problem_abstract_scale}) with $u(t_{0}) = u_{0}$.
\end{teo}

\proof Since  $(A_{\epsilon})_{\beta}$ is sectorial in   $X_{\epsilon}^{\beta}$ with domain  $X_{\epsilon}^{\beta + 1}$ and
  $(H_{\epsilon})_{\beta}: X_{\epsilon}^{\eta} \to X_{\epsilon}^{\beta}$  is locally Lipschitz, the result follows from Theorem
   3.3.3 in \cite{henry_semilinear}. \eproof

 \begin{lema} \label{aproxsemigrupolinear} Let  $\beta$ and $\eta$
 be as in Theorem   \ref{existsollocal} and, additionally,  $\eta < \frac{1}{2}$ and $  0 \leq t \leq T$. Then

  \begin{align*}
\| (e^{-(A_{\epsilon})_{\beta}t} - e^{-(A)_{\beta}t}u)\|_{X^\eta}
 =
 \frac{1}{t^{\eta-\beta}} \Phi(\epsilon),
 \end{align*}
 where
 $\Phi(\epsilon)   $ does not depend on $t$, and
 $\Phi(\epsilon) \to 0  $ as $ \epsilon \to 0$,
 uniformly for  $u$  in bounded sets  of $  X^{\beta} $.
\end{lema}

\proof
By Theorem \ref{boundresolbeta},   there is a common sector
$
S_{b, \frac{\pi}{4}} =
\{ \lambda \in \mathbb{C} : \frac{\pi}{4} \leq |\arg(\lambda-b)| \leq \pi    \} $ and a common constant in the sectorial inequality for the resolvent operators  and the estimates:   $\parallel (\lambda - (A_{\epsilon})_{\beta})^{-1}  \parallel_{\mathcal{L} (X^{\beta}_{\epsilon})} \leq \frac{M}{|\lambda - b|},
\  \,  \ \parallel (\lambda - (A_{\epsilon})_{\beta})^{-1}  \parallel_{\mathcal{L} (X^{\beta}_{\epsilon},X^{\beta+\alpha}_{\epsilon}) } \leq \frac{M |\lambda|^{\alpha}}{|\lambda - b|}\,$ hold for $ \lambda \in S_{b, \frac{\pi}{4}}$, $0\leq \alpha \leq 1$,  with $M $, $b > 0$, independent of $\epsilon$. 

Also, by Lemma \ref{convresolnorm} and Remark
\ref{convresolnormfrac}, if  $\lambda $ is in the sector
$
S_{b, \frac{\pi}{4}}  $, $\lambda \neq b ,  \, \beta > -1/2, \, \eta < 1/2$, then  $\parallel (A_{\epsilon} - \lambda)^{-1} - (A - \lambda)^{-1} \parallel_{\mathcal{L}(X^{\beta}, X^{\eta})} \to 0$ as  $\epsilon \to 0$.
 The convergence is uniform for $\lambda $ in a compact set  of the
 sector.

 Let  $\Gamma$  be a contour in -$S_{b, \phi}$ with  $ arg(\lambda) \to \pm \theta$ as  $|\lambda | \to \infty$ for some   $\theta \in (\frac{\pi}{2}, \pi - \frac{\pi}{4})$. Then

\begin{align} \label{difference_semigroup}
\| (e^{-(A_{\epsilon})_{\beta}t} - e^{-(A)_{\beta}t}u)\|_{X^\eta} & \leq  \frac{1}{2\pi } \int_{\Gamma} |e^{\lambda t}| \, \| ((\lambda + (A_{\epsilon})_\beta)^{-1} - (\lambda + (A)_{\beta})^{-1})u \|_{\eta} \, d \lambda \nonumber \\
 & = \frac{1}{2\pi  } \frac{1}{t^{\eta-\beta}}
  \int_{\Gamma}   t^{\eta-\beta} \, | e^{\lambda t} | \, \|
   ((\lambda + (A_{\epsilon})_\beta)^{-1} - (\lambda + (A)_{\beta})^{-1})u \|_{\eta} \, d  \, \lambda.
\end{align}

   For the integrand in \eqref{difference_semigroup}, we have


   \begin{align*}
   &    \| t^{\eta-\beta} \,  e^{\lambda t}
   ((\lambda + (A_{\epsilon})_\beta)^{-1} - (\lambda + (A)_{\beta})^{-1})u  \|_{\eta}  \\
    &\leq
   | t^{\eta-\beta} \,  e^{\lambda t}|
\left(   \parallel (\lambda + (A_{\epsilon})_\beta)^{-1}u\|_{\eta} + \| (\lambda + (A)_{\beta})^{-1}u  \|_{\eta} \right) \\
  & \leq   2  t^{\eta-\beta} \, | e^{\lambda t} | \frac{M
 |\lambda |^{\eta- \beta} }{ |\lambda - b| }  \|u\|_{\beta}.
   \end{align*}

    This last expression is integrable in $\Gamma$ since,
     with $\mu = \lambda t$, we obtain
    \begin{align*}
     \int_ {\Gamma}  2  t^{\eta-\beta} \, | e^{\lambda t} | \frac{M
 |\lambda |^{\eta- \beta} }{ |\lambda - b| }  \|u\|_{\beta} d \, \lambda & \leq  2 M  \|u\|_{\beta}\int_ {\Gamma} t^{\eta-\beta} \, | e^{\mu} | \frac{
t |  \mu |^{\eta- \beta} }{ t^{\eta- \beta}|\mu - bt| }   \, \frac{1}{t} \, d \, \mu  \\
 & \leq  2 M  \|u\|_{\beta}\int_ {\Gamma}  \, | e^{\mu} | \frac{
 |  \mu |^{\eta- \beta} }{ |\mu - bt| }    \, d \, \mu .
   \end{align*}

    Since
   $ |e^{\lambda t}| \, \| ((\lambda + (A_{\epsilon})_\beta)^{-1} - (\lambda + (A)_{\beta})^{-1})u \|_{\eta}  \to 0$ as $\epsilon \to
    0$, it follows by Lebesgue's Dominated Convergence Theorem that

     \begin{align*}
\| (e^{-(A_{\epsilon})_{\beta}t} - e^{-(A)_{\beta}t}u)\|_{X^\eta}
 = \frac{1}{t^{\eta-\beta}} \Phi(\epsilon),
 \end{align*}
 where
 $\Phi(\epsilon)   $ does not depend on $t$, and
 $\Phi(\epsilon) \to 0  $ as $ \epsilon \to 0$.
  \eproof

\section{Global Existence}

In  this section, we prove that the solutions of the problem  (\ref{problem_abstract_scale}) are globally defined under some additional conditions on  $f$ and  $g$. We also prove the continuity of the solutions with respect to $\epsilon$ at  $\epsilon = 0$, in  $X^{\eta}$, for  $\eta$ in a suitable interval.

 For global existence, we will need the following dissipative conditions on $f$  and $g$:

 There exist constants  $c_{0}, d_{0}, d_{0}^{'}$, with  $d_{0} > d_{0}^{'}$, such that
\begin{equation} \label{hipf2}
\lim_{ |\, u \,| \to \infty} \sup \frac{f(u)}{u} \, \leq \, c_{0},
\end{equation}

\noindent

\begin{equation} \label{hipg2}
\lim_{ |\,u \,| \to \infty} \sup \frac{g(u)}{u} \, \leq \, d_{0}^{'},
\end{equation}

\noindent  and the first eigenvalue of the problem
\begin{equation} \label{autovalor}
\left\{
\begin{array}{lll}
-|Jh_{\epsilon}| h_{\epsilon}^{*}\Delta_{\Omega_{\epsilon}}h_{\epsilon}^{*-1} u + |Jh_{\epsilon}|(a - c_{0})u = \mu u \,\, \text{em} \, \Omega, \\
h_{\epsilon}^{*} \frac{\partial}{\partial N_{\Omega}} h_{\epsilon}^{*-1}u = d_{0} u \, \, \text{em} \, \partial \Omega
\end{array}
\right.
\end{equation}

\noindent is greater than a positive constant $B$ independent of  $\epsilon$, $0 \leq \epsilon \leq \epsilon_{0}$.

\begin{rem}
 The first eigenvalue of  (\ref{autovalor}) satisfies
\begin{equation} \label{primeiroautov}
\lambda_{0} (\epsilon) = \inf_{u \, \in \, H^{1}(\Omega)} \frac{<-|Jh_{\epsilon}| h_{\epsilon}^{*}\Delta_{\Omega_{\epsilon}} h_{\epsilon}^{*-1} u + |Jh_{\epsilon}|(a - c_{0})u , u>}{\parallel u \parallel_{L^{2}(\Omega)}^{2}}.
\end{equation}
 From $ \lambda_{0} \geq B$, we then  obtain
\begin{eqnarray}
B \int_{\Omega} u^{2} dx  &  \leq &  \int_{\Omega}  h_{\epsilon}^{*} \nabla_{\Omega_{\epsilon}} h_{\epsilon}^{*-1} u (x,t) \cdot h_{\epsilon}^{*} \nabla_{\Omega_{\epsilon}} h_{\epsilon}^{*-1} u (x,t) |Jh_{\epsilon} (x)| \, d\,x  \nonumber \\
 &  -  &   d_{0} \int_{\partial \Omega} u^{2}(x,t) |J_{\partial \Omega} h_{\epsilon}(x)| \, dx \nonumber
 +  a \int_{\Omega} u^{2}(x,t) |Jh_{\epsilon}(x)| \, dx \\ &  -  & c_{0} \int_{\Omega} u^{2}(x,t) |Jh_{\epsilon}(x)| \, dx.\label{desigprimeiroautov}
\end{eqnarray}
\end{rem}

\indent To prove global existence, we work in the natural energy space:  $H^{1}(\Omega)$ ($\eta = \frac{1}{2}$ and $\beta < \frac{-1}{4}$).

\begin{lema} \label{estimativasv} Suppose the hypotheses of Theorem   \ref{existsollocal} hold, $f$ and $g$  satisfy the dissipative conditions  (\ref{hipf2}) and  (\ref{hipg2}). Define $V_{\epsilon} : H^{1} (\Omega) \rightarrow \mathbb{R}$ by
\begin{eqnarray} \label{V}
V_{\epsilon} (u) &  =  & \frac{1}{2} \int_{\Omega} h_{\epsilon}^{*} \nabla_{\Omega_{\epsilon}} h_{\epsilon}^{*-1} u (x,t) \cdot h_{\epsilon}^{*} \nabla_{\Omega_{\epsilon}} h_{\epsilon}^{*-1} u (x,t) |Jh_{\epsilon} (x)| \, d\,x + \frac{a}{2} \int_{\Omega}  u(x,t)^{2} |Jh_{\epsilon}(x)| \, dx \nonumber \\
& - & \int_{\Omega} F(u(x,t))\, |Jh_{\epsilon}(x)| \, d\,x - \int_{\partial \Omega} G(\gamma (u)) |J_{\partial \Omega} h_{\epsilon}(x)| \, dS, \, 0 < \epsilon \leq \epsilon_{0},
\end{eqnarray}
\noindent and
\begin{eqnarray*}
V_{0} (u) &  =  & \frac{1}{2} \int_{\Omega} \nabla_{\Omega} u (x,t) \cdot \nabla_{\Omega} u (x,t) dx + \frac{a}{2} \int_{\Omega}  u(x,t)^{2} dx \\ \nonumber
& - & \int_{\Omega} F(u(x,t))\, d\,x - \int_{\partial \Omega} G(\gamma (u)) M_{\pi}(p)dS ,
\end{eqnarray*}
\noindent where  $a$ is a positive number, $F$ and $G : \mathbb{R} \rightarrow \mathbb{R}$ are primitives of $f$ and  $g$, respectively. Then, if $0 \, \leq \, \epsilon \, \leq \, \epsilon_{0}$, $V_{\epsilon}$ is a Lyapunov functional for the problem (\ref{problem_abstract_scale})  and there are  constants $R_{1}$, $R_{2}$,$\tilde{R}_{1}$ and $\tilde{R}_{2}$ such that,
\begin{equation} \label{Vsup}
V_{\epsilon} (u) \, \leq \, R_{1} \parallel u \parallel_{H^{1} (\Omega)}^{2} + \, R_{2}, \,\,\, \forall u \in H^{1}(\Omega),
\end{equation}
and
\begin{equation} \label{Vinf}
V_{\epsilon} (u) \, \geq \, \tilde{R}_{1} \parallel u \parallel_{H^{1} (\Omega)}^{2} - \, \tilde{R}_{2}, \,\,\, \forall u \in H^{1}(\Omega).
\end{equation}
for $ 0 \leq  \epsilon \leq \epsilon_{0}$.

\end{lema}

\proof
Suppose first  $0 < \epsilon \leq \epsilon_{0}$.  $V_{\epsilon}$
 is clearly continuous and  if  $u$ is a solution of (\ref{problem_abstract_scale}) in  $H^{1}(\Omega)$,
\begin{eqnarray*}
\frac{d}{dt} V_{\epsilon} (u(t))  =  \int_{\Omega} h_{\epsilon}^{*} \nabla_{\Omega_{\epsilon}} h_{\epsilon}^{*-1} u (x,t) \cdot h_{\epsilon}^{*} \nabla_{\Omega_{\epsilon}} h_{\epsilon}^{*-1} u_{t} (x,t) |Jh_{\epsilon} (x)| \, d\,x  \\ + a \int_{\Omega}  u(x,t)u_{t}(x,t) |Jh_{\epsilon}(x)| \, d\,x
 -  \int_{\Omega}  f(u(x,t))u_{t}(x,t) \, |Jh_{\epsilon}(x)| \, d\,x  \\ - \int_{\partial \Omega} g(\gamma (u)) u_{t}(x,t) |J_{\partial \Omega} h_{\epsilon}(x)| \, dS \\
 =  - \int_{\Omega} u^2_{t}(x,t)  \, |Jh_{\epsilon}(x)| \, d\,x. \\
\end{eqnarray*}
 Thus  $V_{\epsilon}$ decreases along solutions of (\ref{problem_abstract_scale})  and
\begin{equation*}
\dot V(u) = \frac{d}{dt} V_{\epsilon} (u(t))  = 0,
\end{equation*}
if, and only if,
$u(t)=u_0$ is an equilibrium.

Then  $V_{\epsilon}$, $0 < \epsilon \leq \epsilon_{0}$,
is a Lyapunov functional for the flow generated by
 (\ref{problem_abstract_scale}).

  Now, from  (\ref{hipf2}) there exist  $\varepsilon_{f} > 0$  and   $ k_{f}$ such that  $f(s) \, \leq \, (\varepsilon_{f} + c_{0})|s| + k_{f}$ for  any $s $. Consequently,
\begin{eqnarray*}
\int_{\Omega} F(u)\, |Jh_{\epsilon}| \, d\,x  & = & \int_{\Omega} \left( \int_{0}^{u} f(s) \, |Jh_{\epsilon}| \, ds \right) d\,x \\
& \leq &  \int_{\Omega} \left(\int_{0}^{u} \varepsilon_{f} s + c_{0}s + k_{f} \, ds \right) \, |Jh_{\epsilon}| \, d\,x \\
& \leq & \frac{c_{0}}{2} \int_{\Omega} | \, u \, |^{2}\, |Jh_{\epsilon}| \, d\,x + k_{0}.
\end{eqnarray*}

Also, by (\ref{hipg2}), there exist $\varepsilon_{g} > 0$ and
 $ k_{g}$ such that  $g(s) \, \leq \, (\varepsilon_{g} + d'_{0})|s| + k_{g}$ for  any $s $ and,
\begin{equation*}
\int_{\partial \Omega} G(\gamma (u)) |J_{\partial \Omega} h_{\epsilon}| \, dS \, \leq \, \frac{d_{0}^{'}}{2} \int_{\partial \Omega} | \, \gamma (u) \, |^{2} |J_{\partial \Omega} h_{\epsilon}| \, dS + k_{0}^{'}.
\end{equation*}
Therefore,
\begin{eqnarray*}
V_{\epsilon} (u) & \leq & \frac{1}{2} \int_{\Omega} h_{\epsilon}^{*} \nabla_{\Omega_{\epsilon}} h_{\epsilon}^{*-1} u (x,t) \cdot h_{\epsilon}^{*} \nabla_{\Omega_{\epsilon}} h_{\epsilon}^{*-1} u (x,t) |Jh_{\epsilon} (x)| \, d\,x \\
 & + &  \frac{a}{2} \int_{\Omega} u(x,t)^{2} |Jh_{\epsilon}(x)| \, d\,x
 +   \frac{|\, c_{0} \,|}{2} \int_{\Omega} |\, u \,|^{2}\, |Jh_{\epsilon}(x)| \, d\,x + |\, k_{0} \,| \\
  & +  & \frac{| \, d_{0}^{'} \,|}{2} \int_{\partial \Omega} | \, \gamma (u) \,|^{2} |J_{\partial \Omega} h_{\epsilon}(x)| \, dS + |\, k_{0}^{'} \,| \\
& \leq & \frac{C'}{2} \parallel u \parallel^{2}_{H^{1}(\Omega)} + \, \frac{ K |\, c_{0} \,| }{2} \parallel u \parallel^{2}_{L^{2}(\Omega)} + \, |\, k_{0} \,| + \frac{\overline{K}|\, d_{0}^{'} \,|}{2} \parallel \gamma(u) \parallel^{2}_{L^{2} (\partial \Omega)} + \, | \, k_{0}^{'} \,| \\
& \leq & \frac{C'}{2} \parallel u \parallel^{2}_{H^{1}(\Omega)} + \, \frac{|\, c_{0} \,|}{2} \parallel u \parallel^{2}_{L^{2}(\Omega)} + \, |\, k_{0} \,| + \frac{\overline{K}|\, d_{0}^{'} \,|}{2} \parallel u \parallel^{2}_{L^{2} (\Omega)} + \, | \, k_{0}^{'} \,|
\end{eqnarray*}
where $C'$ is the constant in Remark (\ref{A_epsilon_equiv_norms}), $K$ is a bound for $|Jh_{\epsilon}|$ and
$\overline{K}$, a bound for   $|J_{\partial \Omega}h_{\epsilon}|$.
We conclude that  $V_{\epsilon}(u)$ is bounded by the norm  $\parallel u \parallel_{H^{1}(\Omega)}$  plus constants depending on the nonlinearities 
and independent of   $0 < \epsilon \leq \epsilon_{0}$.

 On the other hand, if   $\lambda_{0}$ is the first eigenvalue of  (\ref{autovalor}) and  $\lambda_{0}(\epsilon) > B$, $0 \leq \epsilon \leq \epsilon_{0}$, by  (\ref{desigprimeiroautov}),
\begin{eqnarray*}
B \int_{\Omega} u^{2} dx & \leq & \int_{\Omega}  h_{\epsilon}^{*} \nabla_{\Omega_{\epsilon}} h_{\epsilon}^{*-1} u (x,t) \cdot h_{\epsilon}^{*} \nabla_{\Omega_{\epsilon}} h_{\epsilon}^{*-1} u (x,t) |Jh_{\epsilon} (x)| \, d\,x \\
 & - &  d_{0} \int_{\partial \Omega} u^{2}(x,t) |J_{\partial \Omega} h_{\epsilon}(x)| \, dx
  +  a \int_{\Omega} u^{2}(x,t) |Jh_{\epsilon}(x)| \, dx
   \\
   &- & c_{0} \int_{\Omega} u^{2}(x,t) |Jh_{\epsilon}(x)| \, dx.
\end{eqnarray*}

\indent Since  $d_{0} > d_{0}^{'}$ we may find  $d_{0}^{''} \neq 0$
such that
\begin{equation*}
d_{0} > d_{0}^{''} > d_{0}^{'} \,\,\, \text{and} \,\,\, \left(\frac{d_{0}}{d_{0}^{''}} - 1 \right) \frac{C_{2} - c_{0}}{2} + \frac{B}{2} > 0,
\end{equation*}

\noindent where $C_{1}$ and $C_{2}$ are constants that limit superiorly the first and the second term in the expression of $A_{\epsilon}u$, respectively. Thus
\begin{align*}
 & V_{\epsilon} (u)  \\
 & \geq  \frac{1}{2} \int_{\Omega} h_{\epsilon}^{*} \nabla_{\Omega_{\epsilon}} h_{\epsilon}^{*-1} u (x,t) \cdot h_{\epsilon}^{*} \nabla_{\Omega_{\epsilon}} h_{\epsilon}^{*-1} u (x,t) |Jh_{\epsilon} (x)| \, d\,x + \frac{a}{2} \int_{\Omega} |\, u \,|^{2} |Jh_{\epsilon}(x)| \, dx \\
& -  \frac{c_{0}}{2} \int_{\Omega} |\, u \,|^{2} |Jh_{\epsilon}(x)| \, dx - \frac{d_{0}^{'}}{2} \int_{\partial \Omega} |\, \gamma (u) \,|^{2} |J_{\partial \Omega} h_{\epsilon}(x)| \, dS - (k_{0} + k_{0}^{'}) \\
& \geq  \frac{1}{2} \int_{\Omega} h_{\epsilon}^{*} \nabla_{\Omega_{\epsilon}} h_{\epsilon}^{*-1} u (x,t) \cdot h_{\epsilon}^{*} \nabla_{\Omega_{\epsilon}} h_{\epsilon}^{*-1} u (x,t) |Jh_{\epsilon} (x)| \, d\,x + \frac{a}{2} \int_{\Omega} |\, u \,|^{2} |Jh_{\epsilon}(x)| \, dx \\
& -  \frac{c_{0}}{2} \int_{\Omega} |\, u \,|^{2} |Jh_{\epsilon}(x)| \, dx - \frac{d_{0}^{''}}{2} \int_{\partial \Omega} |\, \gamma (u) \,|^{2} |J_{\partial \Omega} h_{\epsilon}(x)| \, dS - (k_{0} + k_{0}^{'}) \\
& =  \frac{d_{0}^{''}}{d_{0}} \left[ \frac{d_{0}}{d_{0}^{''}} \frac{1}{2} \int_{\Omega} h_{\epsilon}^{*} \nabla_{\Omega_{\epsilon}} h_{\epsilon}^{*-1} u (x,t) \cdot h_{\epsilon}^{*} \nabla_{\Omega_{\epsilon}} h_{\epsilon}^{*-1} u (x,t)| Jh_{\epsilon} (x)| \, d\,x \right] \\
 & +  \frac{d_{0}^{''}}{d_{0}} \left[  \frac{d_{0}}{d_{0}^{''}} \frac{a}{2} \int_{\Omega} |\, u \,|^{2} |Jh_{\epsilon} (x)| \, dx \right]\\
& +  \frac{d_{0}^{''}}{d_{0}} \left[ - \frac{d_{0}}{d_{0}^{''}}  \frac{c_{0}}{2} \int_{\Omega} |\, u \,|^{2} |Jh_{\epsilon} (x)| \, dx - \frac{d_{0}}{2} \int_{\partial \Omega} |\,\gamma (u) \,|^{2} |J_{\partial \Omega} h_{\epsilon} (x)| \, dS - \frac{d_{0}}{d_{0}^{''}} (k_{0} + k_{0}^{'}) \right] \\
& =  \frac{d_{0}^{''}}{d_{0}} \left[ \left(\frac{d_{0}}{d_{0}^{''}} - 1 \right) \frac{1}{2} \int_{\Omega} h_{\epsilon}^{*} \nabla_{\Omega_{\epsilon}} h_{\epsilon}^{*-1} u (x,t) \cdot h_{\epsilon}^{*} \nabla_{\Omega_{\epsilon}} h_{\epsilon}^{*-1} u (x,t) |Jh_{\epsilon} (x)| \, dx \right] \\
& +  \frac{d_{0}^{''}}{d_{0}} \left[  \frac{1}{2} \int_{\Omega} h_{\epsilon}^{*} \nabla_{\Omega_{\epsilon}} h_{\epsilon}^{*-1} u (x,t) \cdot h_{\epsilon}^{*} \nabla_{\Omega_{\epsilon}} h_{\epsilon}^{*-1} u (x,t) |Jh_{\epsilon} (x)| \, d\,x
\right] \\
 & + \frac{d_{0}^{''}}{d_{0}}  \left[ \left(\frac{d_{0}}{d_{0}^{''}} - 1 \right) \frac{a}{2} \int_{\Omega} |\, u \,|^{2} |Jh_{\epsilon} (x)| \, dx \right]\\
& +  \frac{d_{0}^{''}}{d_{0}} \left[ \frac{a}{2} \int_{\Omega} |\, u \,|^{2} |Jh_{\epsilon} (x)| \, dx - \left(\frac{d_{0}}{d_{0}^{''}} - 1 \right)  \frac{c_{0}}{2} \int_{\Omega} |\, u \,|^{2} |Jh_{\epsilon} (x)| \, dx - \frac{c_{0}}{2} \int_{\Omega} |\, u \,|^{2} |Jh_{\epsilon} (x)| \, d\,x \right] \\
& +  \frac{d_{0}^{''}}{d_{0}} \left[ -\frac{d_{0}}{2} \int_{\partial \Omega} |\, \gamma (u) \,|^{2} |J_{\partial \Omega} h_{\epsilon}(x)| \, dS - \frac{d_{0}}{d_{0}^{''}} (k_{0} + k_{0}^{'}) \right]\\
& \geq  \frac{d_{0}^{''}}{d_{0}} \left[ \left(\frac{d_{0}}{d_{0}^{''}} - 1 \right) \frac{C_{1}}{2} \int_{\Omega} |\nabla u \,|^{2}\, d\,x  \right] \\
 & +  \frac{d_{0}^{''}}{d_{0}} \left[
 \left( \left( \frac{d_{0}}{d_{0}^{''}} - 1 \right) \frac{C_{2}  - c_{0}}{2}
  +  \frac{B}{2} \right) \int_{\Omega} |\, u \,|^{2} dx  - \frac{d_{0}}{d_{0}^{''}} (k_{0} + k_{0}^{'}) \right].
\end{align*}

\noindent  This proves \eqref{Vinf}.


The case $\epsilon = 0$ is similar.  \eproof

\begin{lema} There exist constants $R_{1}^{'}$, $R_{2}^{'}$, such that
\begin{equation*}
R_{1}^{'} |\, V (u) \,| \, \leq \, | \, V_{\epsilon} \, (u) \,| \, \leq \, R_{2}^{'}  |\, V (u) \,|, \,\,\, \forall u \in H^{1}(\Omega),
\end{equation*}
where $V = V_{0}$.
\end{lema}

\proof The result follows from the fact that $|Jh_{\epsilon}| \to 1 $,  uniformly as $\epsilon \to 0$, $|J_{\partial \Omega} h_{\epsilon}|$ can also be bounded, independently of $\epsilon$, both upper and lower, and the superior and inferior  bounds of
$h_{\epsilon}^{*} \nabla_{\Omega_{\epsilon}} h_{\epsilon}^{*-1}$ in $H^1(\Omega)$ are independent of $\epsilon$(see Lemma \ref{self_positive} and Remark \ref{A_epsilon_equiv_norms}).
\eproof

\begin{teo} \label{exist_global} Suppose  $\beta, \, \eta, \, f, \, g $ satisfy the hypotheses of Theorem  \ref{existsollocal} and $f$ and $g$ satisfy (\ref{hipf2}) and (\ref{hipg2}), respectively.
Then, the solutions of problem (\ref{problem_abstract_scale})  are globally defined.
\end{teo}

\proof

\indent Consider first the case  $\eta = \frac{1}{2}$. By Theorem \ref{existsollocal}, there exists, for each $(t_{0}, u_{0}) \in \mathbb{R} \times H^{1} (\Omega)$, $T = T( t_{0}, u_{0})$ such that, the problem (\ref{problem_abstract_scale})  has a solution $u$ in $(t_{0}, t_{0} + T)$ with $u(t_{0}) = u_{0}$. Also, from Lemma
\ref{propfg}, it follows that  $(H_{\epsilon})_{\beta}$  takes bounded sets of  $X^{\eta} = H^{1} (\Omega)$ into bounded sets of  $X^{\beta}$.\\
\indent If  $T < \infty$, by Theorem 3.3.4 of \cite{henry_semilinear}, there exists a sequence  $t_{n} \to T^{-}$ such that  $\parallel u(t_{n}) \parallel_{H^{1} (\Omega)} \to \infty$  and thus  $|V_{\epsilon} (u(t_{n}))| \to \infty$ by  Lemma \ref{estimativasv}, which is impossible, since  $V_{\epsilon}$ decreases along orbits.\\
\indent If  $\eta < \frac{1}{2}$, then  $X^{\frac{1}{2}} \subset X^{\eta}$  and, given $u_{0} \in X^{\eta}$, there exists a local solution $u(t), \, t_{0} < t < T$, by Theorem \ref{existsollocal}. By regularization properties of the semigroup, there exists  $ t_{0} < t' < T$ with  $u_{0}^{'}= T_{\epsilon} (t')u_{0} \in X^{\frac{1}{2}}$  and taking  $u_{0}^{'}$ as the new initial condition we have a global solution which must coincide with our  local solution,  by uniqueness.
 Finally, if
  $\eta > \frac{1}{2}$,   the initial condition
  $u_{0} \in X^{\eta}_{\epsilon} \subset
   X^{\frac{1}{2}} $ so there exists a global solution which belongs to $X^{\eta}_{\epsilon}$, by uniqueness and regularization properties of the flow. \eproof


\begin{lema} \label{limsolxeta} Suppose the initial conditions $u_{0}$ of the problem  (\ref{problem_abstract_scale})  lie in a bounded subset  $B$ of  $X^{\eta}$ and $\eta$ satisfies the conditions of Theorem \ref{existsollocal}. Then, there is a bounded set $C \subset H^1(\Omega)=  X^{1/2}$,  such that   the solution $u_{\epsilon}(t, t_0)$
stays in $C$ for all  $t \geq t_0 $ and   $\epsilon $ sufficiently small.

\end{lema}

\proof  If $ \eta \geq 1/2$, then $X^\eta \subset  H^1(\Omega) $, with continuous inclusion  and the result follows immediately from Theorem \ref{exist_global} and properties (\ref{Vsup}) and (\ref{Vinf}) of the Lyapunov functional.

Suppose that $ \eta < 1/2$ and $ \|u_0\|_{\eta} < C $.
We first show that there exists  $\bar{t} = \bar{t}(C,M) > t_0 $  such that  the  solutions  with initial condition $ u_0$ satisfy
 $\| u(t,u_0) \|_{\eta} \leq 2C$.
   Indeed, by Theorem \ref{propfg}, there exists a   constant $K$ such that $ \| (H_{\epsilon})_{\beta}(u) \parallel_{X^{\beta}} \leq K$, if $ \| u \|_ {\eta} < 2C  $.
Thus, we obtain, while  $\| u(t,u_0) \|_{\eta} \leq 2C$,
\begin{align*}
 & \parallel u_{\epsilon} (t) \parallel_{X^{\eta}}\\
  & \leq  \parallel e^{-(A_{\epsilon})_{\beta} (t - t_{0})} u_{0} \parallel_{X^{\eta}} + \int_{t_{0}}^{t} \parallel e^{-(A_{\epsilon})_{\beta} (t - s)} (H_{\epsilon})_{\beta}(u_{\epsilon}(s)) \parallel_{X^{\eta}} ds \\
& \leq  K_1 \parallel (A_{\epsilon})_{\beta}^{\eta - \beta} e^{-(A_{\epsilon})_{\beta} (t - t_{0})} u_{0} \parallel_{X^{\beta}} +
\int_{t_{0}}^{t} K_1\parallel (A_{\epsilon})_{\beta}^{\eta - \beta} e^{-(A_{\epsilon})_{\beta} (t - s)} (H_{\epsilon})_{\beta}(u_{\epsilon}(s)) \parallel_{X^{\beta}} ds \\
& \leq K_1 C_{\eta - \beta}  e^{-b(t - t_{0})} \parallel u_{0} \parallel_{X^{\eta}} + \int_{t_{0}}^{t} K_1 C_{\eta - \beta} (t - s)^{-(\eta - \beta)} e^{-b(t - s)} \parallel (H_{\epsilon})_{\beta}(u_{\epsilon}(s)) \parallel_{X^{\beta}} ds \\
& \leq K_1C_{\eta - \beta}  e^{-b(t - t_{0})} \parallel u_{0} \parallel_{X^{\eta}} + K K_1 C_{\eta - \beta} \int_{t_{0}}^{t}  (t - s)^{-(\eta - \beta)} e^{-b(t - s)} \, d \,s \\
\end{align*}
where  $C_{\eta - \beta}$  is given by Theorem 1.4.3 in \cite{henry_semilinear}, $b$, such that  $Re(\sigma ((A_{\epsilon})_{\beta})) > b > 0$,  is given by Corollary \ref{sectorialh1} and $K_1$   is the constant given by Theorem \ref{equiv_norms_alpha}

 Since the right hand side of the last inequality  is arbitrarily small, uniformly in $\epsilon$ and $u_0$, if $t$ is close enough to $t_0$, the existence of  $\bar{t}$ as claimed  is proved. Now, choosing $\eta =1/2$,  $\beta> -1/2$ in Theorem \ref{propfg}, we can estimate the $X^{1/2}$ norm at $t=\bar{t}$
  as follows
  \begin{align}\label{est_H1}
\parallel u_{\epsilon} (\bar{t}) \parallel_{X^{1/2}} & \leq  \parallel e^{-(A_{\epsilon})_{\beta} (\bar{t} - t_{0})} u_{0} \parallel_{X^{1/2}} + \int_{t_{0}}^{\bar{t}} \parallel e^{-(A_{\epsilon})_{\beta} (\bar{t} - s)} (H_{\epsilon})_{\beta}(u_{\epsilon}(s)) \parallel_{X^{1/2}} ds \nonumber \\
& \leq  K_1 \parallel (A_{\epsilon})_{\beta}^{1/2- \beta} e^{-(A_{\epsilon})_{\beta} (\bar{t} - t_{0})} u_{0} \parallel_{X^{\beta}} \nonumber \\
& + \int_{t_{0}}^{\bar{t}} K_1\parallel (A_{\epsilon})_{\beta}^{1/2 - \beta} e^{-(A_{\epsilon})_{\beta} (\bar{t} - s)} (H_{\epsilon})_{\beta}(u_{\epsilon}(s)) \parallel_{X^{\beta}} ds \nonumber \\
& \leq  K_1 C_{1/2 - \eta}(\bar{t} - t_{0})^{-(1/2 - \eta)}  e^{-b(\bar{t} - t_{0})} \parallel u_{0} \parallel_{X^{\eta}} \nonumber\\
& + \int_{t_{0}}^{\bar{t}} K_1 C_{1/2 - \beta} (\bar{t} - s)^{-(1/2 - \beta)} e^{-b(\bar{t} - s)} \parallel (H_{\epsilon})_{\beta}(u_{\epsilon}(s)) \parallel_{X^{\beta}} ds \nonumber \\
& \leq K_1C_{1/2 - \eta}(\bar{t} - t_{0})^{-(1/2 - \eta)}  e^{-b(\bar{t} - t_{0})} \parallel u_{0} \parallel_{X^{\eta}} \nonumber \\
& + K K_1 C_{1/2 - \beta} \int_{t_{0}}^{\bar{t}}  (\bar{t} - s)^{-(1/2 - \beta)} e^{-b(\bar{t} - s)} \, d \,s
\end{align}
  Since  $X^{1/2} = H^1(\Omega) $, with equivalent norms, the claim  follows immediately from Theorem \ref{exist_global}.


\begin{teo} \label{contsol}
Suppose the hypotheses of \ref{existsollocal} hold and  $\eta < \frac{1}{2}$, $-\frac{1}{2} < \beta$.  Then, the function $\epsilon \rightarrow u(t,t_{0}, u_{0}, \epsilon) = u_{\epsilon}(t) \in X^{\eta}$ is continuous at  $\epsilon = 0 $ uniformly for $u_{0}$ in bounded sets of  $X^{\eta}$ and  $t_{0} \leq t \leq T$.
\end{teo}

\proof
 By Theorem \ref{propfg}, $H_{\epsilon}$ takes bounded sets of   $X^\eta$ into bounded sets of  $X^\beta$ and is locally Lipschitz with Lipschitz constant $independent$ of $\epsilon$.
	By  Lemma \ref{limsolxeta} $\parallel u_{\epsilon}(s) \parallel_{X^\eta} \, \leq K $ and $\parallel u(s) \parallel_{X^\eta} \, \leq K$ for  $t_{0} \leq s \leq T$. Also, by Lemma \ref{aproxsemigrupolinear},
 $$ \varphi(\epsilon) = 	\parallel (e^{-(A_{\epsilon})_{\beta}(t -t_{0})} - e^{-(A)_{\beta}(t -t_{0})})u_{0} \parallel_{X^\eta} = (t - t_{0})^{\beta - \eta} \Phi(\epsilon)$$
 $$\Phi(\epsilon)
 \to 0  \text{ as }   \epsilon \to 0  \text{ uniformly in bounded sets of }  X^{\beta}.  $$
    Thus,
\begin{align}\label{cont_solution}
\parallel u_{\epsilon} (t) - u(t) \parallel_{X^\eta} \, & \leq \, \parallel (e^{-(A_{\epsilon})_{\beta}(t -t_{0})} - e^{-(A)_{\beta}(t -t_{0})})u_{0} \parallel_{X^\eta} +  \nonumber \\
& +  \int_{t_{0}}^{t} \parallel e^{-(A_{\epsilon})_{\beta}(t-s)} ((H_{\epsilon})_{\beta}(u_{\epsilon}(s)) - (H)_{\beta}(u_{\epsilon} (s))) \parallel_{X^\eta} ds \nonumber \\
& +  \int_{t_{0}}^{t} \parallel e^{-(A_{\epsilon})_{\beta}(t-s)} ((H)_{\beta}(u_{\epsilon}(s)) - (H)_{\beta}(u(s))) \parallel_{X^\eta} ds \nonumber  \\
&  +  \int_{t_{0}}^{t} \parallel (e^{-(A_{\epsilon})_{\beta}(t-s)} - e^{-(A)_{\beta}(t-s)}) (H)_{\beta}(u(s)) \parallel_{X^\eta} ds
 \end{align}

\indent  The first term in \eqref{cont_solution} goes to zero  as $\epsilon$ goes to zero uniformly in $u_0$  by Lemma \ref{aproxsemigrupolinear}

\indent For the second term, writing  $\Delta_{\epsilon} = \sup \{ \, \parallel (H_{\epsilon})_{\beta}(u) - (H)_{\beta}(u) \parallel_{X^\beta} \, | \, \parallel u \parallel_{X^\eta} \, \leq \, K \}$ and applying Theorem \ref{boundresolbeta}, we obtain

\begin{eqnarray*}
& & \int_{t_{0}}^{t}  \parallel e^{-(A_{\epsilon})_{\beta}(t-s)} ((H_{\epsilon})_{\beta}(u_{\epsilon}(s)) - (H)_{\beta}(u_{\epsilon}(s))) \parallel_{X^\eta} ds \\
& \leq & \int_{t_{0}}^{t} \parallel (A_{\epsilon})_{\beta}^{\eta - \beta} e^{-(A_{\epsilon})_{\beta}(t-s)} ((H_{\epsilon})_{\beta}(u_{\epsilon}(s)) - (H)_{\beta}(u_{\epsilon}(s))) \parallel_{X^\beta} ds \\
& \leq & \int_{t_{0}}^{t} C_{\eta - \beta} (t - s)^{\beta -\eta}e^{-b(t- s)} \parallel((H_{\epsilon})_{\beta}(u_{\epsilon}(s)) - (H)_{\beta}(u_{\epsilon}(s))) \parallel_{X^\beta} ds \\
& \leq & C_{\eta - \beta} \Delta_{\epsilon} \int_{t_{0}}^{t} (t - s)^{\beta -\eta}e^{-b(t- s)} ds.
\end{eqnarray*}

\noindent  And $\Delta_{\epsilon}$ goes to zero by Lemma \ref{contfgepsilon}.

For the third term
\begin{eqnarray*}
& & \int_{t_{0}}^{t}  \parallel e^{-(A_{\epsilon})_{\beta}(t-s)} ((H)_{\beta}(u_{\epsilon}(s)) - (H)_{\beta}(u(s))) \parallel_{X^\eta} ds \\
 & \leq & \int_{t_{0}}^{t} \parallel (A_{\epsilon})_{\beta}^{\eta - \beta} e^{-(A_{\epsilon})_{\beta}(t-s)} ((H)_{\beta}(u_{\epsilon}(s)) - (H)_{\beta}(u(s))) \parallel_{X^\beta} ds \\
& \leq & \int_{t_{0}}^{t} \, C_{\eta - \beta} (t - s)^{\beta -\eta} e^{-b(t -s)} \parallel ((H)_{\beta}(u_{\epsilon}(s)) - (H)_{\beta}(u(s))) \parallel_{X^\beta} ds \\
& \leq & C_{\eta - \beta} L \int_{t_{0}}^{t} (t - s)^{\beta -\eta} e^{-b(t -s)} \parallel u_{\epsilon}(s) - u(s) \parallel_{X^\eta} ds.
\end{eqnarray*}

Finally, for the fourth term, since $(H)_{\beta}$ takes bounded sets of $X^{\eta}$ into bounded sets of $X^{\beta}$, using Lemma \ref{limsolxeta},

\begin{align*}
 \int_{t_{0}}^{t}  \parallel (e^{-(A_{\epsilon})_{\beta}(t-s)} - e^{-(A)_{\beta}(t-s)}) (H)_{\beta}(u(s)) \parallel_{X^\eta} ds  & \leq
 \int_{t_{0}}^{t}
   (t - t_{0})^{\beta - \eta} \Phi(\epsilon) \, ds  \\
\leq &\Phi(\epsilon)  \int_{t_{0}}^{t}  (t-t_0)^{\beta- \eta} \,  ds.
 \end{align*}

\indent We therefore have, for $ t_0 <  t < T$ and $u_0$ in a bounded set of $X^\eta$,

\begin{align*}
\parallel u_{\epsilon} (t) - u(t) \parallel_{X^\eta} \, \leq \,
 K_1(\epsilon)  + K_2  \int_{t_{0}}^{t} (t - s)^{\beta -\eta} \parallel u_{\epsilon}(s) - u(s) \parallel_{X^\eta} ds,
\end{align*}

\noindent where $K_1(\epsilon) \to 0 $ as $\epsilon \to 0$. It follows from Gronwall's inequality  ( see \cite{henry_semilinear}) that
 $\parallel u_{\epsilon} (t) - u(t) \parallel_{X^\eta} \to  0$ as  $\epsilon \to 0$ uniformly for  $u_{0}$ in bounded sets of  $X^{\eta}$ and  $0 \leq t \leq T$. \eproof

\section{Existence of global attractors}

We show in this section that the flow generated by (\ref{problem_abstract_scale}) has a global attractor in $X^{\eta}$, for appropriate values of $\eta$, starting with the case
$\eta = \frac{1}{2}$.

\indent Recall that a strongly continuous semigroup $T(t): X \to X$, $t \, \geq \, 0$ is a gradient flow, if \\
\indent (i) Each positive bounded orbit $\{ T(t,x), \  t \geq 0 \}$ is precompact.\vspace{.50cm}

\indent(ii) There exists a Lyapunov functional for   $T(t)$, that is,  a continuous function  $V: X \to \mathbb{R}$ with the following properties:\vspace{.50cm}

\indent(1) $V(x)$ is bounded below,\vspace{.50cm}

\indent(2) $V(x) \to \infty$ as $| \, x \, | \to \infty$,\vspace{.50cm}

\indent(3) $V(T(t)x)$ decreases along orbits of $T(t)$, \vspace{.50cm}

\indent (4) If  $x$ is such that    $V(T(t)x) = V(x)$ for all  $t \in \mathbb{R}$, then $x$ is an equilibrium.

\begin{lema} \label{sistgrad}
 The nonlinear semigroup  $T_{\epsilon, \eta, \beta}(t)$
 generated by (\ref{problem_abstract_scale}) with  $\eta = \frac{1}{2}$ and $-\frac{1}{2} < \beta < - \frac{1}{4}$ in  $X^{\frac{1}{2}} = H^{1}(\Omega)$ is a gradient system.
\end{lema}
\proof From Lemma \ref{estimativasv}, we know that   $V_{\epsilon}$ defined by (\ref{V})  is a Lyapunov functional for
 $T_{\epsilon, \eta, \beta}(t)$, with the required properties. Also, by Lemma \ref{propfg}, $(H_{\epsilon})_{\beta}$ takes bounded sets of  $H^{1}(\Omega)$ into bounded sets of $X^{\beta}$ and, since $(A_{\epsilon})_{\beta} = (A_{\epsilon})^{-\delta} (A_{\epsilon}) (A_{\epsilon})^{\delta}, \, \beta = -\frac{1}{2} + \delta, \, 0 < \delta < \frac{1}{4}$, with $(A_{\epsilon})^{\delta}$ an isometry and   $A_{\epsilon}$
  has compact resolvent by Proposition \ref{limitresolunif}, it follows that   $(A_{\epsilon})_{\beta}$ has compact resolvent.
  Thus, from Theorem 3.3.6 in \cite{henry_semilinear},
  the positive bounded  orbits are precompact.
 This proves the claim. \eproof

\begin{lema}\label{limequih1} Suppose  that  the conditions of the local existence theorem \ref{existsollocal} are met  and $f$,   $g$ satisfy the  dissipative conditions  (\ref{hipf2}), (\ref{hipg2}), respectively. Then,
 the set of equilibria of the system generated by (\ref{problem_abstract_scale})  is bounded in $H^{1}(\Omega)$, uniformly in $ \epsilon$,  $  0 \leq \epsilon \leq \epsilon_0$ .
\end{lema}

\noindent \proof  The equilibria of (\ref{problem_abstract_scale}) satisfy, for   $0 \leq \epsilon \leq \epsilon_0$,
\begin{align} \label{eq1}
0  & =  \int_{\Omega} |Jh_{\epsilon}(x)| (h^{*}_{\epsilon} \Delta_{\Omega_{\epsilon}} h_{\epsilon}^{*-1} u) u  dx - a\int_{\Omega} |\, u \,|^{2} |Jh_{\epsilon}(x)| \, dx + \int_{\Omega} f(u)u |Jh_{\epsilon}(x)| \, dx \nonumber \\
& =  - \int_{\Omega} h^{*}_{\epsilon} \nabla_{\Omega_{\epsilon}} h_{\epsilon}^{*-1} u \cdot h^{*}_{\epsilon} \nabla_{\Omega_{\epsilon}} h_{\epsilon}^{*-1} u |Jh_{\epsilon}(x)| \, dx + \int_{\partial \Omega} g(u)u |J_{\partial \Omega} h_{\epsilon}(x)| \, dS \nonumber \\
& -  a\int_{\Omega} |\, u \,|^{2} |Jh_{\epsilon}(x)| \, dx + \int_{\Omega} f(u)u \, |Jh_{\epsilon}(x)| \, dx.
\end{align}

\indent By  (\ref{hipf2}) and (\ref{hipg2}), for any  $\delta > 0$, there exists constant $K_{\delta}$  such that

\begin{equation*}
f(u)u \, \leq \, (c_{0} + \delta)u^{2} + K_{\delta} \,\,\, \text{and} \,\,\, g(u)u \, \leq \, (d_{0}^{'} + \delta)u^{2} + K_{\delta},
\end{equation*}
\indent Using these inequalities and the bounds  $|Jh_{\epsilon}| \leq K$, $|J_{\partial \Omega} h_{\epsilon}| \leq \overline{K}$, we obtain, from (\ref{eq1}),
\begin{eqnarray} \label{eq2}
& & \int_{\Omega}  h^{*}_{\epsilon} \nabla_{\Omega_{\epsilon}} h_{\epsilon}^{*-1} u \cdot h^{*}_{\epsilon} \nabla_{\Omega_{\epsilon}} h_{\epsilon}^{*-1} u |Jh_{\epsilon}(x)| \, dx \nonumber \\
& \leq & - a\int_{\Omega} |\, u \,|^{2} |Jh_{\epsilon}(x)| \, dx  + (\, c_{0} + \delta \,) \int_{\Omega} u^{2} |Jh_{\epsilon}(x)| \, dx + K_{\delta} K \, |\, \Omega \,| \nonumber \\
& + & ( \, d_{0}^{'} + \delta \,) \int_{\partial \Omega} u^{2} |J_{\partial \Omega}h_{\epsilon} (x)| \, dS + K_{\delta} \overline{K} \, |\, \partial \Omega \, | \nonumber \\
& = & - a\int_{\Omega} |\, u \,|^{2} |Jh_{\epsilon} (x)| \, dx  + (\, c_{0} + \delta \,) \int_{\Omega} u^{2} |Jh_{\epsilon}(x)| \, dx \nonumber \\
& + & (\, d_{0}^{'} + \delta \,) \int_{\partial \Omega} u^{2}|J_{\partial \Omega}h_{\epsilon} (x)| \, dS + K_{\delta} ( \, K |\, \Omega \,| + \overline{K} |\, \partial \Omega \,| ).
\end{eqnarray}

\indent Since the first eigenvalue  $\lambda_{0}(\epsilon)$ of  (\ref{autovalor}) is bounded below by $ B > 0$, $0 \leq \epsilon \leq \epsilon_{0}$, we have
\begin{eqnarray} \label{eq3}
B \int_{\Omega} u^{2} dx & \leq & \int_{\Omega}  h_{\epsilon}^{*} \nabla_{\Omega_{\epsilon}} h_{\epsilon}^{*-1} u (x,t) \cdot h_{\epsilon}^{*} \nabla_{\Omega_{\epsilon}} h_{\epsilon}^{*-1} u (x,t) |Jh_{\epsilon} (x)| \, dx  \\
 & - & d_{0} \int_{\partial \Omega} u^{2}(x,t) |J_{\partial \Omega} h_{\epsilon} (x)| \, dS \nonumber \\
& + & a \int_{\Omega} u^{2}(x,t)| Jh_{\epsilon}(x)| \, dx - c_{0} \int_{\Omega} u^{2}(x,t) |Jh_{\epsilon}(x)| \, dx.
\end{eqnarray}

\indent From \eqref{eq2} and \eqref{eq3}, it follows that
\begin{eqnarray} \label{eq4}
B \int_{\Omega} |\, u \,|^{2} dx & \leq & \delta \int_{\Omega} u^{2} |Jh_{\epsilon}(x)| \, dx  \nonumber \\
 & + & (\, d_{0}^{'} - d_{0} + \delta \,) \int_{\partial \Omega} u^{2} |J_{\partial \Omega}h_{\epsilon}(x)| \, dS + K_{\delta} \, ( \, K | \, \Omega \,| + \overline{K}|\, \partial \Omega \,| \,) \nonumber \\
& \leq & K \delta \int_{\Omega} u^{2} dx + \overline{K}(\, d_{0}^{'} - d_{0} + \delta \,) \int_{\partial \Omega} u^{2} dS + K_{\delta} \, ( \, K | \, \Omega \,| + \overline{K}|\, \partial \Omega \,| \,).
\end{eqnarray}

\indent Now, with  $\delta < \min \{ B/K, d_{0} - d_{0}^{'} \}$  and   $l_{\delta} = \min \{ B - K \delta, \overline{K}((\, d_{0} - d_{0}^{'} \,) - \delta) \}$, it follows from \eqref{eq4} that
\begin{equation} \label{eq5}
\int_{\Omega} u^{2} dx + \int_{\partial \Omega} u^{2} dS \, \leq \, \frac{K_{\delta} \, (\, K |\, \Omega \,| + \overline{K}|\, \partial \Omega \,| \,)}{l_{\delta}}.
\end{equation}

\indent  Substituting \eqref{eq5} in  \eqref{eq2},
\begin{align}  \label{eq6}
\int_{\Omega}  &  h^{*}_{\epsilon} \nabla_{\Omega_{\epsilon}} h_{\epsilon}^{*-1} u \cdot h^{*}_{\epsilon} \nabla_{\Omega_{\epsilon}} h_{\epsilon}^{*-1} u |Jh_{\epsilon}(x)| \, dx \nonumber \\  &  \leq  K(\, c_{0} + \delta \,) \int_{\Omega} u^{2} dx + \overline{K}(\, d_{0}^{'} + \delta \,) \int_{\partial \Omega} u^{2} dS \nonumber \\
& +  K_{\delta} \, (\, K |\, \Omega \,|  + \overline{K}|\, \partial \Omega \,| \,) \nonumber \\
& \leq  (\, K c_{0} + \overline{K} d_{0}^{'} + (K + \overline{K})\delta \,) \frac{K_{\delta} \, (\, K |\, \Omega \,| + \overline{K}|\, \partial \Omega \,| \,)}{l_{\delta}} \nonumber \\
& +  K_{\delta} \, ( \, K |\, \Omega \,| + \overline{K}|\, \partial \Omega \,| \,),
\end{align}

 Now, from Lemma  \ref{self_positive} and \eqref{eq6}

%

\begin{align} \label{eq8}
\int_{\Omega} &  |\, \nabla u \,|^{2} dx \nonumber \\
  & \leq \frac{1}{C_{1}} \left((\, K  c_{0} + \overline{K} d_{0}^{'} + (K + \overline{K})\delta \,) \frac{K_{\delta} \, (\, K |\, \Omega \,| + \overline{K}|\, \partial \Omega \,| \,)}{l_{\delta}}  + K_{\delta} \, ( \, K |\, \Omega \,| + \overline{K}|\, \partial \Omega \,| \,) \right).
 \end{align}

\indent The boundness of equilibria in $H^1(\Omega)$ follows from
 \eqref{eq5} and \eqref{eq8}.  \eproof

\begin{teo} \label{existatratorh1} Suppose that  $f$ and $g$ satisfy the conditions of the local existence theorem \ref{existsollocal}, with  $\eta = \frac{1}{2}$, $-\frac{1}{2} < \beta < -\frac{1}{4}$ and the dissipative hypotheses  (\ref{hipf2}) and (\ref{hipg2}), respectively. Then, for  $0 \, \leq \, \epsilon \, \leq \, \epsilon_{0}$, the flow $T_{\epsilon, \frac{1}{2}, \beta} (t, u) = T_{\epsilon} (t, u)$ generated by  (\ref{problem_abstract_scale})
 has a global attractor in  $X^{\frac{1}{2}} = H^{1}(\Omega)$.
\end{teo}

\proof The flow is gradient and the set of equilibria is bounded so, by Theorem 3.8.5 in \cite{hale}, it remains only to prove that it is also asymptotically smooth. Let $B$ be a non - empty, closed and bounded subset of
$X^{\frac{1}{2}} = H^{1}(\Omega)$, with  $T_{\epsilon}(t, u)B \subset B, \, \forall t \, \geq \, 0$, and   $J = \overline{T_{\epsilon}(1)B}$.  Then,  if $t > 1$, $T_{\epsilon}(t)B = T_{\epsilon}(1)T_{\epsilon}(t - 1)B   \subset T_{\epsilon}(1)B \subset J$, so  $J$ absorbs  $B$. It remains to  show that  $J$ is compact.\\
Suppose $B$ is contained in the ball of radius $N$ in
 $X^{\frac{1}{2}}$.
  Then, if $\Re((A_{\epsilon})_{\beta}) > \omega > 0$,  for any $j = T_{\epsilon}(1)b$  with $b \in B$,
\begin{eqnarray*}
\parallel T_{\epsilon}(1) b  \parallel_{X_{\epsilon}^{\frac{1}{2} + \delta}} & \leq & \parallel e^{-(A_{\epsilon})_{\beta}}b \parallel_{X_{\epsilon}^{\frac{1}{2} + \delta}} + \left| \left| \int_{0}^{1} e^{-(A_{\epsilon})_{\beta} (1 -s)} (H_{ \epsilon})_{\beta} (T_{\epsilon}(s)) \, ds \right| \right|_{X_{\epsilon}^{\frac{1}{2} +\delta}} \\
& \leq & \parallel ((A_{\epsilon})_{\beta})^{\delta} e^{-(A_{\epsilon})_{\beta}} b \parallel_{X^{\frac{1}{2}}}  \\
 & +  & \int_{0}^{1} \parallel ((A_{\epsilon})_{\beta})^{1 - \delta} e^{-(A_{\epsilon})_{\beta} (1 -s)} (H_{\epsilon})_{\beta} (T_{\epsilon}(s))  \parallel_{X^{-\frac{1}{2} + 2\delta}}ds\\
& \leq &  C_{\delta} e^{-\omega} \parallel b \parallel_{X^{\frac{1}{2}}}  \\
 & +  & \int_{0}^{1} C_{1 - \delta} (1 - s)^{- (1 - \delta)} e^{- \omega (1 - s)} \parallel (H_{\epsilon})_{\beta} (T_{\epsilon}(s)) \parallel_{X^{-\frac{1}{2} + 2\delta}} \, ds \\
& \leq &  C_{\delta} e^{-\omega} N + C_{1 - \delta}L \int_{0}^{1} (1 - s)^{-(1 - \delta)} e^{- \omega(1 - s)} ds,
\end{eqnarray*}

\noindent where we applied  Corollary \ref{sectorialh1}  and  $L$ is a local Lipschitz  constant for  $(H_{\epsilon})_{\beta}$.
Thus  $T_{\epsilon}(1)B$ is bounded in $X_{\epsilon}^{\frac{1}{2} + \delta}$ and, therefore, precompact in  $X_{\epsilon}^{\frac{1}{2}}$, by compact embedding. It follows that  its closure  $ J $ is compact, as required. \eproof

\begin{teo} \label{existatratorgeral}
Suppose that  $f$ and $g$ satisfy the conditions of the local existence theorem \ref{existsollocal}  and the dissipative hypotheses  (\ref{hipf2}) and (\ref{hipg2}), respectively. Then, for  $0 \, \leq \, \epsilon \, \leq \, \epsilon_{0}$, the flow $T_{\epsilon, \frac{1}{2}, \beta} (t, u) = T_{\epsilon} (t, u)$ generated by  \eqref{problem_abstract_scale} has a global attractor
$\mathcal{A}_{\epsilon}$  in  $X^{\eta}$. Also  $\mathcal{A}_{\epsilon}$ does not depend  on  $\eta$ or $\beta$.
\end{teo}

\proof
The case $\eta =1/2 $  was proved
in Theorem \ref{existatratorh1}. Let us now first suppose that
 $\eta = \eta_{0} < \frac{1}{2}$  and let  $B \subset X^{\eta} $ be a bounded set. We may suppose that  $B = B_{R}$,  the ball  about
 the origin with radius $ R$. By Theorem  \ref{propfg}, there exists  $N$ such that  $\parallel (H_{\epsilon})_{\beta} (u) \parallel_{X^{\beta}} < N$ if  $\parallel u \parallel_{X^{\eta}} < 2R$. Let  $T(t)u_{0} = u(t, t_{0}, u_{0})$ be the solution of \eqref{problem_abstract_scale}  with  $\eta = \eta_{0}, \, u_{0} \in B_{R}$ and $b$ such   $\Re (\sigma ((A_{\epsilon})_{\beta})) > b$. While  $\parallel u(t, t_{0}, u_{0}) \parallel_{X^{\eta}} < 2R$, we have

\begin{eqnarray} \label{bound_time}
\parallel u(t, t_{0}, u_{0}) \parallel_{X^{\eta}} & \leq & \parallel ((A_{\epsilon})_{\beta})^{\eta - \beta} e^{-(A_{\epsilon})_{\beta}( t - t_{0})} u_{0} \parallel_{X^{\beta}} \nonumber \\
& + & \int_{t_{0}}^{t} \parallel ((A_{\epsilon})_{\beta})^{\eta - \beta} e^{-(A_{\epsilon})_{\beta}(t - s)} (H_{\epsilon})_{\beta} u(s) \parallel_{X^{\beta}} ds \nonumber \\
& \leq & C_{\eta - \beta} (t - t_{0})^{\beta - \eta} e^{-b (t - t_{0})} \parallel u_{0} \parallel_{X^{\beta}} \nonumber \\
& + & N C_{\eta - \beta} \int_{t_{0}}^{t} (t - s)^{(\beta -\eta)}e^{-b (t -s)} ds.
\end{eqnarray}

\noindent Let  $T = \sup \{ \, t \, \geq \, t_{0} \, | \, u(s, t_{0}, u_{0}) \in B_{2R}, \, \forall s \, \leq \, t \, \}$ and $\delta > t_{0}$ such that the last expression in \eqref{bound_time} is smaller than  $2R, \, \forall t \in (t_{0}, \delta)$, and $ \epsilon $ small enough.  It follows that  $T \, \geq \, \delta$ and the solutions with initial conditions in  $B_R \subset X^{\eta}$  remain in $B_{2R}$ for  $t_{0} \, \leq \, t \, \leq \, T$. Now, for
$t$ in this interval
\begin{eqnarray*}
\parallel u(t, t_{0}, u_{0}) \parallel_{X^{\frac{1}{2}}} & \leq & \parallel e^{-(A_{\epsilon})_{\beta} (t - t_{0})} u_{0} \parallel_{X^{\frac{1}{2}}} + \int_{t_{0}}^{t} \parallel e^{-(A_{\epsilon})_{\beta} (t - s)} (H_{\epsilon})_{\beta} u(s) \parallel_{X^{\frac{1}{2}}} ds \\
& \leq & \parallel ((A_{\epsilon})_{\beta})^{\frac{1}{2} - \eta} e^{-(A_{\epsilon})_{\beta}(t - t_{0})} u_{0} \parallel_{X^{\eta}} \\
& + & \int_{t_{0}}^{t} \parallel ((A_{\epsilon})_{\beta})^{\frac{1}{2} - \beta} e^{-(A_{\epsilon})_{\beta} (t - s)} (H_{\epsilon})_{\beta} u(s) \parallel_{X^{\beta}} ds \\
& \leq & C_{\frac{1}{2} -\eta} ( t - t_{0})^{-(\frac{1}{2} - \eta)} e^{-b (t - t_{0})} \parallel u_{0} \parallel_{X^{\eta}} \\
 & +  & NC_{\frac{1}{2} - \beta} \int_{t_{0}}^{t} (t - s)^{-(\frac{1}{2} - \beta)} e^{- b (t - s)} ds,
\end{eqnarray*}

\noindent so $T(t)B_{R}$ is in a bounded set of $ X^{\frac{1}{2}}$ for  $t_{0} < t \, \leq \, T$. By Theorem \ref{existatratorh1}, the global attractor of \eqref{problem_abstract_scale} with $\eta = \frac{1}{2}$ attracts  $T(t)B_{R}$ in the norm  $X^{\frac{1}{2}}$, $t_{0} < t < T$, and therefore it also attracts  the ball $B_{R}$ in the norm of  $X^{\eta}$. Since  $\mathcal{A}_{\epsilon}$ is compact and invariant, it must be the attractor in $X^{\eta}$,  $\eta = \eta_{0}$.\\
\indent Suppose now that $\frac{1}{2} < \eta = \eta_{0} < \beta + 1$. Then, if  $B$ is a bounded set in $X^{\eta}_{\epsilon}$, it is also bounded in $X_{\epsilon}^{\frac{1}{2}} = H^{1} (\Omega)$ so it is attracted by  $\mathcal{A}_{\epsilon}$, the attractor of \eqref{problem_abstract_scale} with  $\eta = \frac{1}{2}$,  by the action of $T_{\epsilon, \frac{1}{2}, \beta}(t) = T_{\frac{1}{2}} (t)$, which coincides with the flow in  $X^{\eta}_{\epsilon}$. We now show that   $T_{\frac{1}{2}}(t)$ is continuous from  $X^{\frac{1}{2}}$ into $X^{\eta}_{\epsilon}$ for  $t > 0$.\\
\indent Let $u_{1}, \, u_{2} \in X^{\frac{1}{2}}$ and  $u_{i} (t, t_{0}, u_{i}) = T_{\frac{1}{2}}(t)u_{i}, \, i=1,2$. Then if $t > t_0$
\begin{align*}
& \parallel u_{1}(t, t_{0}, u_{1}) - u_{2}(t, t_{0}, u_{2}) \parallel_{X^{\eta}_{\epsilon}}  \\ &  \leq   \parallel e^{-(A_{\epsilon})_{\beta}(t - t_{0})} (u_{1} - u_{2}) \parallel_{X^{\eta}_{\epsilon}} \\
& +  \int_{t_{0}}^{t} \parallel e^{-(A_{\epsilon})_{\beta} (t -s)} ((H_{\epsilon})_{\beta}(u_{1}(s)) - (H_{\epsilon})_{\beta} (u_{2}(s))) \parallel_{X^{\eta}_{\epsilon}} ds \\
& \leq  \parallel ((A_{\epsilon})_{\beta})^{\eta - \frac{1}{2}} e^{-(A_{\epsilon})_{\beta} (t - t_{0})} (u_{1} - u_{2}) \parallel_{X^{\frac{1}{2}}} \\
& +  \int_{t_{0}}^{t} \parallel ((A_{\epsilon})_{\beta})^{\eta - \beta} e^{-(A_{\epsilon})_{\beta} (t - s)}((H_{\epsilon})_{\beta} (u_{1}(s)) - (H_{\epsilon})_{\beta} (u_{2}(s))) \parallel_{X^{\beta}} ds \\
& \leq  C_{\eta - \frac{1}{2}} (t - t_{0})^{\frac{1}{2}-\eta } e^{-b (t -t_{0})} \parallel u_{1} - u_{2} \parallel_{X^{\frac{1}{2}}} \\
& +  L_{\beta, \frac{1}{2}} C_{\eta - \beta} \int_{t_{0}}^{t} (t - s)^{\beta - \eta} e^{-b (t- s)} \parallel u_{1}(s) - u_{2}(s) \parallel_{X^{\frac{1}{2}}} ds,
\end{align*}

\noindent where $L_{\beta, \frac{1}{2}}$  is a local Lipschitz constant for $(H_{\epsilon})_{\beta}$.\\
\indent  The continuity then follows from Gronwall's inequality
 (see \cite{henry_semilinear} Theo 7.1.1).



\indent Thus, if  $V_{\delta}$ is a small  neighborhood of  $\mathcal{A}_{\epsilon}$ in $X^{\frac{1}{2}}$ containing  $T_{\frac{1}{2}} (t) (B) = T_{\eta}(t)(B)$, then  $T_{\frac{1}{2}}(1)V_{\delta}$ is a  small neighborhood of  $\mathcal{A}_{\epsilon}$ in $X^{\eta}$ containing  $T_{\frac{1}{2}}(t + 1)B = T_{\eta}(t + 1)B$. Since  $\mathcal{A}_{\epsilon} = T_{\frac{1}{2}}(1) \mathcal{A}_{\epsilon} \in X^{\eta}$, it must be the attractor of   $T_{\eta}(t)$ in   $X^{\eta}$.

\section{Upper semicontinuity of the family of attractors}

In this section, we show that the family of attractors of the flow generated by \eqref{problem_abstract_scale} is upper semicontinuous at $\epsilon =0$.

\begin{teo} \label{bound_attractorh1}
Suppose that  $f$ and $g$ satisfy the conditions of the local existence theorem \ref{existsollocal}  and the dissipative hypotheses  (\ref{hipf2}) and (\ref{hipg2}), respectively. Then, the family  $\mathcal{A}_{\epsilon}$ of attractors for the  flow $T_{\epsilon, \frac{1}{2}, \beta} (t, u) = T_{\epsilon} (t, u)$ generated by  \eqref{problem_abstract_scale}  is uniformly bounded in $H^1(\Omega)$.

\end{teo}
\proof
 By Lemma \ref{limequih1}, for $\epsilon$ sufficiently small,
 the set of equilibria   $E_{\epsilon}$ of \eqref{problem_abstract_scale}  is inside the ball  $B_{R} \subset  H^{1}(\Omega)$, with radius  $R$ independent of $\epsilon$. \\
\indent  Now, for $u \in \mathcal{A}_{\epsilon}$, there exists  $t_{u}$ such that  $u = T(t_{u}) u_{0}$ for some $u_{0} \in B_{R}$ (see \cite{hale} Theorem 3.8.5).\\
\indent The claimed boundedness follows then from the properties of the Lyapunov functional $V$, given by Lemma \ref{estimativasv}
 \eproof

\begin{teo} \label{upper_attractor}
Suppose that  $f$ and $g$ satisfy the conditions of the local existence theorem \ref{existsollocal}  and the dissipative hypotheses  (\ref{hipf2}) and (\ref{hipg2}), respectively and, additionally $ \eta < 1/2$.  Then, the family  $\mathcal{A}_{\epsilon}$  of attractors  for the  flow $T_{\epsilon, \frac{1}{2}, \beta} (t, u) = T_{\epsilon} (t, u)$ generated by  \eqref{problem_abstract_scale}  is uppersemicontinuous at $\epsilon=0$.
\end{teo}

\proof   
By Theorem  \ref{bound_attractorh1}, $\cup_{ 0 \, \leq \, \epsilon \, \leq \, \epsilon_{0} \, } \mathcal{A}_{\epsilon}$ is uniformly bounded  in $H^{1}(\Omega)$, so it is also in a bounded set  $B$ of  $X^{\eta}$, for  $\frac{1}{2} - \varsigma < \eta < \frac{1}{2}$, with  $\varsigma$ sufficiently small so as to satisfy the conditions of Theorem  \ref{existsollocal}.
 Then, for any $\delta > 0$, there exists $T= T_{\delta}$  such that $d( T_{0, \eta, \beta}(T)B, \mathcal{A}_{0}) \leq \frac{\delta}{2}$.
Also, by  Theorem \ref{contsol},   $d(T_{\epsilon, \eta, \beta}(T)B,T_{0, \eta, \beta}(T)B) \leq \frac{\delta}{2}$, if  $\epsilon$ is sufficiently small.
Therefore,  $d( T_{\epsilon, \eta, \beta}(T)B, \mathcal{A}_{0}) \leq \delta$.   Since  $\mathcal{A}_{\epsilon} \subset T_{\epsilon, \eta, \beta}(T)B$, it follows that  $d(\mathcal{A}_{\epsilon}, \mathcal{A}_{0}) \leq \delta$, as we wanted to show. \eproof

\section{Continuity of the equilibria}

  We first show that the family of equilibria of \eqref{problem_abstract_scale} is uppersemicontinuous at $ \epsilon=0$.

\begin{teo} \label{contsupeq}  Suppose the hypotheses of  Theorem \ref{existsollocal} and the dissipative hypotheses \eqref{hipf2}, \eqref{hipg2} are fulfilled and additionally that  $\eta < \frac{1}{2}$.Then, the family of equilibria sets $\{E_{\epsilon} \, | \, 0 \leq \epsilon \leq \epsilon_{0} \, \}$ of  (\ref{problem_abstract_scale}) is uppersemicontinuous in  $X^{\eta}$ at $\epsilon =0$.
\end{teo}

\proof  Let $\epsilon_{n}$ be a sequence of positive numbers, with  $\epsilon_{n} \to  0$ and   $u_{\epsilon_{n}} \in E_{\epsilon_{n}}$. Consider a subsequence, which we still denote by
  $u_{\epsilon_{n}}, $ for simplicity.

\indent Since the family of equilibria is uniformly bounded in   $H^{1}(\Omega)$, by Lemma \ref{limequih1}, it  is a compact  set in    $X^{\eta}$. Therefore, there exists a convergent  subsequence  $u_{\epsilon_{n_{k}}} \to u_{0} $. By  continuity  of the flow in $\epsilon$, uniformly in bounded sets,   given by Theorem \ref{contsol}, we obtain for $t>0$,

$$ T_{\epsilon_{n_{k}}, \beta, \eta} (t) u_{\epsilon_{n_{k}}} =
  T_{\epsilon_{n_{k}}, \beta, \eta} (t) u_{\epsilon_{n_{k}}} -
  T_{0, \beta, \eta} (t) u_{\epsilon_{n_{k}}} +
   T_{0, \beta, \eta} (t) u_{\epsilon_{n_{k}}} \to T_{0, \beta, \eta}(t) u_{0}.$$

   On the other hand,

   $$ T_{\epsilon_{n_{k}}, \beta, \eta} (t) u_{\epsilon_{n_{k}}} =
 u_{\epsilon_{n_{k}}} \to u_{0}. $$

\noindent By uniqueness of the limit  $T_{0, \beta, \eta} u_{0} = u_{0}$, which proves the result. \eproof

 For the lower semicontinuity, we will need additional hypotheses on the nonlinearities, namely

\begin{equation} \label{hip_derivf}
f \in C^{1}(\mathbb{R}, \mathbb{R}), \, \text{ with bounded derivatives}.
\end{equation}

 \begin{equation} \label{hip_derivg} g \in C^{2}(\mathbb{R}, \mathbb{R}) \, \text{ with bounded derivatives}.
\end{equation}

\begin{lema} \label{fderivavel} If $f$ satisfies (\ref{hip_derivf}), $-\frac{1}{2} \leq \beta \leq 0$ and $\beta + 1 > \eta > 0$, then the map  $F : X^{\eta} \, \times \,\mathbb{R} \rightarrow X^{\beta}$ defined by (\ref{deff}) is Gateaux-differentiable in  $u$ with Gateaux differential  $F_{u} (u, \epsilon) = \frac{\partial F }{\partial u}(u, \epsilon)w$ given by
\begin{equation}\label{derf}
\left\langle \frac{\partial F }{\partial u}(u, \epsilon)w , \phi \right\rangle_{\beta, - \beta} = \int_{\Omega} f^{'}(u) \, w  \, \phi \, Jh_{\epsilon} \, dx, \, \forall w \in X^{\eta}, \, \forall \phi \in X^{-\beta}.
\end{equation}
Furthermore, $F_{u} (u, \epsilon)$ is linear and bounded.
\end{lema}

\proof If $u, w \in X^{\eta}$ and  $\phi \in X^{-\beta}$,
\begin{eqnarray}\label{diff_gateauxf}
& & \left|\frac{1}{t} \left\langle F(u + tw, \epsilon) - F(u,\epsilon) - t\frac{\partial F }{\partial u}(u, \epsilon)w , \phi \right\rangle_{\beta, - \beta} \right| \nonumber \\
& = & \left| \frac{1}{t} \int_{\Omega} (f(u + tw) - f(u) - tf^{'}(u)w)\, \phi \, |Jh_{\epsilon}| \, dx \right| \nonumber \\
& \leq & \frac{K}{|t|} \int_{\Omega} |f(u + tw) - f(u) - tf^{'}(u)w| \, | \phi | \, dx \nonumber \\
& \leq & \frac{K}{|t|} \left( \int_{\Omega} |f(u + tw) - f(u) - tf^{'}(u)w|^{2} dx \right)^{\frac{1}{2}} \parallel \phi \parallel_{L^{2}(\Omega)} \nonumber \\
& \leq & K_{1} \frac{K}{|t|} \left( \int_{\Omega} |f(u + tw) - f(u) - tf^{'}(u)w|^{2} dx \right)^{\frac{1}{2}} \parallel \phi \parallel_{X^{-\beta}} \nonumber \\
& \leq & K_{1} K \parallel \phi \parallel_{X^{-\beta}} \left( \int_{\Omega} |(f^{'}(u + \bar{t}w) - f^{'}(u))w|^{2} dx \right)^{\frac{1}{2}},
\end{eqnarray}

\noindent where $K = \sup \{ \, |Jh_{\epsilon} (x)| \, | \, x \in \Omega, 0 \leq \epsilon \leq \epsilon_{0} \}$,
$K_{1}$ is the constant of the embedding  $X^{-\beta} \subset L^{2}(\Omega)$, and we have used the Mean Value Theorem in the last passage  with  $0 < \bar{t} < t$. The integrand in \eqref{diff_gateauxf} is bounded by the integrable function  $ 4\parallel f' \parallel_{\infty}^{2} w^{2}$ and, by (\ref{hip_derivf}) goes to $0$ as  $t \to 0$. Thus, by  the Dominated Convergence Theorem
\begin{equation*}
\lim_{t \, \to \, 0} \frac{F(u + tw, \epsilon) - F(u, \epsilon)}{t} = \frac{\partial F}{\partial u}(u, \epsilon)w \,\, \text{in} \, X^{\beta}, \, \forall u, w \in X^{\eta}.
\end{equation*}

\indent Therefore,  $F$ is Gateaux - differentiable, with Gateaux differential given by   (\ref{derf}).
%
Now, for $u, \, w \in X^{\eta}$ and $\phi \in X^{-\beta}$, we have
\begin{eqnarray*}
\left| \left\langle \frac{\partial F }{\partial u}(u, \epsilon)w , \phi \right\rangle_{\beta, - \beta} \right| & \leq & \int_{\Omega}  |f^{'}(u) \, w  \, \phi \,| \, |Jh_{\epsilon}| \, dx \\
& \leq & K \parallel f^{'} \parallel_{\infty} \int_{\Omega}  |w \, \phi | \, dx \\
& \leq &  K \parallel f^{'} \parallel_{\infty} \, \parallel w \parallel_{L^{2}(\Omega)} \, \parallel \phi \parallel_{L^{2}(\Omega)} \\
& \leq & K K_{1} K_{2}  \parallel f^{'} \parallel_{\infty} \, \parallel w \parallel_{X^{\eta}} \, \parallel \phi \parallel_{X^{-\beta}}.
\end{eqnarray*}

Thus,
\begin{equation} \label{limitderivf}
\left| \left| \frac{\partial F }{\partial u}(u, \epsilon)w \right| \right|_{X^{\beta}} \, \leq K K_{1} K_{2}  \parallel f^{'} \parallel_{\infty} \, \parallel w \parallel_{X^{\eta}},
\end{equation}

\noindent where $K_{2}$ is an embedding constant of  $X^{\eta} \subset L^{2}(\Omega)$, proving that $\frac{\partial F}{\partial u}(u, \epsilon)$ is bounded . \eproof

 We now show that the Gateaux differential of
 $F(u,\epsilon)$ is continuous in  $u$ and in $\epsilon$.


\begin{lema} \label{contderf} If  $f$ satisfies (\ref{hip_derivf}), $-\frac{1}{2} \leq \beta \leq 0$ and  $\beta + 1 > \eta > 0$, the Gateaux derivative of $F(u,\epsilon)$ is continuous with respect to  $u$, that is, the map $u \mapsto \frac{\partial F}{\partial u}(u, \epsilon) \in \mathcal{B} (X^{\eta}, X^{\beta})$ is continuous. Furthermore,   $\epsilon \mapsto \frac{\partial F}{\partial u}(u, \epsilon)$ is continuous, uniformly for   $u \in X^{\eta}$.
\end{lema}

\proof
Let  $u_{n} \to  u \in X^{\eta}$  and  $0 < \tilde{\eta} < \eta \leq \beta + 1$. Then, for  $w \in X^{\tilde{\eta}}, \, \phi \in X^{\beta}$,

\begin{align} \label{cont_Gateaux_f_u}
 & \left| \left\langle \left( \frac{\partial F }{\partial u}(u_{n}, \epsilon) - \frac{\partial F}{\partial u}(u, \epsilon) \right)w , \phi \right\rangle_{\beta, - \beta} \right| \nonumber \\
 & \leq  \int_{\Omega} |\, (f^{'}(u) - f^{'}(u_{n})) \, w \, \phi \, | Jh_{\epsilon}| \, | \, dx \nonumber \\
& \leq  K \left( \int_{\Omega} |\, (f^{'}(u) - f^{'}(u_{n})) \, w \, |^{2} dx \right)^{\frac{1}{2}} \parallel \phi \parallel_{L^{2}(\Omega)} \nonumber \\
& \leq  K K_{1} \parallel \phi \parallel_{X^{-\beta}} \left( \int_{\Omega} | \, (f^{'}(u) - f^{'}(u_{n})) \, w \,|^{2} dx \right)^{\frac{1}{2}},
\end{align}

\noindent where $K = \sup \{ \, |Jh_{\epsilon} (x)| \, | \, x \in \Omega, 0 \leq \epsilon \leq \epsilon_{0} \}$, $K_{1}$ is a constant for the embedding
 $X^{-\beta} \subset L^{2}(\Omega)$. The integrand  in \eqref{cont_Gateaux_f_u}  is bounded by the integrable function $ 4 \parallel f^{'} \parallel^{2}_{\infty} w^{2}$  and goes to $0$ as  $u_{n} \to  u \in X^{\eta}$.  From the Dominated Convergence Theorem, it follows that the sequence $\frac{\partial F}{\partial u}(u_{n}, \epsilon)$ converges strongly to  $\frac{\partial F}{\partial u}(u, \epsilon)$ in  $B(X^{\tilde{\eta}}, X^{\beta})$. Since  $X^{\eta}$ is compactly embedded in  $X^{\tilde{\eta}}$, the convergence in  $\mathcal{B}(X^{\eta}, X^{\beta})$ follows from Lemma 4.5 in \cite{marcone} .

The convergence with respect to $\epsilon$ at $\epsilon = 0$ follows from
\begin{eqnarray*}
\left| \left\langle \left( \frac{\partial F}{\partial u}(u, \epsilon) - \frac{\partial F}{\partial u}(u, 0) \right) w, \phi \right\rangle_{\beta, - \beta} \right| & = & \left| \int_{\Omega} f^{'}(u)w \, \phi \, (|Jh_{\epsilon}| - 1) \, dx \right| \\
& \leq & (|Jh_{\epsilon}| - 1) \int_{\Omega} |f^{'}(u) \, w \, \phi| \, dx \\
& \leq & (|Jh_{\epsilon}| - 1) \parallel f^{'} \parallel_{\infty} \, \parallel w \parallel_{L^{2}(\Omega)} \, \parallel \phi \parallel_{L^{2}(\Omega)} \\
& \leq & (|Jh_{\epsilon}| - 1) K_{1}K_{2} \parallel f^{'} \parallel_{\infty} \, \parallel w \parallel_{X^{\eta}} \, \parallel \phi \parallel_{X^{-\beta}}
\end{eqnarray*}
\noindent with $K_{2}$ the constant of the embedding of $X^{\eta}$ into  $L^{2}(\Omega)$, as in the proof of Lemma \ref{fderivavel} . Since $Jh_{\epsilon} \rightarrow 1$  uniformly in $x$ when $\epsilon \rightarrow 0$,
\begin{equation*}
\left| \left| \left( \frac{\partial F}{\partial u}(u, \epsilon) - \frac{\partial F}{\partial u}(u, 0) \right)  \right| \right|_{\mathcal{L}(X^{\eta}, X^{\beta})} \rightarrow 0, \, \text{when} \, \epsilon \rightarrow 0.
\end{equation*}\\
For other values of $\epsilon$, the convergence is immediate. \eproof


\begin{lema} \label{gderivavel}
If  $g$ satisfies (\ref{hipf2}), $\beta + 1 > \eta > \frac{1}{4}$ and  $-\frac{1}{2} \leq \beta < -\frac{1}{4}$, then the  map  $G : X^{\eta} \times \mathbb{R} \to  X^{\beta}$ defined by (\ref{defg}) is  Gateaux - differentiable  in  $u$ with Gateaux differential  $G_{u}(u, \epsilon) = \frac{\partial G}{\partial u}(u, \epsilon)$ given by
\begin{equation}\label{derg}
\left\langle \frac{\partial G}{\partial u}(u, \epsilon)w, \phi \right\rangle_{\beta, - \beta} = \int_{\partial \Omega} g^{'}(\gamma(u))\, \gamma(w) \, \gamma(\phi) \, |J_{\partial \Omega} h_{\epsilon}| \, d\sigma(x),
\end{equation}
 $ \forall w \in X^{\eta}, \, \forall \phi \in X^{-\beta}, \, 0 < \epsilon \leq \epsilon_{0}$
and
\begin{eqnarray}
\left\langle \frac{\partial G}{\partial u}(u, 0)w, \phi \right\rangle_{\beta, - \beta} & = & \int_{I_{1}} g^{'}(\gamma(u))\, \gamma(w) \, \gamma(\phi) \, M_{\pi}(p) \, d\sigma(x) \nonumber \\
& + &\int_{I_{2}} g^{'}(\gamma(u))\, \gamma(w) \, \gamma(\phi) \,d\sigma(x) \nonumber \\
& + & \int_{I_{3}} g^{'}(\gamma(u))\, \gamma(w) \, \gamma(\phi) \,d\sigma(x) \nonumber \\
& + & \int_{I_{4}} g^{'}(\gamma(u))\, \gamma(w) \, \gamma(\phi) \,d\sigma(x), \, \forall w \in X^{\eta}, \, \forall \phi \in X^{-\beta}.
\end{eqnarray}
Furthermore, $G_{u}(u, \epsilon)$ is linear and bounded uniformly on $u \in X^{\eta}$ and $0 \leq \epsilon \leq \epsilon_0$.
\end{lema}
\proof
 If  $u, w \in X^{\eta}$ and $\phi \in X^{-\beta}$, then
\begin{eqnarray} \label{deriv_Gateauxg}
& & \left| \frac{1}{t} \left\langle G(u + tw, \epsilon) - G(u, \epsilon) - t \frac{\partial G}{\partial u} (u, \epsilon)w, \phi \right\rangle_{\beta, - \beta} \right| \nonumber \\
& = & \left| \frac{1}{t} \int_{\partial \Omega} (g(\gamma(u + tw)) - g(\gamma(u)) - tg^{'}(\gamma(u)) \gamma(w)) \, \gamma(\phi) \, |J_{\partial \Omega} h_{\epsilon}| \, d \sigma(x) \right| \nonumber \\
& \leq & \frac{1}{|t|} \int_{\partial \Omega} |g(\gamma(u + tw)) - g(\gamma(u)) - tg^{'}(\gamma(u)) \gamma(w) |\,|\gamma(\phi)|\,|J_{\partial \Omega} h_{\epsilon}| \, d\sigma(x) \nonumber \\
& \leq & \overline{K} \frac{1}{|t|} \int_{\partial \Omega} |g(\gamma(u + tw)) - g(\gamma(u)) - tg^{'}(\gamma(u)) \gamma(w) |\, |\gamma(\phi)|\, d\sigma(x) \nonumber \\
& \leq & \overline{K} \frac{1}{|t|} \left( \int_{\partial \Omega} |g(\gamma(u + tw)) - g(\gamma(u)) - tg^{'}(\gamma(u)) \gamma(w) |^{2} d \sigma(x) \right)^{\frac{1}{2}} \parallel \gamma(\phi) \parallel_{L^{2}(\partial \Omega)} \nonumber \\
& \leq & \overline{K_{1}} \, \overline{K} \frac{1}{|t|} \left( \int_{\partial \Omega} |g(\gamma(u + tw)) - g(\gamma(u)) - tg^{'}(\gamma(u)) \gamma(w) |^{2} d \sigma(x) \right)^{\frac{1}{2}} \parallel \phi \parallel_{X^{-\beta}} \nonumber \\
& \leq & \overline{K_{1}} \, \overline{K} \left( \int_{\partial \Omega} |(g^{'}(\gamma(u + \bar{t} w)) - g^{'}(\gamma(u))) \gamma(w) |^{2} d \sigma(x) \right)^{\frac{1}{2}} \parallel \phi \parallel_{X^{-\beta}}
\end{eqnarray}
\noindent where  $ \overline{K} = \sup \{\, |J_{\partial \Omega}h_{\epsilon} (x)| \, | \, x \in \partial \Omega, 0 \leq \epsilon \leq \epsilon_{0} \}, \, \overline{K_{1}}$ is the constant given by the continuity of the trace map from  $X^{-\beta}$ into  $L^{2}(\partial\Omega)$ since $- \beta > \frac{1}{4}$  and we have used the Mean Value Theorem in the last passage, with  $0 < \bar{t} < t$. The integrand in \eqref{deriv_Gateauxg} is bounded by the integrable function $ 4\parallel g' \parallel_{\infty}^{2} |\gamma(w)|^{2}$, using \eqref{hip_derivg}, and goes to $0$ as   $t \to 0$.  Thus, from the Lebesgue Dominated Convergence Theorem
\begin{equation*}
\lim_{t \to 0} \frac{G(u + tw, \epsilon) - G(u, \epsilon)}{t} = \frac{\partial G}{\partial u}(u, \epsilon)w \,\, \text{in} \, X^{\beta}, \, \forall u, w \in X^{\eta}.
\end{equation*}

Therefore, $G$  is Gateaux - differentiable with Gateaux differential given by  (\ref{derg}), for  $0 < \epsilon \leq \epsilon_{0}$.
Now, if $u, w \in X^{\eta}$ and $\phi \in X^{-\beta}$,
\begin{eqnarray*}
\left| \left\langle \frac{\partial G}{\partial u}(u, \epsilon)w, \phi \right\rangle_{\beta, -\beta} \right| & \leq &  \int_{\partial \Omega} | g^{'}(\gamma(u))\gamma(w) \gamma(\phi)| \, |J_{\partial \Omega} h_{\epsilon}| \, d \sigma(x)  \\
& \leq & \overline{K} \parallel g^{'} \parallel_{\infty} \int_{\partial \Omega} |\gamma(w) \gamma(\phi)| \, d \sigma(x) \\
& \leq & \overline{K} \parallel g^{'} \parallel_{\infty} \, \parallel \gamma(w) \parallel_{L^{2}(\partial \Omega)} \, \parallel \gamma(\phi) \parallel_{L^{2}(\partial \Omega)} \\
& \leq & \overline{K_{1}} \, \overline{K_{3}} \, \overline{K} \parallel g^{'} \parallel_{\infty} \, \parallel w \parallel_{X^{\eta}} \, \parallel \phi \parallel_{X^{-\beta}}.
\end{eqnarray*}
Thus,
\begin{equation} \label{limitderivg}
 \left| \left| \frac{\partial G}{\partial u} (u,\epsilon) w \right| \right|_{X^{\beta}} \, \leq \overline{K_{1}} \, \overline{K_{3}} \, \overline{K} \parallel g^{'} \parallel_{\infty} \, \parallel w \parallel_{X^{\eta}},
\end{equation}
where  $\overline{K_{3}}$  is the constant from the continuity of the trace map from  $X^{\eta}$ into $L^{2}(\partial \Omega)$
since $\eta > \frac{1}{4}$.
It proves that $\frac{\partial G}{\partial u}(u, \epsilon)$ is 
 bounded as claimed. \eproof

\begin{lema} \label{contderg}
If $g$ satisfies  (\ref{hip_derivg}), $\beta + 1 > \eta > \frac{1}{4}$ and   $-\frac{1}{2} \leq \beta < -\frac{1}{4}$, then the Gateaux  derivative of  $G(u, \epsilon)$ with respect to  $u$ is continuous  in  $u$, that is , the map  $u \mapsto \frac{\partial G}{\partial u}(u, \epsilon) \in \mathcal{B}(X^{\eta}, X^{\beta})$ is continuous. Furthermore, $\epsilon \mapsto \frac{\partial G}{\partial u}(u, \epsilon)$ is continuous, uniformly for  $u  \in X^{\eta}$.
\end{lema}

\proof  Let  $u_{n} \to u \in X^{\eta}$ and  $\frac{1}{4} < \tilde{\eta} < \eta$. Then, if  $w \in X^{\tilde{\eta}}, \, \phi \in X^{-\beta}$,
\begin{eqnarray} \label{cont_gateauxg}
& & \left| \left\langle \left( \frac{\partial G }{\partial u}(u_{n}, \epsilon) - \frac{\partial G}{\partial u}(u, \epsilon) \right)w , \phi \right\rangle_{\beta, - \beta} \right| \nonumber \\
& \leq & \int_{\partial \Omega} |(g^{'}(\gamma(u)) - g^{'}(\gamma(u_{n}))) \, \gamma(w) \, \gamma(\phi) \, | \, |J_{\partial \Omega} h_{\epsilon}| \, d \sigma(x) \nonumber \\
& \leq & \overline{K} \left( \int_{\partial \Omega} |(g^{'}(\gamma(u)) - g^{'}(\gamma(u_{n}))) \, \gamma(w)|^{2} d \sigma(x) \right)^{\frac{1}{2}} \parallel \gamma(\phi) \parallel_{L^{2}(\partial \Omega)} \nonumber \\
& \leq & \overline{K} \, \overline{K_{1}} \parallel \phi \parallel_{X^{-\beta}} \left( \int_{\Omega} |(g^{'}(\gamma(u)) - g^{'}(\gamma(u_{n}))) \, \gamma(w)|^{2} d \sigma(x) \right)^{\frac{1}{2}}
\end{eqnarray}
\noindent where  $\overline{K} , \overline{K_{1}}$ are constants given by the boundness of $|J_{\partial \Omega}h_{\epsilon}|$ for  $0 \leq \epsilon \leq \epsilon_{0}$ and the continuity of the trace map  from  $X^{-\beta}$ into  $L^{2}(\partial \Omega)$. Observe that the integrand in \eqref{cont_gateauxg} is bounded by the integrable function   $4\parallel g^{'} \parallel^{2}_{\infty} |\gamma(w)|^{2}$
and goes to zero  as  $u_{n} \to  u \in X^{\eta}$. Therefore, by Lebesgue Dominated Convergence Theorem, the sequence  of operators  $\frac{\partial G}{\partial u}(u_{n}, \epsilon)$ converges strongly  to  $\frac{\partial G}{\partial u}(u, \epsilon)$ in  $B(X^{\tilde{\eta}}, X^{\beta})$. Since  $X^{\eta}$ is compactly embedded in   $X^{\tilde{\eta}}$,  it follows from Lemma 4.5 in \cite{marcone} that they converge also in the operator norm in $\mathcal{B}(X^{\eta}, X^{\beta})$.

With respect to the convergence in $\epsilon$, it is  sufficient, for our purposes,  to prove it for $\epsilon = 0$. We have
\begin{eqnarray*}
\left\langle \frac{\partial G}{\partial u}(u, \epsilon)w, \phi \right\rangle_{\beta, - \beta} & = & \int_{\partial \Omega} g^{'}(\gamma(u))\, \gamma(w) \, \gamma(\phi) \, |J_{\partial \Omega} h_{\epsilon}| \, d\sigma(x) \\
& = & \int_{I_{1}} g^{'}(\gamma(u))\, \gamma(w) \, \gamma(\phi) \, |J_{\partial \Omega} h_{\epsilon}| \, d\sigma(x) \nonumber \\
& + & \int_{I_{2}} g^{'}(\gamma(u))\, \gamma(w) \, \gamma(\phi) \, |J_{\partial \Omega} h_{\epsilon}| \, d\sigma(x) \\
& + & \int_{I_{3}} g^{'}(\gamma(u))\, \gamma(w) \, \gamma(\phi) \, |J_{\partial \Omega} h_{\epsilon}| \, d\sigma(x) \nonumber \\
& + & \int_{I_{4}} g^{'}(\gamma(u))\, \gamma(w) \, \gamma(\phi) \, |J_{\partial \Omega} h_{\epsilon}| \, d\sigma(x) \\
& = & \left\langle \left(\frac{\partial G_{1}}{\partial u} + \frac{\partial G_{2}}{\partial u} + \frac{\partial G_{3}}{\partial u} + \frac{\partial G_{4}}{\partial u} \right) (u, \epsilon)w, \phi \right\rangle_{\beta, - \beta},
\end{eqnarray*}

\noindent with  $G_{i}$ in the portion $I_{i}$, $1 \leq i \leq 4 $, of the  boundary.
The expression of $J_{\partial \Omega} h_{\epsilon} $ in each portion of the boundary is given in  \eqref{expressionJbound}.
%
%

 We first show that
\begin{equation*}
\left| \left| \left( \frac{\partial G_{i}}{\partial u}(u,\epsilon) - \frac{\partial G_{i}}{\partial u}(u,0) \right) \right| \right|_{\mathcal{L}(X^{\eta}, X^{\beta})} \to 0, \,\, \text{as} \,\, \epsilon \to 0,
\end{equation*}
for $i=2,3,4 $,   uniformly for  $u$ in bounded sets of  $X^\eta$.
This is immediate for $i=3, 4$, since there is no dependence in $\epsilon$. In $I_{2}$, we have
\begin{align*}
 & \left| \left\langle \left( \frac{\partial G_{2}}{\partial u}(u,\epsilon) - \frac{\partial G_{2}}{\partial u}(u,0) \right)w, \phi \right\rangle_{\beta, - \beta} \right|  \\  &  \leq  \int_{0}^{1}  |g^{'}(\gamma(u(1, x_{2}))) \gamma(w) \gamma(\phi)|(|Jh_{\epsilon}| - 1) dx_{2} \\
& \leq  (|Jh_{\epsilon}| - 1) \int_{0}^{1}  |g^{'}(\gamma(u(1, x_{2}))) \gamma(w) \gamma(\phi)| dx_{2} \\
& \leq  (|Jh_{\epsilon}| - 1) \parallel g' \parallel_{\infty} \parallel \gamma(w) \parallel_{L^{2}(\partial \Omega)} \parallel \gamma(\phi) \parallel_{L^{2}(\partial \Omega)} \\
& \leq  (|Jh_{\epsilon}| - 1) \overline{K}_{1} \overline{K}_{3}\parallel g' \parallel_{\infty} \parallel w \parallel_{X^{\eta}} \parallel \phi \parallel_{X^{-\beta}},
\end{align*}
 where  $\overline{K}_{1}, \overline{K}_{3}$ are constants due to the continuity of the trace function from   $X^{-\beta}$ into  $L^{2}(\partial \Omega)$ and   $X^{\eta}$ into  $L^{2}(\partial \Omega)$, respectively. Since  $Jh_{\epsilon} \to  1$ uniformly
 when $\epsilon \to 0 $, the convergence in this portion is proved.

  We now consider  the term
 in   $I_{1} = \{(x_{1}, 1) : 0 \leq x_{1} \leq 1 \}$, the oscillating portion of the boundary.  Let  $\frac{1}{4} < \tilde{\eta}  < \eta$ and $ -\frac{1}{2} \leq \beta < \beta' < \frac{-1}{4} $.  Then, if $w \in X^{\tilde{\eta}}$, $ \phi \in
 X^{\beta'} $,

\begin{align} \label{boundG1}
 & \left|
\left\langle \left( \frac{\partial G_{1}}{\partial u} (u, \epsilon)  - \frac{\partial G_{1}}{\partial u} (u, 0 )  \right) w  , \phi \right\rangle_{\beta', -\beta'} \right|  \nonumber \\
  & \leq  | \int_{0}^{1}  g^{'}(\gamma(u(x_{1}, 1))) \gamma(w) \gamma(\phi) \sqrt{1 + cos^{2}(x_{1}/ \epsilon)}  \, dx  \nonumber \\
   &  -
  \int_{0}^{1}  g^{'}(\gamma(u(x_{1}, 1))) \gamma(w) \gamma(\phi)
 M_{\pi} (p)  \, dx |
   \nonumber    \\
 & \leq   \|g'\|_{\infty} \left|
\left\langle \left(  T_ \epsilon - T_0    \right) w , \phi \ \right\rangle_{\beta', -\beta'} \right|,
 \end{align}
 where

\begin{align*}
M_{\pi} (p) &= \frac{1}{\pi} \int_{0}^{\pi} \sqrt{1 + cos^{2}(y)} \, dy, \\
 \left\langle   T_ \epsilon  w , \phi \right\rangle_{\beta', -\beta'}  & =   \int_{0}^{1}  \gamma(w) \gamma(\phi) \sqrt{1 + cos^{2}(x_{1}/ \epsilon)}  \, dx,  \\
  \left\langle   T_0  w , \phi \right\rangle_{\beta', -\beta'}  & =
\int_{0}^{1}  \gamma(w) \gamma(\phi)
M_{\pi} (p)  \, dx.
 \end{align*}

\indent The function  $p(x_{1}) = \sqrt{1 + cos^{2}(x_{1})}$ is periodic of period  $\pi$. From the Theorem of Convergence in the Average   (see Theorem 2.6 in  \cite{cionarescu} ), $p_{\epsilon} (x_{1}) = p(x_{1}/\epsilon)$ satisfies  $p_{\epsilon} \rightharpoonup M_{\pi} (p)$ weak* in  $L^{\infty}(0,1)$.

Therefore,

 \begin{equation} \label{convT}  \langle  T_ \epsilon  w , \phi \rangle_{\beta', -\beta'}  \to
  \langle  T_0  w , \phi \rangle_{\beta', -\beta'} \ ,
 \end{equation}
   and thus

  \begin{equation*}
\left\langle \frac{\partial G_{1}}{\partial u} (u, \epsilon) w, \phi \right\rangle_{\beta', -\beta'}  \to \left\langle \frac{\partial G_{1}} {\partial u} (u, 0)w, \phi \right\rangle_{\beta', -\beta'},
\end{equation*}
uniformly for $u \in  X^{\eta}$.

Now, by Lemma \ref{gderivavel} (with $g\equiv 1$),  $ T_{\epsilon} $ is bounded  as a map from $X^{\tilde{\eta}}$ to $X^{ \beta'} $,  uniformly in $\epsilon$.  Since the embedding $X^{\beta^{'}} \subset X^{\beta}$ is compact, it follows that, for any sequence
 $\epsilon_n \to 0$ and $w \in X^{\tilde{\eta}}$, there is a subsequence, which we still denote by $\epsilon_n $, and $z \in X^{\beta} $  such that
 $$ T_{\epsilon_n} w \to  z \text{ in  }  X^{\beta}.$$

 By \eqref{convT},
  $ T_{\epsilon} w \to  T_{0} w \text{  weakly in  }  X^{\beta'}, \text{ as }  \epsilon \to 0 , $
  and we conclude that

 $$ T_{\epsilon} w \to  T_{0} w \text{ in  }  X^{\beta}, \text{ as }  \epsilon \to 0 , $$
 for any $w \in X^{\tilde{\eta}}$,
 that is, $ T_{\epsilon}  \to T_0$ strongly as maps from
   $X^{\tilde{\eta}}$  to $X^{\beta}$.

  Since the embedding  $ X^{\eta} \subset
 X^{\tilde{\eta}} $ is compact, it follows from Lemma 4.5 in \cite{marcone} that $ T_{\epsilon} \to  T_{0}  \text{ in  }  \mathcal{L}(X^{\eta}, X^{\beta}) $.   Using \eqref{boundG1}, we conclude that

 \begin{equation*}
\frac{\partial G_{1}}{\partial u} (u, \epsilon)  \to
 \frac{\partial G_{1}} {\partial u} (u, 0)  \text{ in  }  \mathcal{L}(X^{\eta}, X^{\beta} )
\end{equation*}
uniformly for $u \in  X^{\eta}$. \eproof

\begin{lema} \label{derh} Suppose  $f$ and $g$ satisfy the hypotheses \eqref{hip_derivf} and \eqref{hip_derivg}, respectively, $\beta + 1 > \eta > \frac{1}{4}$ and $-\frac{1}{2} \leq \beta < - \frac{1}{4}$. Then, the map   $(H_{\epsilon})_{\beta} = (F_{\epsilon})_{\beta} + (G_{\epsilon})_{\beta} : X^{\eta} \rightarrow X^{\beta}$ is continuously Fr\'echet differentiable with respect to  $u$  and  its derivative $((H_{\epsilon})_{\beta})_{u}$ is bounded in $u$, uniformly for  $0\leq \epsilon \leq \epsilon_0$ and  continuous in $\epsilon$, at  $\epsilon = 0$, uniformly for  $u \in X^{\eta}$.
\end{lema}

\proof
By Lemmas \ref{contderf} and \ref{contderg},  the map $H_{\epsilon}$ is continuously Gateaux differentiable, so the claim follows from well known results (see, for instance, Proposition 2.8 in \cite{rall}). The boundness of $((H_{\epsilon})_{\beta})_{u}$ follows from the inequalities    (\ref{limitderivf}) and  (\ref{limitderivg}) and its continuity with respect to $\epsilon$ is proved in the mentioned lemmas. \eproof

\indent We are now in a position to prove the continuity of equilibria of \eqref{problem_abstract_scale}.

\begin{teo}\label{contequi} Suppose that $1/4 < \eta < 1/2$,   $f$ and $g$
satisfy the conditions for global existence, Theorem \ref{exist_global}, and the hypotheses \eqref{hip_derivf}, \eqref{hip_derivg}  and the equilibria of \ref{problem_abstract_scale} are all hyperbolic. Then, the family of equilibria    $\{ E_{\epsilon} \, | \, 0 \, \leq \,\epsilon \, \leq \, \epsilon_{0} \}$ of \eqref{problem_abstract_scale} is lower semicontinuous in  $X^{\eta}$, at $\epsilon = 0$.
\end{teo}

\proof
Observe that  $u \in X^{\eta}$ is an equilibrium of \eqref{problem_abstract_scale} if and only if
$
(A_{\epsilon})_{\beta}u = (H_{\epsilon})_{\beta}u  \Longleftrightarrow  -(A_{\epsilon})_{\beta}u + (H_{\epsilon})_{\beta}u = 0
 \Longleftrightarrow  -(A_{\epsilon})_{\beta}^{-1} (A_{\epsilon})_{\beta}u + (A_{\epsilon})_{\beta}^{-1} (H_{\epsilon})_{\beta}u = 0
 \Longleftrightarrow  -u + (A_{\epsilon})_{\beta}^{-1} (H_{\epsilon})_{\beta}u = 0.
$
 We then wish to   apply the Implicit Function Theorem to the function  $Z = - I + (A_{\epsilon})_{\beta}^{-1} (H_{\epsilon})_{\beta} : X^{\eta} \times \mathbb{R} \rightarrow X^{\eta}$ to show that near any equilibrium of \eqref{problem_abstract_scale}, with $\epsilon =0$, there exists an equilibria of \eqref{problem_abstract_scale} for $\epsilon$ sufficiently small  .
 To this end, we need to check that
  $Z $ and its partial derivative $\frac{\partial Z}{\partial u}$ are both continuous in $u$ and  $\epsilon$ and $\frac{\partial Z}{\partial u}(e,0)$ is invertible, for any equilibrium $ e \in X^{\eta}$.
(see \cite{loomis}, Theorem 9.3). The continuity of $Z$ in both variables  follows from Corollary \ref{limitresolunif}, Theorem \ref{convresolnorm}, Theorem \ref{propfg} and  Lemma \ref{contfgepsilon}. The continuity of $\frac{\partial Z}{\partial u}$ follows from Lemma \ref{derh}.

Now, we have $ \frac{\partial Z}{\partial u} = \frac{\partial (I + (A_{\epsilon})_{\beta}^{-1}(H_{\epsilon})_{\beta})}{\partial u}$ and
\begin{eqnarray*}
\frac{\partial ((A_{\epsilon})_{\beta}^{-1}(-(A_{\epsilon})_{\beta} + (H_{\epsilon})_{\beta}))}{\partial u} & = & \frac{\partial (A_{\epsilon})_{\beta}^{-1}}{\partial u} (-(A_{\epsilon})_{\beta} + (H_{\epsilon})_{\beta})\frac{\partial (-(A_{\epsilon})_{\beta} + (H_{\epsilon})_{\beta})}{\partial u} \\
& = & (A_{\epsilon})_{\beta}^{-1}(-(A_{\epsilon})_{\beta} + (H_{\epsilon})_{\beta})\frac{\partial (-(A_{\epsilon})_{\beta} + (H_{\epsilon})_{\beta})}{\partial u}.
\end{eqnarray*}

\noindent  Since the equilibria are all hyperbolic, that is,  $\frac{\partial (-(A_{\epsilon})_{\beta} + (H_{\epsilon})_{\beta})}{\partial u} (e,0)$  is invertible, it follows that $ \frac{\partial Z}{\partial u}(e,0)$ is also invertible. \eproof

\section{Lower semicontinuity of the attractors}

 In this section, we finally prove the lower semicontinuity of the attractors for the family of flows generated by \eqref{problem_abstract_scale}.

\begin{lema} \label{propr} Suppose that  $f$ and $g$ satisfy the hypotheses of Theorem  \ref{contequi} and $e(\epsilon)$ is a hyperbolic equilibrium  of \eqref{problem_abstract_scale}, for $0 \leq \epsilon \leq \epsilon_0$. Then,  the map  $(H_\epsilon)_{\beta}  = (F_{\epsilon})_{\beta} + (G_{\epsilon})_{\beta} : X^{\eta} \to  X^{\beta}$, $\frac{1}{2} > \eta > \frac{1}{4}$, $-\frac{1}{2} \leq \beta < - \frac{1}{4}$,
satisfies
$ (H_{\epsilon})_{\beta} (e(\epsilon) + w, \epsilon) = (A_{\epsilon})_{\beta} e(\epsilon) + ((H_{\epsilon})_{\beta})_{u}(e(\epsilon), \epsilon)w + r(w, \epsilon),  $  with  $r(w, \epsilon)= o( \|w\|_{{X}^{\beta}})$, and
\begin{itemize}
\item[i)] $ \parallel r(w,0) - r(w, \epsilon) \parallel_{X^{\beta}} \, \to 0$ as $\epsilon \to 0$,  uniformly for $w$ in bounded sets,
\item[ii)] $\parallel r(w_{1}, \epsilon) - r(w_{2}, \epsilon) \parallel_{X^\beta} \, \leq k(w_1,w_2) \parallel w_{1} - w_{2} \parallel_{X^{\eta}}, \,\,$  with $ k(w_1,w_2)  \to 0$ as $ \parallel w_{1}  \parallel_{X^{\eta}}, \parallel w_{2} \parallel_{X^{\eta}} \to  0 \,$ uniformly for  $0 \leq \epsilon \leq \epsilon_0.$
\end{itemize}

\end{lema}

\proof Let $\epsilon_0$ be such that $e(\epsilon)$ is a hyperbolic equilibrium of \eqref{problem_abstract_scale}, for $\epsilon \leq \epsilon_0.$   Since
 $(H_{\epsilon})_{\beta}$ is  Frechet differentiable,  we have
 $(H_{\epsilon})_{\beta} (e(\epsilon) + w, \epsilon)  =  (H_{\epsilon})_{\beta}(e(\epsilon), \epsilon) + ((H_{\epsilon})_{\beta})_{u}(e(\epsilon), \epsilon)w + r(w, \epsilon)$ and, since $e(\epsilon)$ is an equilibrium,
 $(H_{\epsilon})_{\beta}(e(\epsilon), \epsilon) = (A_{\epsilon})_{\beta}e(\epsilon)$.
Now, if  $w \in X^{\eta}$,
\begin{eqnarray}
 \parallel r(w,0) - r(w, \epsilon) \parallel_{X^{\beta}} & = & \parallel (H)_{\beta} (e(0) + w, 0) - (H)_{\beta}(e(0), 0) - ((H)_{\beta})_{u}(e(0), 0)w  \nonumber \\
& - & (H_{\epsilon})_{\beta} (e(\epsilon) + w, \epsilon) + (H_{\epsilon})_{\beta}(e(\epsilon), \epsilon) + ((H_{\epsilon})_{\beta})_{u}(e(\epsilon), \epsilon)w \parallel_{X^\beta} \nonumber \\
& \leq & \parallel (H_{\epsilon})_{\beta} (e(0) + w, 0) - (H_{\epsilon})_{\beta}(e(0), 0) \nonumber \\
& - & (H_{\epsilon})_{\beta} (e(\epsilon) + w, \epsilon) + (H_{\epsilon})_{\beta}(e(\epsilon), \epsilon) \parallel_{X^{\beta}}
\nonumber  \\
& + & \parallel ((H_{\epsilon})_{\beta})_{u}(e(0), 0)w - ((H_{\epsilon})_{\beta})_{u}(e(\epsilon), \epsilon)w \parallel_{X^\beta}. \nonumber
\end{eqnarray}

 The item i) follows then from the continuity of
 $ (H_{\epsilon})_{\beta}$  and  $((H_{\epsilon})_{\beta})_u$  in $X^{\eta}$ given by Theorem \ref{propfg} and Lemma \ref{derh}, respectively, and in $\epsilon$, by Lemma 9.6, and the continuity of the equilibria, given by Theorems \ref{contsupeq} and  \ref{contequi}.

ii) If  $w_{1}, w_{2} \in X^{\eta}, \, 0 \leq \epsilon \leq \epsilon_{0}$,

\begin{align}
\parallel r(w_{1}, \epsilon) - r(w_{2}, \epsilon) \parallel_{X^\beta} & =  \parallel (H_{\epsilon})_{\beta} (e(\epsilon) + w_{1}, \epsilon) - (H_{\epsilon})_{\beta}(e(\epsilon), \epsilon) - ((H_{\epsilon})_{\beta})_{u}(e(\epsilon), \epsilon)w_{1} \nonumber \\
& -  (H_{\epsilon})_{\beta} (e(\epsilon) + w_{2}, \epsilon) + (H_{\epsilon})_{\beta}(e(\epsilon), \epsilon) + ((H_{\epsilon})_{\beta})_{u}(e(\epsilon), \epsilon)w_{2} \parallel_{X^\beta} \nonumber \\
& \leq  \parallel F(e(\epsilon) + w_{1}, \epsilon) - F(e(\epsilon), \epsilon) - F_{u}(e(\epsilon), \epsilon)w_{1} \nonumber \\
& -  F(e(\epsilon) + w_{2}, \epsilon) + F(e(\epsilon), \epsilon) + F_{u}(e(\epsilon), \epsilon)w_{2} \parallel_{X^\beta} \label{estimateF1}\\
& +  \parallel G(e(\epsilon) + w_{1}, \epsilon) - G(e(\epsilon), \epsilon) - G_{u}(e(\epsilon), \epsilon)w_{1} \nonumber  \\
& -  G(e(\epsilon) + w_{2}, \epsilon) + G(e(\epsilon), \epsilon) + G_{u}(e(\epsilon), \epsilon)w_{2} \parallel_{X^\beta}. \label{estimateG1}
\end{align}

We estimate \eqref{estimateF1}:

\begin{align*}
 & |  < F(e(\epsilon) + w_{1}, \epsilon) - F(e(\epsilon), \epsilon) - F_{u}(e(\epsilon), \epsilon)w_{1}  \\
  & -  F(e(\epsilon) + w_{2}, \epsilon) + F(e(\epsilon), \epsilon) + F_{u}(e(\epsilon), \epsilon)w_{2}, \phi >_{\beta, - \beta}| \\
 &\leq  \int_{\Omega} | (f(e(\epsilon) + w_{1}) - f(e(\epsilon)) - f'(e(\epsilon))w_{1} -  f(e(\epsilon) + w_{2}) + f(e(\epsilon))  \\ & \quad \quad \quad + f'(e(\epsilon))w_{2}) \, \phi |  \, |Jh_{\epsilon}| \, dx \\
& \leq  K \left( \int_{\Omega} |(f'(e(\epsilon) + \xi_{x}) - f'(e(\epsilon))) (w_{1} - w_{2})|^{2} \, dx \right)^{\frac{1}{2}} \parallel \phi \parallel_{L^{2}(\Omega)} \\
& \leq  K K_{1} \left( \int_{\Omega} |(f'(e(\epsilon) + \xi_{x}) - f'(e(\epsilon))) (w_{1} - w_{2})|^{2} \, dx \right)^{\frac{1}{2}} \parallel \phi \parallel_{X^{-\beta}} \\
& \leq K K_{1} \left( \int_{\Omega} (f'(e(\epsilon) + \xi_{x}) - f'(e(\epsilon)))^{2p} \, dx \right)^{\frac{1}{2p}} \, \left( \int_{\Omega} (w_{1} - w_{2})^{2q} dx \right)^{\frac{1}{2q}} \parallel \phi \parallel_{X^{-\beta}} \\
& \leq K K_{1} \left( \int_{\Omega} (f'(e(\epsilon) + \xi_{x}) - f'(e(\epsilon)))^{2p}\, dx \right)^{\frac{1}{2p}} \parallel w_{1} - w_{2} \parallel_{L^{2q}(\Omega)}  \parallel \phi \parallel_{X^{-\beta}} \\
& \leq  K K_{1} K_{6} \left( \int_{\Omega} (f'(e(\epsilon) + \xi_{x}) - f'(e(\epsilon)))^{2p}\, dx \right)^{\frac{1}{2p}} \parallel w_{1} - w_{2} \parallel_{X^{\eta}}  \parallel \phi \parallel_{X^{-\beta}} \\
\end{align*}
where $K = \sup \{|Jh_{\epsilon}(x)| \, | \, x \in \Omega, \, 0 \leq \epsilon \leq \epsilon_{0} \}$,  $\xi_{x}$ comes from the application of the Mean Value Theorem, for each  $x \in \Omega, \, w_{1}(x) \leq \xi_{x} \leq w_{2}(x) \,\, \text{or} \,\, w_{2}(x) \leq \xi_{x} \leq w_{1}(x)$, 
$K_{1}$ is a constant for  the  embedding of $X^{-\beta}$ into $L^{2}(\Omega)$ and  $K_{6}$ is a  constant for the embedding of $X^{\eta}$ into  $L^{2q}$ with  $q < \frac{1}{1 - 2\eta}$, $p$ and $q$ are conjugate exponents.
%

Now, the integrand above is bounded by   $ 4^{p}\parallel f' \parallel_{\infty}^{2p}$, due to  (\ref{hip_derivf}) and goes to  $0$ a.e.  as  $\parallel w_{1} \parallel_{X^{\eta}} $,
$\parallel w_{2} \parallel_{X^{\eta}} \to 0$.
Thus, the integral converges to  $0$  by  Lebesgue  Dominated Convergence Theorem.

 As for \eqref{estimateG1}, we have
\begin{align*}
&  | < G(e(\epsilon) + w_{1}, \epsilon) - G(e(\epsilon), \epsilon) - G_{u}(e(\epsilon), \epsilon)w_{1} \\
& -  G(e(\epsilon) + w_{2}, \epsilon) + G(e(\epsilon), \epsilon) + G_{u}(e(\epsilon), \epsilon)w_{2}, \phi >_{\beta, - \beta} | \\
& \leq  \int_{\partial \Omega} | (g(\gamma(e(\epsilon) + w_{1})) - g(\gamma(e(\epsilon))) - g'(\gamma(e(\epsilon)))\gamma(w_{1}) \\
&  \quad \quad \quad -   g(\gamma(e(\epsilon) + w_{2})) + g(\gamma(e(\epsilon))) + g'(\gamma(e(\epsilon)))\gamma(w_{2})) \, \gamma(\phi)| \, |J_{\partial \Omega} h_{\epsilon}| \,  d\sigma(x) \\
& \leq  \int_{\partial \Omega} | (g'(\gamma(e(\epsilon) + \xi_{x})) - g'(\gamma(e(\epsilon)))) \gamma(w_{1} - w_{2}) \gamma(\phi)| \, |J_{\partial \Omega} h_{\epsilon}| \,  d\sigma(x) \\
& \leq  \overline{K} \left( \int_{\partial \Omega} |(g'(\gamma(e(\epsilon) + \xi_{x})) - g'(\gamma(e(\epsilon)))) \gamma(w_{1} - w_{2})|^{2}  d\sigma(x) \right)^{\frac{1}{2}} \parallel \gamma(\phi) \parallel_{L^{2} (\partial \Omega)} \\
& \leq  \overline{K} \, \overline{K_{1}} \left( \int_{\partial \Omega} |(g'(\gamma(e(\epsilon) + \xi_{x})) - g'(\gamma(e(\epsilon)))) \gamma(w_{1} - w_{2})|^{2}  d\sigma(x) \right)^{\frac{1}{2}} \parallel \phi \parallel_{X^{-\beta}} \\
& \leq  \overline{K} \, \overline{K_{1}} \left( \int_{\partial \Omega} (g'(\gamma(e(\epsilon) + \xi_{x})) - g'(\gamma(e(\epsilon))))^{2\bar{p}}  d\sigma(x) \right)^{\frac{1}{2\bar{p}}} \,  \\
 & \left( \int_{\partial \Omega} \gamma(w_{1} - w_{2})^{2\bar{q}} dx \right)^{\frac{1}{2\bar{q}}} \parallel \phi \parallel_{X^{-\beta}} \\
& \leq  \overline{K} \, \overline{K_{1}} \left( \int_{\partial \Omega} (g'(\gamma(e(\epsilon) + \xi_{x})) - g'(\gamma(e(\epsilon))))^{2\bar{p}} dx \right)^{\frac{1}{2\bar{p}}} \parallel \gamma(w_{1} - w_{2}) \parallel_{L^{2\bar{q}}(\partial \Omega)}  \parallel \phi \parallel_{X^{-\beta}} \\
& \leq  \overline{K} \, \overline{K_{1}} \, \overline{K_{6}} \left( \int_{\partial \Omega} (g'(\gamma(e(\epsilon) + \xi_{x})) - g'(\gamma(e(\epsilon))))^{2\bar{p}} dx \right)^{\frac{1}{2\bar{p}}} \parallel w_{1} - w_{2} \parallel_{X^{\eta}}  \parallel \phi \parallel_{X^{-\beta}} \\
\end{align*}
where  $\overline{K} = \sup \{ \, |J_{\partial \Omega} h_{\epsilon}(x)|, \, x \in \partial \Omega, \, 0 \leq \epsilon \leq \epsilon_{0} \}$, $\overline{K}_{1}$ is due to the continuity of the trace map $\gamma$ from $X^{-\beta}$ into $L^{2}(\partial \Omega)$, when  $-\beta > \frac{1}{4}$, $\overline{K_{6}}$ is due to the continuity of the trace map from  $X^{\eta} \subset H^{2\eta}(\Omega)$ into  $L^{2\bar{q}}(\partial \Omega)$, when $\bar{q} < \frac{1}{2 - 4\eta}$, $\xi_{x}$ comes from the application of the Mean Values Theorem, for each $x \in \Omega, \, w_{1}(x) \leq \xi_{x} \leq w_{2}(x) \,\, \text{or} \,\, w_{2}(x) \leq \xi_{x} \leq w_{1}(x)$,  $\bar{p}$ and  $\bar{q}$ are conjugated exponents. 

%
%

Now, the integrand above is bounded by   $ 4^{\bar{p}}\parallel g' \parallel_{\infty}^{2\bar{p}}$, due to  (\ref{hip_derivg}) and goes to  $0$ a.e.  as  $\parallel w_{1} \parallel_{X^{\eta}} $,
$\parallel w_{2} \parallel_{X^{\eta}} \to 0$.
Thus, the integral converges to  $0$  by the Lebesgue  Dominated Convergence Theorem. This completes the proof of item ii). \eproof

 Let  $\epsilon_0 $ be such that  $e(\epsilon)$  is a hyperbolic equilibrium of \eqref{problem_abstract_scale}  for  $0 \leq \epsilon \leq \epsilon_0 $. Then, the linearized operator  $L(\epsilon) = (A_{\epsilon})_{\beta} - ((H_{\epsilon})_{\beta})_{u} (e(\epsilon), \epsilon)$ is an isomorphism. We now state some of its additional properties, starting with an auxiliary result.

 \begin{lema} \label{HAbounded} If  $f, g, \beta, \eta$ satisfy the conditions of Lemma \ref{derh} then $((H_{\epsilon})_{\beta})_{u} (u, \epsilon)$ is
 $(A_{\epsilon})_{\beta}$ bounded, that is,

\begin{equation} \label{Abounded}
\parallel ((H_{\epsilon})_{\beta})_{u} (u, \epsilon)w \parallel_{X^{\beta}} \leq C \zeta \parallel (A_{\epsilon})_{\beta} w \parallel_{X^{\beta}} + CC' \zeta^{\frac{-(\eta - \beta)}{1 - (\eta - \beta)}} \parallel w \parallel_{X^{\beta}}, \,\,\, \forall \zeta > 0,
\end{equation}\\
\noindent for any $w \in D((A_{\epsilon})_{\beta})$, where  $C$ and $C'$ are constants independent of $\epsilon$.
\end{lema}

\proof   By Lemma \ref{derh},  $((H_{\epsilon})_{\beta})_{u} (u, \epsilon) : X^{\eta} \rightarrow X^{\beta}$ is bounded, uniformly in $\epsilon$.  By  Theorem 1.4.4 in \cite{henry_semilinear}, we have

\begin{eqnarray*}
\parallel ((H_{\epsilon})_{\beta})_{u} (u, \epsilon)w \parallel_{X^{\beta}} & \leq & C \parallel w \parallel_{X^{\eta}} \\
& = & C \parallel (A_{\epsilon})_{\beta}^{\eta - \beta} w \parallel_{X^{\beta}} \\
& \leq & C \zeta \parallel (A_{\epsilon})_{\beta} w \parallel_{X^{\beta}} + C C' \zeta^{\frac{-(\eta - \beta)}{1 - (\eta - \beta)}} \parallel w \parallel_{X^{\beta}}
\end{eqnarray*}

\noindent  where $C'$ depends only on the sector and on the constant in the resolvent inequality for $(A_{\epsilon})_{\beta}$, and   $\zeta > 0$  is an arbitrarily small constant.  \eproof

\begin{teo} \label{linearizadoset}
 The linearized operator  $L(\epsilon) = (A_{\epsilon})_{\beta} - ((H_{\epsilon})_{\beta})_{u} (e(\epsilon), \epsilon) : X^{\beta + 1} \subset X^\beta \rightarrow X^{\beta}$, $-\frac{1}{2} \leq \beta < -\frac{1}{4}$, $0 \leq \epsilon \leq \epsilon_{0}$, is self-adjoint, has compact resolvent  and  there is a common sector and a common constant in the resolvent inequality for the family $\{ L(\epsilon) \}_{0 \,\leq \, \epsilon \, \leq \, \epsilon_{0}}$.
\end{teo}

\proof The operator  $(A_{\epsilon})_{\beta} : X^{\beta + 1} \subset X^{\beta} \to X^{\beta}$, $-\frac{1}{2} \leq \beta \leq 0$, is self-adjoint and $ \|(A_{\epsilon})_{\beta}u\|_{\beta}
 \geq b \| u \|_{\beta+1} $, with $ b> 0$, by Lemma \ref{self_positive}  and the fact, already observed, that  $(A_{\epsilon})_{\beta+ 1/2} : X^{\beta}
 (X^{\beta+1} ) \to
  X^{-1/2} ( X^{+1/2} ) $  is an isometry.
   On the other hand,     $((H_{\epsilon})_{\beta})_{u} (u, \epsilon) : X^{\eta} \rightarrow X^{\beta}$, $\beta + 1 \geq \eta > \frac{1}{4}$, is symmetric and $ (A_{\epsilon})_{\beta}$  bounded, uniformly in $\epsilon$, by Lemma \ref{HAbounded}.   It follows that
    \begin{equation} \label{lowerboundL}
     \|L_{\epsilon} + \rho u\|_{\beta}
 \geq  b' \| u \|_{\beta+1},
\end{equation}
 for some constants $\rho$ and $b'$.

Therefore, $ L_{\epsilon} $   is (uniformly) lower bounded and   self-adjoint and, thus,  sectorial, with  common sector and constant in the resolvent inequality.

 The compactness of the resolvent follows from \eqref{lowerboundL}   and the compactness of the embedding of $ X^{\beta+1} \hookrightarrow  X^{\beta}$. \eproof

\begin{teo}\label{convreslinea}  If  $\lambda \in \cap \, \{ \rho(L(\epsilon)) : 0 \leq \epsilon \leq \epsilon_{0} \}$, $\frac{1}{2} > \eta > \frac{1}{4}$, $-\frac{1}{2} < \beta < -\frac{1}{4}$, then  $ (\lambda I - L(\epsilon))^{-1}  \to  (\lambda I - L(0))^{-1}$  in ${ \mathcal{L}( X^{\beta}, X^{\eta})}$  as  $\epsilon \to 0 $.
\end{teo}

\proof We write
\begin{equation} \label{formresol}
(\lambda I - L(\epsilon))^{-1} = (\lambda I - (A_{\epsilon})_{\beta})^{-1} ( I - (-((H_{\epsilon})_{\beta})_{u} (e(\epsilon), \epsilon))(\lambda I - (A_{\epsilon})_{\beta})^{-1})^{-1}.
\end{equation}

\noindent  Since
 $ (\lambda I - (A_{\epsilon})_{\beta})^{-1}  \to
 (\lambda I - A_{\beta})^{-1} $ in  ${ \mathcal{L}( X^{\beta}, X^{\eta})}$  by Remark \ref{convresolnormfrac}, it remains to prove that
$  ( I - (-((H_{\epsilon})_{\beta})_{u} (e(\epsilon), \epsilon))(\lambda I - (A_{\epsilon})_{\beta})^{-1})^{-1} \to
 ( I - (-((H_{\epsilon})_{\beta})_{u} (e(0), 0))(\lambda I - (A_{\beta}))^{-1})^{-1}$ in  ${ \mathcal{L}( X^{\beta}, X^{\beta})}$.  By a well known result, it is  sufficient to prove that the operators  $  ( I - (-((H_{\epsilon})_{\beta})_{u} (e(\epsilon), \epsilon))(\lambda I - (A_{\epsilon})_{\beta})^{-1}) $ are bounded, converge to
 $( I - (-((H_{\epsilon})_{\beta})_{u} (e(0), 0))(\lambda I - (A_{\beta}))^{-1})$ in  ${ \mathcal{L}( X^{\beta}, X^{\beta})}$ and their inverses are also bounded.

  The boundeness follows from Lemma \ref{derh} and  Corollary \ref{limitresolunif}, and  the convergence, from   Remark \ref{convresolnormfrac}  and Lemma \ref{derh}. \eproof

\begin{lema}
Let  $X,Y$ be Banach spaces,  $T(\epsilon): X \to Y $, $0\leq \epsilon \leq \epsilon_0 $,  a family
of operators in $\mathcal{L}(X,Y)$ and    $R_{\epsilon}(\lambda)$   denote their resolvent operators.
 Then we have the following identity, for $\lambda$ and
   $\lambda_0$  in their resolvent set,

\begin{equation} \label{identresolv}
R_{\epsilon}(\lambda) - R_0(\lambda) = (1 + (\lambda - \lambda_{0})R_{\epsilon}(\lambda)) (R_{\epsilon}(\lambda_{0}) - R_0 (\lambda_{0}))(1 + (\lambda - \lambda_{0})R_{0}(\lambda))
\end{equation}\\
\proof The proof is obtained by applying the first resolvent identity

\begin{equation*}
R_{\epsilon}(\lambda_{1}) - R_{\epsilon}(\lambda_{2}) = (\lambda_{1} - \lambda_{2}) R_{\epsilon}(\lambda_{1})R_{\epsilon}(\lambda_{2})
\end{equation*}
to each operator.  \eproof
\end{lema}

 \begin{cor} \label{convreslinea_compact}
 The convergence obtained in \ref{convreslinea} is uniform for
 $\lambda $ in a compact set
 of $ \cap \, \{ \rho(L(\epsilon)) : 0 \leq \epsilon \leq \epsilon_{0} \}$.
 \end{cor}
 \proof  The result is an easy consequence of \eqref{identresolv}
  and can also be obtained by checking that the boundness and convergences in  the proof of Theorem
   \ref{convreslinea} are uniform in compact sets of the resolvent.
   \eproof

From Theorem \ref{linearizadoset}, the spectral set  $\sigma(L(\epsilon))$ is a discrete set of real eigenvalues of finite multiplicity, with $ + \infty$ as its only accumulation point.  Also, by Theorem \ref{convreslinea}, $\sigma(L(\epsilon))$ is continuous at $\epsilon=0$ (see \cite{teschl}, theorem 6.38).

 Consider the decomposition
 $X^{\beta} =  X_{1}(\epsilon) +  X_{2}(\epsilon)$  corresponding to the spectral sets $\sigma_{1}= \sigma(L(\epsilon))\cap \{ \Re(\lambda) < 0 \}$ and  $ \sigma_{2} = \sigma ((L(\epsilon)) \cap \{ \Re(\lambda) > 0 \}$. Let  $E_{1}= E_{1}(0), \, E_{2} = E_{2}(0)$ be the corresponding spectral projections, that is

\begin{equation*}
E_{1}(0) x = \frac{1}{2\pi} \int_{\Gamma} (\lambda - L(0))^{-1} x d\lambda,
\end{equation*}
where  $\Gamma$  is a closed continuous curve in left half complex plane containing $\sigma_{1}(L(0))$ in its interior. Due to the convergence of the spectrum,  it follows that  $\sigma_{1}(L(\epsilon))$ is also in the interior of  $\Gamma$ for  $\epsilon$ sufficiently small. Thus, we can write:

\begin{equation*}
E_{1}(\epsilon) x = \frac{1}{2\pi} \int_{\Gamma} (\lambda - L(\epsilon))^{-1} x d\lambda.
\end{equation*}

\begin{teo} \label{convprojesp}
 The spectral projection $E_{1}(\epsilon)$ is continuous at  $\epsilon = 0$ in $\mathcal{L}(X^{\beta}, X^{\eta})$ and,  $E_{2}(\epsilon)$ is continuous at $\epsilon=0$ in $\mathcal{L}(X^{\beta})$ and  in $\mathcal{L}(X^{\eta})$, $\frac{1}{2} > \eta > \frac{1}{4}$, $-\frac{1}{2} < \beta < -\frac{1}{4}$.
\end{teo}

\proof We have
\begin{align*}
 ||  E_{1}(\epsilon) - E_{1}(0) \parallel_{\mathcal{L}(X^{\beta}, X^{\eta})} &  \leq  \, \frac{1}{2\pi} \int_{\Gamma} \parallel ((\lambda - L(\epsilon))^{-1} - (\lambda - L(0))^{-1}) \parallel_{\mathcal{L}(X^{\beta}, X^{\eta})} d|\lambda| \\
& \leq  \frac{C({\epsilon})}{2\pi} \int_{\Gamma}  d|\lambda| \\
& \to 0 \,\,\, \text{as} \,\,\, \epsilon \to 0,
\end{align*}

\noindent where $C(\epsilon)$  comes from the convergence of the resolvents of the linearized operators obtained in Corollary  \ref{convreslinea_compact}.

\indent Since $E_{2}(\epsilon) = I - E_{1}(\epsilon)$, the result follows. \eproof

Suppose that the critical point   is hyperbolic and consider the operators  $L(\epsilon)_{1} = L(\epsilon)|_{X_{1}}$ and  $L(\epsilon)_{2} = L(\epsilon)|_{X_{2}}$.  Then  $L(\epsilon)_{1} : X_{1} \rightarrow X_{1}$ is bounded with  $\sigma (L(\epsilon)_{1}) = \sigma_{1}$  and  $L(\epsilon)_{2} : D(L(\epsilon)) \cap X_{2} \rightarrow X_{2}$ is sectorial (self-adjoint), with  $\sigma(L(\epsilon)_{2}) = \sigma_{2}$.
 We will need some estimates for the corresponding semigroups. Observe first that, by Lemma \ref{HAbounded}, the operator
  $((H_{\epsilon})_{\beta})_{u} (u, \epsilon)$ is  $(A_{\epsilon})_{\beta}$ bounded.  It follows, by easy computations, that
the fractional spaces defined by these operators are the same, with equivalent norms, uniformly in $\epsilon$.  We therefore, obtain the following estimates for the semigroups
\begin{eqnarray} \label{semig2}
\parallel e^{-L(\epsilon)_{2}t} \parallel_{\mathcal{L}(X^{\beta})} \leq M e^{-b t},  \  \parallel L(\epsilon) e^{-L(\epsilon)_{2}t} \parallel_{\mathcal{L}(X^{\beta})} \leq M t^{-1} e^{-bt } \, , \nonumber \\ \parallel e^{-L(\epsilon)_{2}t} \parallel_{\mathcal{L}(X^{\beta}, X^{\eta})} \leq N t^{-\eta + \beta} e^{-bt}, \, t \geq 0.
\end{eqnarray}

\begin{equation} \label{semig1}
\parallel e^{-L(\epsilon)_{1}t} \parallel_{\mathcal{L}(X^{\beta})} \leq M e^{ b t} \, , \, \parallel L(\epsilon) e^{-L(\epsilon)_{1}t} \parallel_{\mathcal{L}(X^{\beta})} \leq M e^{bt } \, , \, \parallel e^{-L(\epsilon)_{1}t} \parallel_{\mathcal{L}(X^{\beta}, X^{\eta})} \leq N e^{bt} \, , \, t \leq 0.
\end{equation}
where $b$, $M$, $N$ are constants independent of $\epsilon$.

We now show the convergence of the semigroups, when $\epsilon$ goes to $0$

\begin{lema} \label{difseml12}Under the conditions of Theorem \ref{convreslinea}, we have
 \begin{align*}
 \| e^{-L(\epsilon)_{2}t} E_{2}(\epsilon)  - e^{-L(0)_{2}t} E_{2}(0)  \parallel_{\mathcal{L} (X^{\beta}, X^{\eta})} \to 0 \text{ as }
  \epsilon \to 0, \\
  \| e^{-L(\epsilon)_{1}t} E_{1}(\epsilon)  - e^{-L(0)_{1}t} E_{1}(0)  \parallel_{\mathcal{L} (X^{\beta}, X^{\eta})} \to 0 \text{ as }
  \epsilon \to 0.
 \end{align*}

\end{lema}

\proof It is enough to prove the first inequality, the other is similar. We have
\begin{align*}
&    \left| \left| e^{-L(\epsilon)_{2}t} E_{2}(\epsilon)  - e^{-L(0)_{2}t} E_{2}(0)  \right| \right|_{\mathcal{L} (X^{\beta}, X^{\eta})} \\
& =  \, \left| \left| \frac{1}{2\pi i} \int_{\Gamma} ((\lambda + L(\epsilon))^{-1} E_{2}(\epsilon)  - (\lambda + L(0))^{-1} E_{2}(0)) e^{\lambda t} d\lambda \right| \right|_{\mathcal{L}(X^{\beta}, X^{\eta})} \nonumber \\
& \leq \frac{1}{2\pi} \int_{\Gamma} \parallel ((\lambda + L(\epsilon))^{-1} - (\lambda + L(0))^{-1}) \parallel_{\mathcal{L} (X^{\beta}, X^{\eta})} \parallel E_{2}(0) \parallel_{\mathcal{L} (X^{\beta}, X^{\beta})} |e^{\lambda t}| \, |d\lambda|  \nonumber \\
& +  \frac{1}{2\pi} \int_{\Gamma} \parallel (\lambda + L(\epsilon))^{-1} \parallel_{\mathcal{L} (X^{\beta}, X^{\eta})} \parallel E_{2}(\epsilon) - E_{2}(0) \parallel_{\mathcal{L} (X^{\beta}, X^{\beta})} |e^{\lambda t}| \, |d\lambda|  \nonumber \\
& \leq  \frac{1}{2\pi} \int_{\Gamma} \parallel (I + (\lambda - \lambda_{0}) (\lambda - L(\epsilon))^{-1}) \cdot \\
   &   \quad  ( (\lambda_{0} - L(\epsilon))^{-1} - (\lambda_{0} - L(0))^{-1}) (I + (\lambda - \lambda_{0}) (\lambda - L(0))^{-1}) \parallel_{\mathcal{L}(X^{\beta}, X^{\eta})}  \nonumber \cdot \\
  &    \quad  \parallel E_{2}(0) \parallel_{\mathcal{L} (X^{\beta}, X^{\beta})} |e^{\lambda t}| \, |d\lambda|  \nonumber \\
& +  \frac{1}{2\pi} \int_{\Gamma} \parallel (\lambda + L(\epsilon))^{-1} \parallel_{\mathcal{L} (X^{\beta}, X^{\eta})} \parallel E_{2}(\epsilon) - E_{2}(0) \parallel_{\mathcal{L} (X^{\beta}, X^{\beta})} |e^{\lambda t}| \, |d\lambda|  \nonumber \\
& \leq  \frac{1}{2\pi} \int_{\Gamma} \parallel (I + (\lambda - \lambda_{0}) (\lambda - L(\epsilon))^{-1}) \parallel_{\mathcal{L}(X^{\eta}, X^{\eta})} \, \cdot \\
 &   \quad \parallel ( (\lambda_{0} - L(\epsilon))^{-1} - (\lambda_{0} - L(0))^{-1}) \parallel_{\mathcal{L}(X^{\beta}, X^{\eta})}  \cdot  \nonumber \\
&   \quad  \parallel (I + (\lambda - \lambda_{0}) (\lambda - L(0))^{-1}) \parallel_{\mathcal{L}(X^{\beta}, X^{\beta})} \parallel E_{2}(0) \parallel_{\mathcal{L} (X^{\beta}, X^{\beta})} |e^{\lambda t}| \, |d\lambda|  \nonumber \\
& +  \frac{1}{2\pi} \int_{\Gamma} \parallel (\lambda + L(\epsilon))^{-1} \parallel_{\mathcal{L} (X^{\beta}, X^{\eta})} \parallel E_{2}(\epsilon) - E_{2}(0) \parallel_{\mathcal{L} (X^{\beta}, X^{\beta})} |e^{\lambda t}| \, |d\lambda|  \nonumber \\
& \leq  \frac{C(\epsilon)}{2\pi} e^{-bt} \int_{\Gamma} (1 + |\lambda - \lambda_{0}| D) (1 + |\lambda - \lambda_{0}| \bar{D}) \parallel E_{2}(0) \parallel_{\mathcal{L} (X^{\beta}, X^{\beta})} |e^{(\lambda +b) t}| \, |d\lambda|  \nonumber \\
& +  \frac{C'(\epsilon)}{2\pi} e^{-bt} \int_{\Gamma} D |e^{(\lambda + b) t}| \, |d\lambda|  \nonumber \\
& =  C_{2}(\epsilon) e^{-bt} \label{diferencadois}
\to  0 \,\,\, \text{ as } \,\,\, \epsilon \to 0, \nonumber
\end{align*}

\noindent where $\lambda_{0} \in  \Gamma$ is fixed,  $D$, $\bar{D}$ are constants independent of  $\epsilon$ and  $C(\epsilon)$, $C'(\epsilon)$ converge to $0$ as   $\epsilon$ goes to $0$ by Theorems    \ref{convreslinea} and   \ref{convprojesp}. \eproof

 The following result is proved in \cite{kato}, Theorem I-4.10.

\begin{lema} \label{opT} Suppose the conditions of Theorem   \ref{convprojesp} are met.  Then, for
 $0\leq \epsilon \leq \epsilon_0$, with $\epsilon_0$ sufficiently small and  $i = 1,2$, there are isomorphisms  ${T_{i}(\epsilon)}$
 in   $X^{\beta}$ such that   $T_{i}(\epsilon)X_{i}(0) = X_{i}(\epsilon), T_{i}(0) = I$ and the map $\epsilon \to T_{i}(\epsilon) \in \mathcal{L}(X^{\beta})$ is continuous at $\epsilon =0$.
\end{lema}

We are now in a position to prove the central result of this section.

\begin{teo}\label{contvar} Suppose satisfied the hypotheses of Theorem \ref{contequi}. Then,  there are constants  $\rho > 0$  and  $M \geq 1$ such that, for  $0 \leq \epsilon \leq \epsilon_0$,
 with $\epsilon_0$ sufficiently small,

\begin{itemize}
\item[i)] There exists a local invariant  manifold (the stable manifold) for the flow
generated by (\ref{problem_abstract_scale}) at  $e(\epsilon)$:
\begin{align*}
 &  W_{loc}^{s} (e(\epsilon)) \\
  &  =  \left\lbrace e(\epsilon) + z_{0} \in X^{\eta} \, : \, \parallel E_{2}(0) z_{0} \parallel_{X^{\eta}} \, \leq \frac{\varrho}{2M} , \, \parallel z(t, t_{0}, z_{0}, \epsilon) \parallel_{X^{\eta}} \, \leq \rho \, \text{for} \, t \geq t_{0} \right\rbrace,
\end{align*}
\noindent where  $z(t, t_{0}, z_{0}, \epsilon)$ is the solution of the linearized equation
\begin{equation} \label{eqlinear}
z_{t} + L(\epsilon)z = r(z, \epsilon)
\end{equation}
\noindent for  $t \geq t_{0}$  with initial value $z_{0}$. When  $z_{0} + e(\epsilon) \in W_{loc}^{s}(e(\epsilon))$, $\parallel z(t, t_{0}, z_{0}, \epsilon) \parallel_{X^{\eta}} \to 0$ as $t \to \infty$.
\item[ii)]  There exists a local invariant   manifold (the unstable manifold) for the flow
generated by (\ref{problem_abstract_scale}) at  $e(\epsilon)$:
\begin{align*}   &   W_{loc}^{u}  (e(\epsilon))   \\
 & =
    \left\lbrace   e(\epsilon) +  z_{0} \in X^{\eta} \, : \, \parallel E_{1}(0) z_{0} \parallel_{X^{\eta}} \, \leq \frac{\varrho}{2M} , \, \parallel z(t, t_{0}, z_{0}, \epsilon) \parallel_{X^{\eta}} \, \leq \rho \, \text{for } \, t \leq t_{0} \right\rbrace,
\end{align*}
\noindent where  $z(t, t_{0}, z_{0}, \epsilon)$ is the solution of the  linearized equation  (\ref{eqlinear}) for  $t \leq t_{0}$  with initial value  $z_{0}$. When $z_{0} + e(\epsilon) \in W_{loc}^{u}(e(\epsilon))$, $\parallel z(t, t_{0}, z_{0}, \epsilon) \parallel_{X^{\eta}} \to   0$ as  $t \to  -\infty$.
\item[iii)] If  $\beta(O,Q) = \sup_{o \, \in \, O} \inf_{q \, \in \, Q} \parallel q - o \parallel_{X^{\eta}}$ para $O, Q \subset X^{\eta}$, then
\begin{equation*}
\beta( W^{s}_{loc}(e(\epsilon)), W^{s}_{loc}(e(0)) ), \beta( W^{s}_{loc}(e(0)), W^{s}_{loc}(e(\epsilon)) )
\end{equation*}
\begin{equation*}
\beta( W^{u}_{loc}(e(\epsilon)), W^{u}_{loc}(e(0)) ), \beta( W^{u}_{loc}(e(0)), W^{u}_{loc}(e(\epsilon)) )
\end{equation*}
\noindent converge to zero  when  $\epsilon$ goes to  zero.
\end{itemize}
\end{teo}

\proof
  The first two claims follow from Theorem  5.2.1 in
 \cite{henry_semilinear},

 Now, given $\varrho > 0$, let

\begin{equation*}
Y_{0} = \left\lbrace a \in X_{2}(0) \, ; \, \parallel a \parallel_{X^{\eta}} \, \leq  \, \frac{\varrho}{2M} \right\rbrace
\end{equation*}

\noindent and

\begin{equation*}
Z_{0} = \{ z:[ t_{0}, \infty] \rightarrow X^{\eta} \, ; \, \text{z is continuous}, \, \sup \parallel z(t) \parallel_{X^{\eta}} \leq \varrho, \, E_{2}(0) z(t_{0}) = a, \, a \in Y_{0} \},
\end{equation*}
with the distance  $ \|z_1 - z_2\| := \sup \parallel z_1(t) -z_2 (t) \parallel_{X^{\eta}}.$

\indent For each $a \in Y_{0}$, consider the map $G_{a} : Z_{0} \times V \rightarrow Z_{0}$ defined by

\begin{eqnarray*}
G_{a}(z, \epsilon)(t) = e^{-L(\epsilon)_{2}(t - t_{0})} T_{2}(\epsilon)a + \int_{t_{0}}^{t} e^{-L(\epsilon)_{2} (t - s)} E_{2}(\epsilon) r(z(s), \epsilon) ds  \\ - \int_{t}^{\infty} e^{-L(\epsilon)_{1}(t - s)} E_{1}(\epsilon) r(z(s),\epsilon) ds.
\end{eqnarray*}

 We will show that
 $G_{a} ( \cdot, \epsilon) : Z_{0} \to Z_{0}$ is a contraction,  with contraction constant uniform in $\epsilon$, for  $a \in Y_{0}$.

  Using  (\ref{semig2}) and (\ref{semig1}) and   Theorem  \ref{propr}, item ii), we obtain

\begin{eqnarray*}
& & \parallel G_{a}(z_{1}, \epsilon)(t) - G_{a}(z_{2}, \epsilon)(t) \parallel_{X^{\eta}} \\
& \leq & \int_{t_{0}}^{t} \parallel e^{-L(\epsilon)_{2} (t - s)} E_{2}(\epsilon) (r(z_{1}(s), \epsilon) - r(z_{2}(s), \epsilon)) \parallel_{X^{\eta}} ds \\
& + & \int_{t}^{\infty} \parallel e^{-L(\epsilon)_{1}(t - s)} E_{1}(\epsilon) (r(z_{1}(s),\epsilon) - r(z_{2}(s),\epsilon)) \parallel_{X^{\eta}} ds \\
& \leq & \int_{t_{0}}^{t} N(t-s)^{-\eta + \beta} e^{-b(t - s)} \parallel E_{2}(\epsilon) \parallel_{\mathcal{L}(X^{\beta})} \parallel r(z_{1}(s), \epsilon) -  r(z_{2}(s), \epsilon) \parallel_{X^{\beta}} ds \\
& + & \int_{t}^{\infty} N e^{-b(t - s)} \parallel E_{1}(\epsilon) \parallel_{\mathcal{L}(X^{\beta})} \parallel r(z_{1}(s), \epsilon) -  r(z_{2}(s), \epsilon)  \parallel_{X^{\beta}} ds \\
& \leq & \int_{t_{0}}^{t} N(t-s)^{-\eta + \beta} e^{-b(t - s)} \parallel E_{2}(\epsilon) \parallel_{\mathcal{L}(X^{\beta})} k(\varrho) \parallel z_{1}(s) -  z_{2}(s) \parallel_{X^{\eta}} ds \\
& + & \int_{t}^{\infty} N e^{-b(t - s)} \parallel E_{1}(\epsilon) \parallel_{\mathcal{L}(X^{\beta})} k(\varrho) \parallel z_{1}(s) -  z_{2}(s)  \parallel_{X^{\eta}} ds.
\end{eqnarray*}

\noindent Thus

\begin{eqnarray*}
\sup_{t \, \in \, [t_{0}, \infty]} \parallel G_{a}(z_{1}, \epsilon)(t) - G_{a}(z_{2}, \epsilon)(t) \parallel_{X^{\eta}} \\
  \leq
  k(\varrho)  \lbrace \parallel E_{2}(0) \parallel_{\mathcal{L}(X^{\beta})} \int_{0}^{\infty} N(t - s)^{- \eta + \beta} e^{-b(t - s)} ds \\
  +  \parallel E_{1}(0)
  \parallel_{\mathcal{L}(X^{\beta})}
\int_{0}^{\infty} N e^{-b(t - s)} ds \rbrace
 \| z_1 -z_2  \|_{X^{\eta}},
\end{eqnarray*}
 where $k(\rho) $ is the constant given  in Theorem  \ref{propr}, item ii),  which goes to $0$ as $\rho \to 0$.
So we may  choose $\varrho > 0$  small enough  such that $G_a$ becomes a contraction for $\epsilon$ small. Thus, for any   $a \in Y_{0}$ and $\epsilon \leq \epsilon_0$, there exists a fixed point  $z(t,t_{0}, a, \epsilon)$, which is a solution of  $z_{t} + L(\epsilon)z = r(z, \epsilon)$, $t \geq t_{0}$, with initial value

\begin{equation*}
z(t_{0}, t_{0}, a, \epsilon) = T_{2}(\epsilon)a - \int_{t_{0}}^{\infty} e^{-L(\epsilon)_{1}(t_{0} - s)} E_{1}(\epsilon) r(z(s,t_{0}, a , \epsilon),\epsilon) ds.
\end{equation*}

Now, for $\epsilon $ small enough, we have

\begin{eqnarray*}
& || & G_{a}(z, 0) (t) - G_{a} (z, \epsilon)(t) \parallel_{X^{\eta}}\\
& \leq & \parallel e^{-L(\epsilon)_{2}(t - t_{0})} E_{2}(\epsilon) (T_{2}(\epsilon) - T_{2}(0))a \parallel_{X^{\eta}} \\
& + & \parallel  (e^{-L(\epsilon)_{2}(t - t_{0})} E_{2}(\epsilon) - e^{-L(0)_{2}(t - t_{0})} E_{2}(0)) T_{2}(0)a \parallel_{X^{\eta}} \\
& + & \left| \left| \int_{t_{0}}^{t} (e^{-L(\epsilon)_{2}(t - s)} E_{2}(\epsilon) - e^{-L(0)_{2}(t - s)} E_{2}(0)) r(z(s), \epsilon) ds \right| \right|_{X^{\eta}} \\
& + & \left| \left| \int_{t_{0}}^{t} e^{-L(0)_{2}(t - s)} E_{2}(0)(r(z(s), \epsilon) - r(z(s), 0)) ds \right| \right|_{X^{\eta}} \\
& + & \left| \left| \int_{t}^{\infty} (e^{-L(\epsilon)_{1}(t - s)} E_{1}(\epsilon) - e^{-L(0)_{1}(t - s)} E_{1}(0)) r(z(s), \epsilon) ds \right| \right|_{X^{\eta}} \\
& + & \left| \left| \int_{t}^{\infty} e^{-L(0)_{1}(t - s)} E_{1}(0)(r(z(s), \epsilon) - r(z(s), 0)) ds \right| \right|_{X^{\eta}} \\
& \leq & M e^{-b(t - t_{0})} \parallel E_{2}(\epsilon) \parallel_{\mathcal{L}(X^{\eta})} \, \parallel (T_{2}(\epsilon) - T_{2}(0)) \parallel_{\mathcal{L}(X^{\eta})} \parallel a \parallel_{X^{\eta}} + C_{2} (\epsilon)e^{-b( t - t_{0})} \parallel a \parallel_{X^{\eta}} \\
& + & \sup_{ \parallel z \parallel_{X^{\eta}} \, \leq \, \varrho} \parallel (r(z(s), \epsilon) - r(z(s), 0)) \parallel_{X^{\beta}} ( \parallel E_{2}(0) \parallel_{\mathcal{L}(X^{\beta})} \int_{0}^{\infty} N(t - s)^{-\eta + \beta}e^{-b(t-s)} ds \\
& + & \parallel E_{1}(0) \parallel_{\mathcal{L}(X^{\beta})} \int_{0}^{\infty} Ne^{-b(t-s)} ds ) \\
& + & C_{2}(\epsilon) \varrho k(\varrho) \left( \int_{0}^{\infty}   e^{-b(t -s)} ds \right) + C_{1}(\epsilon) \varrho k(\varrho) \left( \int_{0}^{\infty}   e^{-b(t -s)} ds \right),
\end{eqnarray*}

\noindent with $\parallel (T_{2}(\epsilon) - T_{2}(0)) \parallel_{\mathcal{L}(X^{\eta})}$, $\sup_{ \parallel z \parallel_{X^{\eta}} \, \leq \, \varrho} \parallel (r(z(s), \epsilon) - r(z(s), 0)) \parallel_{X^{\beta}}$, $C_{1}(\epsilon)$, $C_{2}(\epsilon)$  going to zero, when  $\epsilon \to 0$,  by Lemma \ref{propr}, Lemma \ref{difseml12}, Lemma \ref{opT} and the inequalities \eqref{semig2} and
 \eqref{semig1}. Therefore

\begin{equation*}
\sup_{t \, \in \,[0, \infty]} \parallel G_{a}(z, 0) (t) - G_{a} (z, \epsilon)(t) \parallel_{X^{\eta}} \leq C(\epsilon),
\end{equation*}

\noindent with  $C(\epsilon) \to 0$ as  $\epsilon \to 0$  and $G_{a}$ is continuous at  $\epsilon = 0$,
  uniformly for   $ a \in Y_{0}$.
  By the Fixed Point Theorem with parameters,  $z(t, t_{0}, a, \epsilon)$ is continuous at $\epsilon = 0 \in V$.

\indent For each  $\epsilon \leq \epsilon_0$, we have, by Theorem  5.2.1 in \cite{henry_semilinear}, that  $W^{s}_{loc} (e(\epsilon))$  is the image of the map $\phi_{\epsilon} : Y_{0} \to X^{\eta}$,
 given by

\begin{equation*}
\phi_{\epsilon} (a) = T_{2}(\epsilon)a - \int_{t_{0}}^{\infty} e^{-L(\epsilon)_{1}(t_{0} - s)} E_{1}(\epsilon) r(z(s, t_{0}, a, \epsilon),\epsilon) ds,
\end{equation*}

\noindent where  $z(s, t_{0}, a, \epsilon)$ is the solution of  $z_{t} + L(\epsilon)z = r(z, \epsilon)$, $t \geq t_{0}$, with initial value    $z(t_{0}, t_{0}, a, \epsilon) = \phi_{\epsilon} (a) = G_{a}(z, \epsilon)(t_{0})$. It follows that

\begin{equation*}
\parallel \phi_{0}(\cdot) - \phi_{\epsilon} (\cdot) \parallel_{X^{\eta}} \to 0 \,\,\, \text{as} \,\,\, \epsilon \to 0 \,\,\, \text{uniformly in } \,\,\, Y_{0}.
\end{equation*}

\noindent Since  $W^{s}_{loc} (e(\epsilon))$ is the image of the map  $\phi_{\epsilon}$,

\begin{equation*}
\beta( W^{s}_{loc}(e(\epsilon)), W^{s}_{loc}(e(0)) ) = \sup_{a \, \in \, Y_{0}} \inf_{b \, \in \, Y_{0}} \parallel \phi_{\epsilon} (a) - \phi_{0}(b) \parallel_{X^{\eta}},
\end{equation*}
%
it follows that

\begin{equation*}
\beta( W^{s}_{loc}(e(\epsilon)), W^{s}_{loc}(e(0)) ) \leq \sup_{a \, \in \, Y_{0}} \parallel \phi_{\epsilon} (a) - \phi_{0}(a) \parallel_{X^{\eta}} \to 0,
\end{equation*}

\noindent as   $ \epsilon \to  0$ in $V$. 

 We are  now in a position to prove the main result of this section

\begin{teo} \label{semicontinfatratores}
Suppose the hypotheses of  Theorem \ref{existatratorgeral} hold. Then, the family of attractors  $\{ \, \mathcal{A}_{\epsilon} \, | \, 0 \leq \epsilon \leq \epsilon_{0} \, \}$ of the flow $T_{\epsilon, \eta, \beta} (t,u)$, generated by (\ref{problem_abstract_scale})
 in  $X^{\eta}$, with  $\frac{1}{4} < \eta < \frac{1}{2}$ is lower semicontinuous  at $\epsilon = 0$.

\end{teo}

 \proof We just need to check the conditions of Theorem 4.10.8 in   \cite{hale} are met;

\begin{itemize}
\item[(H1)] $T_{0}(t)$ is a  $C^{1}$ gradient system, asymptotically compact, and the orbits of bounded sets are bounded. \\
This follows from Lemma  \ref{sistgrad}  and Theorem  \ref{existatratorh1}.

\item[(H2)]  The set of equilibria $E_{0}$ of  $T_{0}(t)$
is bounded in $X^{\beta} $.\\ This is proved in Lemma \ref{limequih1}.
\item[(H3)] Each element of  $E_{0}$ is hyperbolic.\\ This is part of our hypotheses.
\item[(H4)]   There exists a neighborhood  $U_{0}$ of  $\mathcal{A}_{0}$ which is attracted by
the the global attractor   $\mathcal{A}_{\epsilon}$ of  $T_{\epsilon}(t)$, for each  $\epsilon \neq 0$ .\\ Follows from  Theorems  \ref{existatratorgeral} and  \ref{bound_attractorh1}.
\item[(H5)] If  $E_{\epsilon}$ is the set of equilibria of  $T_{\epsilon}(t)$, then there exists a neighborhood  $W_{0}$ of  $E_{0}$ such that  $W_{0} \cap E_{\epsilon} = \{ \phi_{1, \epsilon}, ..., \phi_{N, \epsilon}\}$, where each  $\phi_{j, \epsilon}, \, j = 1, ..., N$, is hyperbolic and  $\phi_{j, \epsilon} \to \phi_{j}$ as  $\epsilon \to 0$.\\
 Follows from Theorems \ref{contsupeq} and \ref{contequi}.
\item[(H6)] $d(W^{u}_{loc}(\phi_{j}), W^{u}_{loc}(\phi_{j, \epsilon})) \to 0$ when  $\epsilon \to 0$.\\
Follows from Theorem \ref{contvar}.
\item[(H7)] For any  $\eta > 0, \, t_{0} > 0, \tau_{0} > 0$,
there exist  $\delta^{*}$ and $\epsilon_{0}^{*}$ such that
\begin{equation*}
\parallel T_{\epsilon}(t)y_{\epsilon} - T_{0}(t)x \parallel_{X^{\beta}} \, \leq \, \eta, \,\,\, \text{for} \,\,\, t_{0} \leq t \, \leq \, \tau_{0},
\end{equation*}
\noindent if  $x \in \mathcal{A}_{0}$, $y_{\epsilon} \in \mathcal{A}_{\epsilon}$, $\parallel y_{\epsilon} - x \parallel_{X^{\beta}} \, \leq \, \delta^{*}$. \\
Follows from Theorems  \ref{contsol}  and \ref{bound_attractorh1}.\end{itemize}

\eproof

\end{document}